\documentclass[onecolumn,final]{elsart3p}

\usepackage{amssymb}
\usepackage{amsfonts}
\usepackage{amsmath}
\usepackage{graphicx}
\usepackage{bm,bbm}
\usepackage{comment}
\usepackage{rotating}
\usepackage{color}
\usepackage{multirow}
\usepackage{multicol}
\usepackage{numprint}
\usepackage{overpic}

\usepackage[numbers,sort]{natbib}

\theoremstyle{plain}
\newtheorem{theorem}{Theorem}[section]

\newtheorem{lemma}{Lemma}[section]
\newtheorem{remark}{Remark}[section]

\theoremstyle{remark}

\newcommand{\fs}{f^s}
\newcommand{\Ps}{P^s}
\newcommand{\Ms}{M^s}

\renewcommand{\cv}{\mathbf{c}}

\newcommand{\Cv}   {\mathbf{C}}
\newcommand{\Cvs}  {\mathbf{C}^s}

\newcommand{\C}   {C}
\newcommand{\Cs}  {\C^s}

\newcommand{\Csn}  {\Cs_{n}}

\newcommand{\matA}{\mathbbm{A}}
\newcommand{\matB}{\mathbbm{B}}

\newcommand{\matD}{\mathbbm{D}}

\newcommand{\ABS}[1]{|#1|}

\newcommand{\ENERGY}{\mathcal{E}}

\newcounter{numbi}
\newcounter{numbii}

\newcommand{\vs}{v}
\newcommand{\Es}{E}

\newcommand{\ms}{m^s}
\newcommand{\qs}{q^s}

\newcommand{\Js}{J^s}

\newcommand{\dx}{dx}
\newcommand{\dv}{dv}
\newcommand{\dt}{dt}
\newcommand{\ds}{ds}

\newcommand{\Ensk}{\mathcal{E}^s_{\footnotesize{\mbox{\it kin}}}}

\newcommand{\JJ}{J}

\newcommand{\abs}[1]{\left|#1\right|}

\newcommand{\Diffv}[1]{\delta_v\big[ #1 \big]}

\newcommand{\norm}[1]{|\!|#1|\!|}

\def\trait #1 #2 #3 {\vrule width #1pt height #2pt depth #3pt}
\def\fin{\hfill
        \trait .3 5 0
        \trait 5 .3 0
        \kern-5pt
        \trait 5 5 -4.7
        \trait 0.3 5 0
\medskip}
\newcommand{\ENDPROOF}{\fin}
\newcommand{\BEGINPROOF}{\emph{Proof}.~~}

\newcommand{\prty}[2]{\pi_{#1,#2}}
\newcommand{\va}{\vs_a}
\newcommand{\vb}{\vs_b}

\newcommand{\NL}{N_L}
\newcommand{\NF}{N_F}
\newcommand{\NH}{N_H}

\newcommand{\Mk}{m_k}
\newcommand{\Epot}{\mathcal{E}_{\footnotesize{\mbox{\it pot}}}}

\newcommand{\Csnk}{\Cs_{n,k}}
\newcommand{\Cvk} {\mathbf{C}^{s}_k}
\newcommand{\Jsk} {\Js_{k}}
\newcommand{\img} {i}
\newcommand{\STAR}{\star}

\newcommand{\FLARGE}[1]{\begin{large}#1\end{large}}

\newcommand{\EOD}{\end{document}}

\usepackage[normalem]{ulem}

\begin{document}

\begin{frontmatter}
  \title{A Legendre-Fourier spectral method with exact 
    conservation laws for the Vlasov-Poisson
    system}
  
  
  \author[lanl,IMATI] {G. Manzini}        %
  \author[lanl] {, G.~L. Delzanno}  %
  \author[lanl] {, J. Vencels}      %
  \author[eth]  {, and S. Markidis}
  
  
  \address[lanl]{T-5 Applied Mathematics and Plasma Physics Group, Los Alamos National Laboratory,
    Los Alamos, NM 87545, USA}

  \address[eth] {HPCViz Department, KTH Royal Institute of Technology, Stockholm, Sweden}

  \address[IMATI]{ Istituto di Matematica Applicata e Tecnologie %
    Informatiche, Consiglio Nazionale delle Ricerche (IMATI-CNR), \\
    via Ferrata 1, I -- 27100 Pavia, Italy, %
  }

  \begin{abstract}
    We present the design and implementation of an $L^2$-stable
    spectral method for the discretization of the Vlasov-Poisson model
    of a collisionless plasma in one space and velocity dimension.
    The velocity and space dependence of the Vlasov equation are
    resolved through a truncated spectral expansion based on Legendre
    and Fourier basis functions, respectively.
    The Poisson equation, which is coupled to the Vlasov equation,
    is also resolved through a Fourier
    expansion.
    The resulting system of ordinary differential equation is
    discretized by the implicit second-order accurate Crank-Nicolson
    time discretization.
    The non-linear dependence between the Vlasov and Poisson equations
    is iteratively solved at any time cycle by a Jacobian-Free
    Newton-Krylov method.
    In this work we analyze the structure of the main conservation
    laws of the resulting Legendre-Fourier model, e.g., mass,
    momentum, and energy, and prove that they are exactly satisfied in
    the semi-discrete and discrete setting.
    The $L^2$-stability of the method is ensured by discretizing the
    boundary conditions of the distribution function at the boundaries
    of the velocity domain by a suitable penalty term.
    The impact
      of the penalty term on the conservation properties is
      investigated theoretically and numerically. An implementation of
      the penalty term that does not affect the conservation of mass,
      momentum and energy, is also proposed and studied.
    A collisional term is introduced in the discrete model to
      control the filamentation effect, but does not affect the
      conservation properties of the system.
    Numerical results on a set of standard test problems
    illustrate the performance of the method.
  \end{abstract}
  

\end{frontmatter}

\renewcommand{\arraystretch}{1.}

\raggedbottom




\newcommand{\BIB}{[\textbf{add citation}]{~}}

\section{Introduction}

%

Collisionless magnetized plasmas are described by the kinetic
(Vlasov-Maxwell) equations and are characterized by high
dimensionality, anisotropy and a wide variety of spatial and temporal
scales \cite{goldston}, thus requiring the use of sophisticated
numerical techniques to capture accurately their rich non-linear
behavior.

In general terms, there are three broad classes of methods devoted to
the numerical solution of the kinetic equations.
dimensional phase space via macro-particles that evolve according to
Newton's equations in the self-consistent electromagnetic field
\cite{birdsall,hockney}.
PIC is the most widely used method in the plasma physics community
because of its robustness and relative simplicity.
The well known statistical noise associated with the macro-particles
implies that PIC is really effective for problems where a low
signal-to-noise ratio is acceptable.
The Eulerian-Vlasov methods discretize the phase space with a six
dimensional computational mesh \cite{cheng&knorr,sonnendruker,filbet}.
As such they are immune to statistical noise but they require
significant computational resources and this is perhaps why their
application has been mostly limited to problems with reduced
dimensionality.
For reference, storing a field in double precision on a mesh with
$10^{12}$ cells requires about $8$ terabytes of memory.
A third class of methods, called transform methods, is spectral and is
based on an expansion of the velocity part of the distribution
function in basis functions (typically Fourier or Hermite), leading to
a truncated set of moment equations for the expansion coefficients
\cite{armstrong,engelmann,klimas,holloway96,schumer}.
Similarly to Eulerian-Vlasov methods, transform methods might be
resource-intensive if the convergence of the expansion series is slow.

In recent years there seems to be a renewed interest in Hermite-based
spectral methods.
Some reasons for this can be attributed to the advances in high
performance computing and to the importance of simpler, reduced
kinetic models in elucidating aspects of the complex dynamics of
magnetized plasmas \cite{parker,loureiro}.
Another reason is that (some form of) the Hermite basis can unify
fluid (macroscopic) and kinetic (microscopic) behavior into one
framework \cite{camporeale0,vencels,delzanno,Vencels:2016}.
Thus, it naturally enables the 'fluid/kinetic coupling' that might be
the (inevitable) solution to the multiscale problem of computational
plasma physics and is a very active area of research
\cite{markidis,daldorff}.

The Hermite basis is defined by the Hermite polynomials with a
Maxwellian weight and is therefore closely linked to Maxwellian
distribution functions.
Two kinds of basis have been proposed in the literature (differing in
regard to the details of the Maxwellian weight): symmetrically- and
asymmetrically-weighted \cite{holloway96,schumer}.
The former features $L^2$-stability but conservation laws for total
mass, momentum and energy are achieved only in limited cases (i.e.,
they depend on the parity of the total number of Hermite modes, on the
presence of a velocity shift in the Hermite basis, ...).
The latter features exact conservation laws in the discrete and the
connection between the low-order moments and typical fluid moments,
but $L^2$-stability is not guaranteed \cite{camporeale1,schumer}.
Earlier works pointed out that a proper choice of the velocity shift
and the scaling of the Maxwellian weight (free parameters of the
method) is important to improve the convergence properties of the
series \cite{holloway96,camporeale0}.
Indeed, the optimization of the Hermite basis is a crucial aspect of
the method, which however at this point does not yet have a definitive
solution.

One could of course envision a different spectral approach which
considers a full polynomial expansion without any weight or free
parameter.
While any connection with Maxwellians is lost, such expansion could be
of interest in presence of strong non-Maxwellian behavior and
eliminates the optimization problem.
The Legendre polynomials are a natural candidate in this case, because
of their orthogonality properties.
They are normally applied in some preferred coordinate system (for
instance spherical geometry) to expose quantities like angles that are
defined on a bounded domain.
Indeed, Legendre expansions are very popular in neutron transport
\cite{sanchez} and some application in kinetic plasma physics can be
found for electron transport described by the Boltzmann equation
\cite{sosov}.
Surprisingly, however, we have not found any example in the context of
collisionless kinetic theory and in particular for the Vlasov-Poisson
system.

The main contribution of the present paper is the formulation,
development and successful testing of a spectral method for the one
dimensional Vlasov-Poisson model of a plasma based on a Legendre
polynomial expansion of the velocity part of the plasma distribution
function.
The expansion is applied directly in the velocity domain, which is
assumed to be finite.
It is shown that the Legendre expansion features many of the
properties of the asymmetrically-weighted Hermite expansion: the
structure of the equations is similar, the low-order moments
correspond to the typical moments of a fluid, and conservation laws
for the total mass, momentum and energy (in weak form, as defined in
Sec. 4) can be proven.
It also features properties of the symmetrically-weighted Hermite
expansion: $L^2$-stability is also achieved by introducing a penalty
on the boundary conditions in weak form.
This strategy is inspired by the Simultaneous Approximation Strategy
(SAT) technique
~\cite{Svard-Nordstrom:2014,Kreiss-Scherer:1974,Kreiss-Scherer:1977,Strand:1994,Carpenter-Gottlieb-Abarbanel:1994,Nordstrom-Forsberg-Adamsson-Eliasson:2003,Carpenter-Nordstrom-Gottlieb:1999,Mattsson-Nordstrom:2004}.

The paper is organized as follows. 
In Sec.~\ref{sec:vlasov} the Vlasov-Poisson equations for a plasma are
introduced together with the spectral discretization: the velocity
part of the distribution function is expanded in Legendre polynomials
while the spatial part is expressed in terms of a Fourier series.
The time discretization is handled via a second-order accurate
Crank-Nicolson scheme.
In Sec.~\ref{sec:L2:stability} the SAT technique is used to enforce
the $L^2$-stability of the numerical scheme.
In Sec.~\ref{sec:conservation:laws} conservation laws for the total
mass, momentum and energy are derived theoretically.
Numerical experiments on standard benchmark tests (i.e., Landau
damping, two-stream instabilities and ion acoustic wave) are performed
in Sec.~\ref{sec:numerical}, proving numerically the stability of the
method and the validity of the conservation laws. Conclusions are
drawn in Sec.~\ref{sec:conclusions}.

\section{The Vlasov-Poisson system and the Legendre-Fourier
  approximation}
\label{sec:vlasov}

We consider the Vlasov-Poisson model for a collisionless plasma of
electrons (labeled ``$e$'') and singly charged ions (``$i$'') evolving
under the action of the self-consistent electric field $E$.
The behavior of each particle species $s\in\{e,i\}$ with mass $\ms$
and charge $\qs$ is described at any time $t\geq 0$ in the phase space
domain $[0,L]\times[\va,\vb]$ by the distribution function
$\fs(x,\vs,t)$, which evolves according to the \emph{Vlasov equation}:
\begin{align}
  \frac{\partial\fs}{\partial t} + v\frac{\partial\fs}{\partial x} +
  \frac{\qs}{\ms}E\,\frac{\partial\fs}{\partial v} &= 0.
  \label{eq:Vlasov}
\end{align}
We assume the physical space to be periodic in $x$, so that no
boundary condition for $\fs$ is necessary at $x=0$ and $x=L$, and that
suitable boundary conditions, e.g., $\fs(x,\va,t)=\fs(x,\vb,t)=0$, are
provided for $\fs$ at the velocity boundaries $\vs=\va$ and $\vs=\vb$
for any time $t\geq 0$ and any spatial position $x\in[0,L]$.
We also assume that an initial solution $\fs(x,\vs,0)=\fs_0(x,\vs)$ is
given at the initial time $t=0$.

\medskip
\begin{remark}
  \label{remark:zero:velocity:bcs}
  If the initial solution $\fs_0(x,\vs)$ has a compact support in the
  phase space domain $[0,L]\times[\va,\vb]$, then $\fs$ has also a
  compact support at any time $t>0$.
  Moreover, the size of the support may increase in time in a
  controlled way, cf.~\cite{Wollman:1982,Glassey-1996}.
  In such a case, it holds that $\fs(x,\va,t)=\fs(x,\vb,t)=0$ until
  the size of the support equals the size of the velocity domain.
  This condition can be used to determine the final time at which a
  plasma simulation based on this numerical model is valid.
\end{remark}

\medskip 
In the Vlasov-Poisson system, the electric field $E(x,t)$ is the
solution of the \emph{Poisson equation}:
\begin{align}
  \epsilon_0\frac{\partial\Es}{\partial x} = \rho
  \qquad
  \text{for~}x\in(0,L),\,t\geq 0,
  \label{eq:Poisson}
\end{align}
where $\epsilon_0$ is the dielectric constant and
\begin{align}
  \rho(x,t) = \sum_{s\in\{e,i\}} q_s\int_{-\infty}^{\infty}\fs(x,v,t)\,dv
  \label{eq:rho:def}
\end{align}
is the \emph{total charge density} of the plasma.
By taking the time derivative of the Poisson equation and using the
continuity equation
\begin{align*}
  \frac{\partial\rho}{\partial t} + \frac{\partial\JJ}{\partial x} = 0,
\end{align*}
where $\JJ(x,t)$ is the \emph{total current density} defined as
\begin{align}
  \JJ(x,t) = \sum_{s\in\{e,i\}} \qs \int_{-\infty}^{\infty}\vs\fs(x,v,t)\,\dv,
  \label{eq:Js:def}
\end{align}
we obtain the \emph{Ampere equation}
\begin{align}
  \epsilon_0\frac{\partial\Es}{\partial t} + \JJ = C_{A},
  \label{eq:Ampere:equation:def}
\end{align}
where $C_{A}$ is a suitable constant factor.

\medskip
\begin{remark}
  The Ampere equation can be used with the Vlasov equation instead of
  the Poisson equation to obtain the Vlasov-Ampere formulation.
  In the continuum setting, the two formulations are equivalent
  in the one-dimensional electrostatic case without any external
  electric field as the one considered in this work.
\end{remark}

\subsection{Velocity integration using Legendre expansion}
Consider the infinite set of \emph{Legendre polynomials}
$\{L_n(\eta)\}_{n=0}^{\infty}$, which are recursively defined for
$\eta\in[-1,1]$ by~\cite[Chapters 8, 22]{Abramowitz}:
\begin{align}
  L_0(\eta)=1, \quad 
  L_1(\eta)=\eta, \quad
  (2n+1) \eta L_{n}(\eta)=(n+1)L_{n+1}(\eta)+nL_{n-1}(\eta) \textrm{~~for~} n\geq 1,
  \label{eq:legendre:def}
\end{align}
and normalized as follows
\begin{align}
  \int_{-1}^{1}L_m(\eta)L_n(\eta)\ds = \frac{2}{2n+1}\delta_{mn}.
 \label{eq:orthogonality}
\end{align}
We remap the Legendre polynomials onto the velocity range $[\va,\vb]$
through the linear transformation
$\vs(\eta)=(\va+\vb)\slash{2}+\eta(\vb-\va)\slash{2}$.
Let $\eta(\vs)=( 2\vs -(\va+\vb) )\slash{(\vb-\va)}$ be the inverse
mapping from $[\va,\vb]$ to $[-1,1]$.
The $n$-th Legendre polynomial is given by
$\phi_n(\vs)=\sqrt{2n+1}L_{n}( \eta(\vs) )$, where the scaling factor in
front of $L_{n}$ is chosen to satisfy the orthogonality relation:
\begin{align}
  \int_{\va}^{\vb}\phi_m(\vs)\phi_n(\vs)\dv = (\vb-\va) \delta_{mn}.
  \label{eq:Legendre:orthogonality}
\end{align}
The first derivative of the Legendre polynomials is given by
\begin{align*}
  \frac{d L_n(\eta)}{d\eta} =\sum_{i=0}^{n-1} \prty{n}{i} (2i+1)L_{i}(\eta)
  \qquad\eta\in[-1,1],
\end{align*}
where $\prty{n}{i}$ is a switch that takes value $0$ if $n-i$ is even,
and $1$ if $n-i$ is odd.
Using the chain rule and adjusting the normalization factor, we obtain
the first derivative of the translated and rescaled Legendre
polynomials:
\begin{align}
  \frac{d\phi_n}{\dv} 
  = \sqrt{2n+1}\frac{d}{\dv}L_{n}(\eta(\vs))
  = \sum_{i=0}^{n-1} \sigma_{n,i} \phi_{i}(\vs)
  \quad\textrm{~with~}\quad
  \sigma_{n,i} = \prty{n}{i}\,\frac{2\sqrt{ \big(2n+1\big)\,\big(2i+1\big) }}{\vb-\va}.
  \label{eq:legendre:first:derivative}
\end{align}
The recursion relations that are used to expand the Vlasov equation on
the Legendre basis are reported in appendix A for completeness.
 
\medskip 
Consider the spectral decomposition of the distribution function $\fs$
on the basis of Legendre polynomials given by
\begin{align}
  \fs(x,\vs,t) = \sum_{n=0}^{\infty}\Csn(x,t)\phi_n(\vs)
  \quad\textrm{~where~}\quad
  \Csn(x,t) = \frac{1}{\vb-\va} \int_{\va}^{\vb}\fs(x,\vs,t)\phi_n(\vs)\dv
  \label{eq:legendre:decomposition}
\end{align}
in the Vlasov equation~\eqref{eq:Vlasov}. 
The boundary conditions $\fs(x,\va,t)=\fs(x,\vb,t)=0$ are not exactly
satisfied since they are imposed in weak form and the polynomials
$\phi_n(\vs)$ are not zero at the velocity boundaries.

A possible way to circumvent this issue is to consider the modified
basis functions given by $\tilde{\phi}_{n}=\phi_{n+2}-\phi_{n}$ for
$n\geq 0$.
From the properties of the Legendre polynomials, it readily follows
that $\tilde{\phi}_{n}(\va)=\tilde{\phi}_{n}(\vb)=0$ for each $n$ and
the expansion of $\fs$ on this set of functions will automatically
satisfied the homogeneous conditions at the boundary of the velocity
range.
Nonetheless, we verified numerically in the first stages of this work
that this approach may yield an unstable method and the numerical
instability cannot be fixed as there is no mechanism that allows us to
control the growth of the absolute value of the Legendre coefficients
$\Csn$.
Another possible choice is to consider
$\tilde{\phi}_n(v)=\phi_{n+2}(v)-1$
for even $n\geq0$ and $\tilde{\phi}_n(v)=\phi_{n+2}(v)-v$ for odd 
$n\geq1$.
Although we have not implemented this second basis, 
a common characteristic of these choices is the loss of orthogonality,
which we suspect may influence negatively the stability properties of the
method.
The alternative approach that we consider hereafter is to integrate by
parts the velocity term in the Vlasov equations.
This strategy allows us to set the boundary conditions in weak form,
and, then, to introduce a penalty term to enforce the $L^2$ stability
of the method through the boundary conditions (see
Section~\ref{sec:L2:stability}).
To this end, we substitute~\eqref{eq:legendre:decomposition}
into~\eqref{eq:Vlasov}, we multiply the resulting equation by
$\phi_{n}$ and integrate over $[\va,\vb]$.
Then, we use the recursion
formulas~\eqref{eq:intg:f}-\eqref{eq:intg:v2f} and the orthogonality
relation~\eqref{eq:orthogonality} and we obtain the following system
of partial differential equations for the Legendre coefficients
$\Csn(x,t)$:
\begin{align}
  \frac{\partial\Cs_{n}}{\partial t} 
  +
  \sigma_{n+1}\frac{\partial\Cs_{n+1}}{\partial x} +
  \sigma_{n}  \frac{\partial\Cs_{n-1}}{\partial x} +
  \overline{\sigma}\frac{\partial\Cs_{n}}{\partial x}
  - 
  \frac{\qs}{\ms}E\left( 
    \sum_{i=0}^{n-1}\sigma_{n,i}\Cs_{i} 
    - \Diffv{\fs\phi_n}_{\va}^{\vb}
  \right)
   = 0
   \quad\textrm{for~}n\geq 0,
  \label{eq:legendre:system}
\end{align}
where conventionally $\Cs_{-1}=0$,
\begin{align}
  \overline{\sigma}= \frac{\va+\vb}{2},\qquad
  \sigma_{n} = 
  \begin{cases}
  0 & \textrm{for~}n=0,\\[0.5em]
  \displaystyle\frac{\vb-\va}{2}\,\frac{n}{\sqrt{(2n+1)(2n-1)}}
  & \textrm{for~}n\geq 1,
  \end{cases}
  \label{eq:legendre:sigma:def}
\end{align}
and
\begin{align}
  \Diffv{\fs\phi_n}_{\va}^{\vb} = 
  \frac{\fs(x,\vb,t)\phi_n(\vb)-\fs(x,\va,t)\phi_n(\va)}{\vb-\va},
  \label{eq:delta_v:def}
\end{align}
is the boundary term resulting from an integration by parts of the
integral term that involves the velocity derivative.
The derivation of the coefficients $\sigma_{n}$ and
$\overline{\sigma}$ can be found in appendix A.
If the distribution $\fs(x,\vs,t)$ has compact support in $\vs$, the
homogeneous boundary conditions at $\vs=\va$ and $\vs=\vb$ are imposed
in weak form by assuming that $\Diffv{\fs\phi_n}_{\va}^{\vb}$
in~\eqref{eq:delta_v:def} is zero.
However, since this term plays a major role in establishing the
conservation laws and ensuring the $L^2$ stability of the
discretization method, we will consider it in all the further
developments and in the analysis of the next sections.

We truncate the spectral expansion of $\fs$ after the first $\NL$
\emph{Legendre modes} by assuming that $\Cs_{n}=0$ for $n\geq\NL$ and
we approximate the distribution function by the finite summation:
\begin{align}
  \fs(x,v,t)\approx
  \fs_{L}(x,v,t) = \sum_{n=0}^{\NL-1}\Cs_{n}(x,t)\phi_n(\vs).
\end{align}
The evolution of each coefficient $\Csn$ with $n\leq\NL-1$ is still
given by~\eqref{eq:legendre:system}.
To ease the notation, we will drop the subindex $L$ in $\fs_{L}$ by
tacitly assuming that all the quantities containing $\fs$ are indeed
numerical approximations dependent on the first $\NL$ modes of the
truncated series.

\medskip 
Let $\Cv$ be the vector that contains all the coefficients $\Csn$ for
$n\in[0,\NL-1]$, i.e., $(\Cvs)_{n}=\Cs_{n}$, and
${\bm\phi(v)=(\phi_0(v),\phi_1(v),\ldots,\phi_{\NL-1}(v))^T}$ the
vector containing the values of the Legendre shape functions evaluated
at $v$.
It holds that $\fs(x,\vs,t)={\bm\phi(v)}^T\Cvs(x,t)$.
System~\eqref{eq:legendre:system} can be rewritten in the
non-conservative vector form:
\begin{align}
  \frac{\partial\Cvs}{\partial t} 
  +\matA\frac{\partial\Cvs}{\partial x}
  -\frac{\qs}{\ms}E\left(
  \matB\Cvs
  -\Diffv{\fs{\bm\phi}}_{\va}^{\vb}
  \right)
   = 0,
  \label{eq:vlasov:legendre:non-conservation}
\end{align}
where
\begin{align}
  \left(\matA\frac{\partial\Cvs}{\partial x}\right)_{n} &=
  \sigma_{n+1}\frac{\partial\Cs_{n+1}}{\partial x} +
  \sigma_{n}  \frac{\partial\Cs_{n-1}}{\partial x} +
  \overline{\sigma}\frac{\partial\Cs_{n}}  {\partial x},
  \label{eq:legendre:matA:def}\\[0.5em]
  (\matB\Cvs)_{n} &= \sum_{i=0}^{n-1}\sigma_{n,i}\Cs_{i},
  \label{eq:legendre:matB:def}
\end{align}
and $\big(\Diffv{\fs{\bm\phi}}_{\va}^{\vb}\big)_{n}=\Diffv{\fs\phi_{n}}_{\va}^{\vb}$.
Since $\matA$ is a constant matrix, it follows that
\begin{align}
  \matA\frac{\partial\Cvs}{\partial x} = \frac{\partial}{\partial x}(\matA\Cvs)
  \quad\textrm{with}\quad
  \left(\matA\Cvs\right)_{n} =  
  \sigma_{n+1}\Cs_{n+1} +
  \sigma_{n}  \Cs_{n-1} +
  \overline{\sigma}\Cs_{n}.
\end{align}
Therefore, system~\eqref{eq:legendre:system} also admits the
conservative form:
\begin{align}
  \frac{\partial\Cvs}{\partial t} + \frac{\partial}{\partial x}\big(\matA\Cvs\big) 
  -\frac{\qs}{\ms}\,E\,\left( \matB\Cvs
  -\Diffv{\fs{\bm\phi}}_{\va}^{\vb}
  \right) = 0,
  \label{eq:vlasov:legendre:conservation}
\end{align}
where $\matA$ is a real and symmetric matrix with real
eigenvalues and eigenvectors.

\subsection{Space integration using Fourier expansion}
We expand each Legendre coefficient $\Csn(x,t)$ on the first $2\NF+1$
functions of the Fourier basis $\psi_k(x)=e^{\frac{2\pi\img}{L} kx}$
(for $k\in[-\NF,\NF]$) as follows
\begin{align}
  \Csn(x,t) = \sum_{k=-\NF}^{\NF}\Cs_{n,k}(t)\psi_k(x),
  \label{eq:fourier:decomposition}
\end{align}
where each coefficient $\Cs_{n,k}(t)$ is a complex function of time
$t$.
The Fourier basis functions satisfy the orthogonality relation
\begin{align}
  \int_{0}^{L}\psi_k(x)\psi_{k'}(x)\dx = L\delta_{k+k',0}.
  \label{eq:Fourier:orthogonality}
\end{align}
Substituting~\eqref{eq:fourier:decomposition}
in~\eqref{eq:legendre:system} and
using~\eqref{eq:Fourier:orthogonality}, we derive the system for the
coefficients $\Cs_{n,k}$, which reads as:
\begin{align}
  \frac{d\Cs_{n,k}}{\dt} 
  +\left(\frac{2\pi\img}{L}k\right)\,
  \Big(
  \sigma_{n+1}\Cs_{n+1,k} +
  \sigma_{n}  \Cs_{n-1,k} +
  \overline{\sigma}\Cs_{n,k}
  \Big)
  -\frac{\qs}{\ms}
  \left[
    \Es\STAR
    \left(
      \sum_{i=0}^{n-1}\sigma_{n,i}\Cs_{i}
      -\Diffv{\fs\phi_n}_{\va}^{\vb}
    \right)
  \right]_{k} = 0,
  \label{eq:legendre:fourier:system}
\end{align}
for $n=0,\ldots,\NL-1$ and $-\NF\leq k\leq\NF$, and where $\STAR$
denotes the convolution integral and $\big[\,\cdot\,\big]_{k}$ denotes
the $k-{th}$ mode of the Fourier expansion of the argument inside the
square brackets.
Explicit formulas for these quantities are given below.
We also recall that if $g(x)$ and $h(x)$ are two given real functions
of $x\in[0,L]$ and $\{g_k\}$ and $\{h_k\}$ the coefficients of their
Fourier expansion on the basis functions $\psi_k$, then the $k$-th
Fourier mode of the convolution product $g\STAR h$ is given by
  $[g\STAR h]_k = \sum_{k'=-\NF}^{\NF}g_{k'} h_{k-k'}$.

\medskip
The Poisson equation for the electric field is similarly transformed
by using~\eqref{eq:fourier:decomposition} and the Fourier expansion of
the electric field
\begin{align}
  E(x,t) = \sum_{k=-\NF}^{\NF}\Es_{k}(t)\psi_k(x)
  \label{eq:E:Fourier}
\end{align}
into~\eqref{eq:Poisson} to obtain
\begin{equation}
  \epsilon_{0}\left(\frac{2\pi\img}{L}k\right)\,E_k(t) 
  = (\vb-\va)\,\sum_{s\in\{e,i\}}\qs\Cs_{0,k}(t).
  \label{eq:legendre:fourier:Poisson:1}
\end{equation}
For $k=0$ the equation above becomes
\begin{align*}
  \sum_{s\in\{e,i\}}\qs\Cs_{0,0}(t) = 0,
\end{align*}
which, according to the hypothesis of neutrality of the plasma,
expresses the fact that the total charge in the system is zero.
This implies that we can set the $0$-th Fourier mode of the electric
field to zero, i.e.,
\begin{align*}
  L\,E_0(t) 
  = \int_{0}^{L}\Es(x,t)\dx = 0.
\end{align*}

\medskip
For convenience of notation, we introduce the vector $\Cvk$ that
contains the $k$-th Fourier coefficients $\Csnk$ for all the Legendre
modes $n\in[0,\NL-1]$, i.e., $(\Cvk)_{n}=\Cs_{n,k}$.
System~\eqref{eq:legendre:fourier:system} can be rewritten in the
vector form:
\begin{align}
  \frac{d\Cvk}{\dt} 
  +\left(\frac{2\pi\img}{L}k\right)\matA\Cvk
  -\frac{\qs}{\ms}\left[E\STAR\left(
      \matB\Cvs
      -\Diffv{\fs{\bm\phi}}_{\va}^{\vb}
    \right)\right]_k
   = 0,
   \label{eq:legendre:fourier:compact}
\end{align}
where
\begin{align}
  \left(\matA\Cvk\right)_{n} &=
  \sigma_{n+1}\Cs_{n+1,k} +
  \sigma_{n}  \Cs_{n-1,k} +
  \overline{\sigma}\Cs_{n,k}
  \label{eq:legendre:fourier:matA:def}\\[0.5em]
  (\matB\Cvk)_{n} &= \sum_{i=0}^{n-1}\sigma_{n,i}\Cs_{i,k}
  \label{eq:legendre:fourier:matB:def}\\[0.5em]
  \big(\Diffv{\fs{\bm\phi}}_{\va}^{\vb}\big)_{n}
  &
  =\Diffv{\fs\phi_{n}}_{\va}^{\vb}
  =\frac{1}{\vb-\va}\Big(\fs(x,\vb,t)\phi_{n}(\vb) - \fs(x,\va,t)\phi_{n}(\va)\Big).
\end{align}
Note the vector expressions:
\begin{align*}
  \Big[\Es\STAR\matB\Cv\Big]_{k} &= \sum_{k'=-\NF}^{\NF}\Es_{k'}\big[\matB\Cv\big]_{k-k'}\\[0.5em]
  \Big[\Es\STAR\Diffv{\fs{\bm\phi}}_{\va}^{\vb}\Big]_{k} &= 
  \sum_{k'=-\NF}^{\NF}\Es_{k'}
  \sum_{n'=0}^{\NL-1}\Cs_{n',k-k'}(t)\big(\phi_{n'}(\vb){\bm\phi(\vb)}-\phi_{n'}(\va){\bm\phi(\va)}\big)
\end{align*}
and for the $n$-th Legendre components:
\begin{align*}
  \Big[\Es\STAR\matB\Cv\Big]_{n,k} &= \sum_{k'=-\NF}^{\NF}\Es_{k'}\sum_{i=0}^{n-1}\sigma_{n,i}\Cs_{i,k-k'}(t)\\[0.5em]
  \Big[\Es\STAR\Diffv{\fs{\bm\phi}}_{\va}^{\vb}\Big]_{n,k} &= 
  \sum_{k'=-\NF}^{\NF}\Es_{k'}
  \sum_{n'=0}^{\NL-1}\Cs_{n',k-k'}(t)\big(\phi_{n'}(\vb)\phi_n(\vb)-\phi_{n'}(\va)\phi_n(\va)\big).
\end{align*}

Consider the current density of species $s$ given
by~\eqref{eq:Js:def}.
We apply the Legendre decomposition~\eqref{eq:legendre:decomposition},
the Fourier decomposition~\eqref{eq:fourier:decomposition} and we
use~\eqref{eq:intg:vf} to obtain the Legendre-Fourier representation
of the total current density:
\begin{align}
  \Jsk(t) = \qs(\vb-\va)L\,\big( \sigma_{1}\Cs_{1,k}(t) + \overline{\sigma}\Cs_{0,k}(t) \big),
  \qquad-\NF\leq k\leq\NF.
  \label{eq:legendre:fourier:current:density:b}
\end{align}
Taking the derivative in time of~\eqref{eq:legendre:fourier:Poisson:1}
and using~\eqref{eq:legendre:fourier:system} with $n=0$, we obtain the
Fourier representation of Ampere's equation:
\begin{align}
  \epsilon_0\left(\frac{2\pi\img}{L}k\right)\frac{d\Es_{k}}{\dt} 
  &= (\vb-\va)L\,\sum_{s\in\{e,i\}}\qs\frac{d\Cs_{0,k}}{\dt} \nonumber\\[0.5em]
  &= -(\vb-\va)L\,\sum_{s\in\{e,i\}}\qs
  \left(
    \left(\frac{2\pi\img}{L}k\right)\Big(\sigma_1\Cs_{1,k}+\overline{\sigma}\Cs_{0,k}\Big)
    +\frac{\qs}{\ms}\left[\Es\STAR\Diffv{\fs\phi_0}_{\va}^{\vb}\right]_{k}
  \right).
  \label{eq:legendre:fourier:Ampere:00}
\end{align}
For $k\neq 0$ and using
definition~\eqref{eq:legendre:fourier:current:density:b} we
reformulate Ampere's equation as
\begin{align}
  \epsilon_0\frac{d\Es_{k}}{\dt} 
  = -\sum_{s\in\{e,i\}}\big( \Jsk(t) + Q^s_k(t) \big)
  \label{eq:legendre:fourier:Ampere}
\end{align}
where 
\begin{align}
  Q^s_k = \left(\frac{2\pi\img}{L}k\right)^{-1}\,
  (\vb-\va)L\,\frac{(\qs)^2}{\ms}\left[\Es\STAR\Diffv{\fs\phi_0}_{\va}^{\vb}\right]_k.
  \label{eq:legendre:fouriere:Ampere:boundary:condition}
\end{align}
For $k=0$, the Fourier decomposition of Ampere's
equation~\eqref{eq:Ampere:equation:def} gives the consistency
condition $C_{A}=\JJ_{0}$, the zero-th Fourier mode of the total
current density $\JJ$.
  
\subsection{Collisional term}
To control the filamentation effect, we modify
system~\eqref{eq:legendre:fourier:compact} by introducing the
artificial collisional operator $\mathcal{C}(\Cvk)$ in the right-hand
side \cite{camporeale1}:
\begin{align}
  \frac{d\Cvk}{\dt} 
  +\left(\frac{2\pi\img}{L}k\right)\matA\Cvk
  -\frac{\qs}{\ms}\left[E\STAR\left(
      \matB\Cvs
      -\Diffv{\fs{\bm\phi}}_{\va}^{\vb}
    \right)\right]_k
  = \mathcal{C}(\Cvk).
   \label{eq:legendre:fourier:compact:collisional}
\end{align}
Consider the diagonal matrix $\matD^{s}_{\nu}$ whose $n$-th diagonal
entry is given by:
\begin{align}
  D^{s}_{n} := {\matD^{s}_{\nu}}_{|nn} = -\nu_s\frac{n(n-1)(n-2)}{(\NL-1)(\NL-2)(\NL-3)}
  \qquad n\geq 0,
  \label{eq:collisional:operator}
\end{align}
and where $\nu_s$ is an artificial diffusion coefficient whose value
can be different from species to species.
Then, the collisional term is given by
$\mathcal{C}(\Cvk)=\matD_{\nu}\Cvk$.
The effect of this operator is to damp the highest-modes of the
Legendre expansion, thus reducing the filamentation and avoiding
recurrence effects.
This operator is designed to be zero for $n=0,1,2$, in order not to
have any influence on the conservation properties of the method.

\newcommand{\Csb}{\overline{C}^s}
\newcommand{\Esb}{\overline{E}^s}
\newcommand{\fsb}{\overline{f}^s}

\newcommand{\Ct}{C^{\tau}}
\newcommand{\Et}{E^{\tau}}
\newcommand{\ft}{f^{\tau}}

\newcommand{\Ctt}{C^{\tau+1}}
\newcommand{\Ett}{E^{\tau+1}}
\newcommand{\ftt}{f^{\tau+1}}

\subsection{Crank-Nicolson time integration}
\label{sec:time:integration}
Let $\Delta t$ be the time step, $\tau$ the time index, and each
quantity superscripted by $\tau$ as taken at time $t^{\tau}=\tau\Delta
t$, e.g.,
$\Et=\Es(\cdot,t^{\tau})$,
$f^{s,\tau}=\fs(\cdot,\cdot,t^{\tau})$
$C^{s,\tau}_{n,k}=\Cs_{n,k}(t^{\tau})$, etc.
We advance the Legendre-Fourier coefficients $\Cs_{n,k}(t)$ in time by
the Crank-Nicolson time marching scheme \cite{cn}.
Omitting the superscript ``$s$'' in $C^{s,\tau}_{n,k}$ and $f^{s,\tau}$ to
ease the notation, Vlasov equation~\eqref{eq:legendre:fourier:system}
for each species 
and any Legendre-Fourier coefficient becomes:
\begin{align}
  &\frac{\Ctt_{n,k}-\Ct_{n,k}}{\Delta t} 
  +\frac{\pi\img}{L}k\,
  \Big(
  \sigma_{n+1}\big(\Ctt_{n+1,k} + \Ct_{n+1,k}\big) +
  \sigma_{n}  \big(\Ctt_{n-1,k} + \Ct_{n-1,k}\big) +
  \overline{\sigma}\big(\Ctt_{n,k} + \Ct_{n,k}\big)
  \Big)
  \nonumber\\[0.5em]
  &\qquad
  -\frac{\qs}{4\ms}
  \left[
    \big(\Ett+\Et\big)\STAR
    \left(
      \sum_{i=0}^{n-1}\sigma_{n,i}\big(\Ctt_{i}+\Ct_{i}\big)
      -\gamma^s\Diffv{ \big(\ftt+\ft\big)\phi_n }_{\va}^{\vb}
    \right)
  \right]_{k} = \mathcal{C}\left(\frac{1}{2}\Ctt_{n,k} + \frac{1}{2}\Ct_{n,k}\right).
  \label{eq:legendre:fourier:system:time}
\end{align}
Equation~\eqref{eq:legendre:fourier:system:time} provides an implicit
and non-linear system for the Legendre-Fourier coefficients
$\Cs_{n,k}(t)$ as each electric field mode $E_{k}^{\tau+1}$ for $k\neq
0$ depends on the unknown coefficient $\Cs_{0,k}(t^{\tau+1})$ that
must be evaluated at the same time $t^{\tau+1}$.
In practice, we apply a Jacobian-free Newton-Krylov
solver~\cite{knoll} to search for the minimizer of the residual given
by~\eqref{eq:legendre:fourier:system:time}.

\medskip
Consider the difference of the Fourier representation of Poisson's equation~\eqref{eq:legendre:fourier:Poisson:1} at 
times $t^{\tau}$ and $t^{\tau+1}$
\begin{align}
  \epsilon_{0}\left(\frac{2\pi\img}{L}k\right)\,
  \big(\Es^{\tau+1}_k-\Es^{\tau}_k\big) 
  = (\vb-\va)\,\sum_{s\in\{e,i\}}\qs
  \big(\Ctt_{0,k}-\Ct_{0,k}\big).
    \label{eq:Ampere:discrete:00}
\end{align}
By setting $n=0$ in~\eqref{eq:legendre:fourier:system:time},
recalling that $\phi_0=1$ and noting that the collisional term does not give 
any contribution, we find that
\begin{align}
  \frac{\Ctt_{0,k}-\Ct_{0,k}}{\Delta t}
  &+ \frac{\pi\img}{L}k\,
  \Big(
  \sigma_{1}\big(\Ctt_{1,k} + \Ct_{1,k}\big) +
  \overline{\sigma}\big(\Ctt_{1,k} + \Ct_{1,k}\big)
  \Big)
  \nonumber\\[0.5em]
  &-\frac{\qs}{4\ms}
  \left[
    \big(\Ett+\Et\big)\STAR
    \gamma^s\Diffv{ \big(\ftt+\ft\big) }_{\va}^{\vb}
  \right]_{k}=0.
  \label{eq:Ampere:discrete:05}
\end{align}
Using~\eqref{eq:Ampere:discrete:05} in~\eqref{eq:Ampere:discrete:00}
yields the discrete analog of Ampere's equation that is consistent with 
the full Crank-Nicolson based discretization of the Vlasov-Poisson system:
\begin{align}
  \epsilon_{0}\big(\Es^{\tau+1}_k-\Es^{\tau}_k\big)
  = -\frac{\Delta t}{2}\sum_{s\in\{e,i\}}
  \big(\Js_{k}(t^{\tau+1})+\Js_{k}(t^{\tau})\big)
  +\Delta t\mathcal{B}^{Amp}_{k}
  \label{eq:full:discrete:Ampere}
\end{align}
where we have introduced the explicit symbol
\begin{align}
  \mathcal{B}^{Amp}_{k} = 
  -\sum_{s\in\{e,i\}}
  \frac{\qs}{4\ms}\left(\frac{2\pi\img}{L}k\right)^{-1}
  \left[
    \big(\Ett+\Et\big)\STAR
    \gamma^s\Diffv{ \big(\ftt+\ft\big) }_{\va}^{\vb}
  \right]_{k}
  \qquad \textrm{for~}k\neq 0
  \label{eq:full:discrete:Ampere:BC}
\end{align}
to denote the boundary terms related to the behavior of 
all the distribution functions of the plasma species at the 
boundaries of the velocity domain.
In Section~\ref{sec:conservation:laws} we make use
of~\eqref{eq:full:discrete:Ampere} and~\eqref{eq:full:discrete:Ampere:BC}
to characterize the conservation of the total energy.

\section{Enforcing $L^2$ stability}
\label{sec:L2:stability}
The distribution function $\fs$ solving the Vlasov equation satisfies
the so-called $L^p$-stability property for $p\geq 1$.
To see this, just multiply equation~\eqref{eq:Vlasov} by
$p\fs(x,\vs,t)^{p-1}$ and integrate over the phase space domain
$[0,L]\times[\va,\vb]$.
Assuming that the velocity range is sufficiently large for having
$\fs(x,\va,t)=\fs(x,\vb,t)=0$, a simple calculation shows that
$d\slash{\dt}\norm{\fs(\cdot,\cdot,t)}_{L^p(\Omega)}^p=0$.
This property is particularly useful for $p=2$, which implies the
\emph{$L^2$ stability} of the method (sometimes called also ``energy
stability'' in the literature).
To derive a relation for the $L^2$ stability of the Legendre-Fourier
method, we need the result stated by the following lemma.
The proof of the lemma requires a few lengthy calculations and is
reported, for the sake of completeness, in appendix C.

\medskip
\begin{lemma}
  \label{lemma:L2stab:useful}
  Let $\Cvk$ be the vector containing the Legendre coefficients of the
  $k$-th Fourier mode of the distribution function $\fs$, $\Es$ the
  electric field and $\matB$ the matrix of coefficients defined
  in~\eqref{eq:legendre:matB:def}.
  Then, it holds that:
  \begin{align}
    2\sum_{k=-\NF}^{\NF}(\Cvk)^{\dagger}\Big[\Es\STAR\matB\Cvs\Big]_k 
    = \left[\Es\STAR\Diffv{(\fs)^2}_{\va}^{\vb}\right]_{0}
    = \sum_{k=-\NF}^{\NF}(\Cvk)^{\dagger}\Big[\Es\STAR\Diffv{\fs\bm\phi}_{\va}^{\vb}\Big]_k,
    \label{eq:legendre:fourier:L2stab:useful}
  \end{align}
  where $[\,\ldots\,]_{0}$ denotes the zero-th Fourier mode of the
  argument inside the brackets, and $\dagger$ denotes the conjugate
  transpose.
  All terms in~\eqref{eq:legendre:fourier:L2stab:useful} are real
  numbers.
\end{lemma}

\medskip
The $L^2$ stability of the Legendre-Fourier method depends on the
behavior of the distribution function $\fs$ at the boundaries
$\vs=\va$ and $\vs=\vb$.
This result is stated by the following theorem.

\medskip
\begin{theorem}
  \label{theo:L2stab}
  The coefficients of the Legendre-Fourier decomposition have the
  property that:
  \begin{align}
    \frac{d}{dt}\sum_{n=0}^{\NL-1}\sum_{k=-\NF}^{\NF}\abs{\Cs_{n,k}(t)}^2 
    = -\frac{\qs}{\ms}\left[\Es\STAR\Diffv{(\fs)^2}_{\va}^{\vb}\right]_{0}
    -2\sum_{n=0}^{\NL-1}\abs{D^{s}_{n}}\sum_{k=-\NF}^{\NF}\abs{\Cs_{n,k}(t)}^2.
    \label{eq:L2stab}
  \end{align}
\end{theorem}
\BEGINPROOF
Multiply~\eqref{eq:legendre:fourier:compact:collisional} from the left by
$(\Cvk)^{\dagger}$, the conjugate transpose of $\Cvk$, to obtain:
\begin{align*}
  (\Cvk)^{\dagger}\frac{d\Cvk}{\dt} 
  +\left(\frac{2\pi\img}{L}k\right)\,(\Cvk)^{\dagger}\matA\Cvk
  -\frac{\qs}{\ms}\,(\Cvk)^{\dagger}\left[E\STAR\left(
      \matB\Cvs
      -\Diffv{\fs{\bm\phi}}_{\va}^{\vb}
    \right)\right]_k
  = (\Cvk)^{\dagger}\matD^{s}_{\nu}\Cvk.
\end{align*}
Add to this equation 
its conjugate transpose.
Matrix $\matA$ is real and symmetric, $(\Cvk)^{\dagger}\matA\Cvk$ is
real and the spatial term cancels out from the equation.
Summing over the Fourier index $k$ we end up with:
\begin{align}
  \frac{d}{\dt}\sum_{k=-\NF}^{\NF}\abs{\Cvk}^2 
  = 2\frac{\qs}{\ms}\textsf{Re}\sum_{k=-\NF}^{\NF}
    (\Cvk)^{\dagger}\left[\,E\STAR
    \left(
      \matB\Cvs
      -\Diffv{\fs{\bm\phi}}_{\va}^{\vb}
    \right)\,
  \right]_{k}
  +2\textsf{Re}\big((\Cvk)^{\dagger}\matD^{s}_{\nu}\Cvk\big).
  \label{eq:legendre:fourier:L2stab:thm:10}
\end{align}
Since $\matD_{\nu}$ is a diagonal matrix with negative real entries $D^{s}_{n}$ and
$(\Cvk)^{\dagger}\matD^{s}_{\nu}\Cvk$ is a real quantity it holds that
\begin{align}
  \textsf{Re}\big((\Cvk)^{\dagger}\matD^{s}_{\nu}\Cvk\big) = 
  -\sum_{n=0}^{\NL-1}\abs{D^{s}_{n}}\sum_{k=-\NF}^{\NF}\abs{\Cs_{n,k}(t)}^2.
\end{align}
The assertion of the theorem follows by applying the result of
Lemma~\ref{lemma:L2stab:useful}.
\ENDPROOF

If $\fs(x,\va,t)=\fs(x,\vb,t)=0$ at the velocity boundaries (see
Remark~\ref{remark:zero:velocity:bcs}), then at any instant $t>0$ the
time derivative in~\eqref{eq:legendre:fourier:L2stab:thm:10} is
negative due to the collisional term and we have that
\begin{align*}
  \sum_{n=0}^{\NL-1}\sum_{k=-\NF}^{\NF}\abs{\Cs_{n,k}(t)}^2 \leq
  \sum_{n=0}^{\NL-1}\sum_{k=-\NF}^{\NF}\abs{\Cs_{n,k}(0)}^2.
\end{align*}
Note that in absence of the collisional term (take $\nu^s=0$
in~\eqref{eq:collisional:operator}) the time derivative is exactly
zero and $\ABS{\Cs_{n,k}(t)}$ is constant.
We refer to this property as \emph{the $L^2$ stability} because the
orthogonality of the Legendre and Fourier basis functions implies that
\begin{align}
  \norm{\fs(\cdot,\cdot,t)}_{L^2(\Omega)}^2 =
  \int_{0}^{L}\int_{\va}^{\vb}\abs{\fs(x,\vs,t)}^2\dv\dx 
  = (\vb-\va)L\,\sum_{n=0}^{\NL-1}\sum_{k=-\NF}^{\NF}\abs{\Cs_{n,k}(t)}^2,
  \label{eq:L2:stability:fs}
\end{align}
(see appendix B), from which we immediately find the $L^2$ stability
of the distribution function $\fs(x,\vs,t)$.
However, $\fs(x,\va,t)$ and $\fs(x,\vb,t)$ can be different than zero
and in general they are non zero since the Legendre
polynomials are globally defined on the whole domain and are non zero at the
velocity boundaries.
If the right-hand side of~\eqref{eq:L2stab} becomes positive, the
collisional term may be not enough to control the other term in the
right-hand side of~\eqref{eq:L2stab}.
Therefore, the method may become unstable and the time integration of
$\fs$ is arrested.

According to~\cite{Svard-Nordstrom:2014} we can enforce the stability
of the method by introducing the boundary conditions
$\fs(x,\va,t)=\fs(x,\vb,t)=0$ in weak form in the right-hand side of
system~\eqref{eq:legendre:fourier:compact:collisional} through the
penalty coefficient $\gamma^s$.
To this end, we modify
system~\eqref{eq:legendre:fourier:compact:collisional} as follows:
\begin{align}
  \frac{d\Cvk}{\dt} 
  +\left(\frac{2\pi\img}{L}k\right)\matA\Cvk
  -\frac{\qs}{\ms}\left[E\STAR\left(
      \matB\Cvs
      -\gamma^s\Diffv{\fs{\bm\phi}}_{\va}^{\vb}
    \right)\right]_k
  = \matD_{\nu}\Cvk.
  \label{eq:legendre:fourier:modified}
\end{align}
By suitably choosing the value of the penalty we minimize
or set equal to zero the term in the right-hand side
of~\eqref{eq:L2stab} that may cause the numerical instability.
This result is presented in the following theorem.

\medskip
\begin{theorem}
  \label{theo:L2-stab:modified}
  The modified form~\eqref{eq:legendre:fourier:modified} of the
  Legendre-Fourier method for solving the Vlasov-Poisson system is
  $L^2$-stable for $\gamma^s=1\slash{2}$ and any $\nu^s\geq 0$.
  The coefficients of the Legendre-Fourier decomposition have the
  property that:
  \begin{align}
    \frac{d}{dt}\sum_{n=0}^{\NL-1}\sum_{k=-\NF}^{\NF}\abs{\Cs_{n,k}(t)}^2 
    = 
    -2\sum_{n=0}^{\NL-1}D^{s}_{n}\sum_{k=-\NF}^{\NF}\abs{\Cs_{n,k}(t)}^2.
    \label{eq:stability:time:variation}
  \end{align}
\end{theorem}
\BEGINPROOF
Repeating the proof of Theorem~\ref{theo:L2stab} yields:
\begin{align}
  \frac{d}{dt}\sum_{n=0}^{\NL-1}\sum_{k=-\NF}^{\NF}\abs{\Cs_{n,k}(t)}^2 
  &= 
  2\frac{\qs}{\ms}\textsf{Re}
  \sum_{k=-\NF}^{\NF}
  (\Cvk)^{\dagger}\left[\,E\STAR
    \left(
      \matB\Cvs
      -\gamma^s\Diffv{\fs{\bm\phi}}_{\va}^{\vb}
    \right)\,
  \right]_{k}\nonumber\\[0.5em]
  &\,
  -2\sum_{n=0}^{\NL-1}\abs{D^{s}_{n}}\sum_{k=-\NF}^{\NF}\abs{\Cs_{n,k}(t)}^2.
  \label{eq:stability:time:variation}
\end{align}
Due to~\eqref{eq:legendre:fourier:L2stab:useful}, the first term of
the right-hand side of~\eqref{eq:stability:time:variation} is zero
(when the coefficient that multiplies $\gamma^s$ is non zero) by
setting
\begin{align*}
  \gamma^s = \frac
  { \sum_{k=-\NF}^{\NF}(\Cvk)^{\dagger}\left[\,E\STAR\matB\Cvs\right]_{k} } 
  { \sum_{k=-\NF}^{\NF}(\Cvk)^{\dagger}\left[\,E\STAR\Diffv{\fs{\bm\phi}}_{\va}^{\vb}\right]_{k} }
  = \frac{1}{2}.
\end{align*}
The assertion of the theorem is then proved by noting that any choice
of $\nu^s\geq 0$ in the collisional term makes the time derivative
non-positive.
\ENDPROOF

  The coefficient $\gamma^s$ in~\eqref{eq:legendre:fourier:modified}
  affects also the first three moment equations and eventually
  perturbs the conservation properties of the Vlasov-Poisson system.
  We may overcome this issue by considering the modified system
  \begin{align}
    \frac{d\Cvk}{\dt} 
    +\left(\frac{2\pi\img}{L}k\right)\matA\Cvk
    -\frac{\qs}{\ms}\left[E\STAR\left(
        \matB\Cvs
        -\matD_{\gamma^s}\Diffv{\fs{\bm\phi}}_{\va}^{\vb}
      \right)\right]_k
    = \matD_{\nu}\Cvk,
    \label{eq:legendre:fourier:modified:b}
  \end{align}
  where the penalty $\gamma^s$ is introduced through the diagonal
  matrix
  $\matD_{\gamma^s}=\textsf{diag}(0,0,0,\gamma^s,\ldots,\gamma^s)$ and
  does not change the conservation properties of the method.
  The penalty $\gamma^s$ can be determined at any time cycle by the
  formula:
  \begin{align*}
    \gamma^s = \frac
    { \sum_{k=-\NF}^{\NF}(\Cvk)^{\dagger}\left[\,E\STAR\matB\Cvs\right]_{k} } 
    { \sum_{n=3}^{\NL-1}\sum_{k=-\NF}^{\NF}\overline{C}^{s}_{nk}\left[\,E\STAR\Diffv{\fs{\bm\phi}}_{\va}^{\vb}\right]_{n,k} },
  \end{align*}
  where $\overline{C}^s_{n,k}$ is the conjugate of $\Csnk$, and the
  result of Theorem~\ref{theo:L2-stab:modified} still holds.
  Alternatively, we can apply $\gamma^s=\frac{1}{2}$ to all the
  Legendre modes except the first three, i.e., for $n=0,1,2$.
  This option is simpler to implement and computationally less
  expensive, but may not fix the stability issue of the method
  completely.
  Instead of equation~\eqref{eq:stability:time:variation}, it holds
  that
  \begin{align}
    \frac{d}{dt}\sum_{n=0}^{\NL-1}\sum_{k=-\NF}^{\NF}\abs{\Cs_{n,k}(t)}^2 
    =
    \frac{1}{2}\sum_{n=0}^{2}\sum_{k=-\NF}^{\NF}\overline{C}^s_{n,k}\left[\,E\STAR\Diffv{\fs{\bm\phi}}_{\va}^{\vb}\right]_{n,k}
    -2\sum_{n=0}^{\NL-1}\abs{D^{s}_{n}}\sum_{k=-\NF}^{\NF}\abs{\Cs_{n,k}(t)}^2,
    \label{eq:stability:time:variation:gr2}
  \end{align}
  and the first term in the right-hand side may still be a source of
  instability if it has the wrong sign.
  Nonetheless, if the dissipative effect of the collisional term
  in~\eqref{eq:stability:time:variation:gr2} is strong enough the scheme
  will remain stable.
  We investigated the effectiveness of this latter strategy in the
  numerical experiments of section~\ref{sec:numerical}.

\section{Conservation laws}
\label{sec:conservation:laws}
The Vlasov-Poisson model in the continuum setting is characterized by
the exact conservation of mass, momentum and energy.
The spectral discretization that is proposed in the previous section
reproduces these conservation laws in the discrete setting.
It turns out that the discrete analogs of the conservation of mass, momentum and energy 
depends on the variation in time of the Legendre-Fourier coefficients $\Cs_{n,k}$ for
$n=0,1,2$ and $k=0$, i.e., $\Cs_{0,0}$, $\Cs_{1,0}$, and $\Cs_{2,0}$.
The contribution of the second term
in~\eqref{eq:legendre:fourier:system} is zero when $k=0$ and the
transformed equation for the coefficients $\Cs_{n,0}$ (including the
stabilization factor $\gamma^s$ of Section~\ref{sec:L2:stability})
becomes:
\begin{align}
  \frac{d\Cs_{n,0}}{\dt} =
  \frac{\qs}{\ms}\left[\Es\STAR\left(\sum_{i=0}^{n-1}\sigma_{n,i}\Cs_{i}
      -\gamma^s\Diffv{\fs\phi_n}_{\va}^{\vb}\right)\right]_{0}.
  \label{eq:legendre:fourier:reduced}
\end{align}
In particular, we have:
\begin{align}
  \textrm{for~}n=0, k=0:\qquad & 
  \frac{d\Cs_{0,0}}{\dt} = -\gamma^s\frac{\qs}{\ms}\left[\Es\STAR\Diffv{\fs\phi_0}_{\va}^{\vb}\right]_{0},
  \label{eq:legendre:fourier:dert_C00}\\[1.em]
  \textrm{for~}n=1, k=0:\qquad & 
  \frac{d\Cs_{1,0}}{\dt} =  \frac{\qs}{\ms}\left[\Es\STAR\Big(\sigma_{1,0}\Cs_{0}-\gamma^s\Diffv{\fs\phi_1}_{\va}^{\vb}\Big)\right]_{0},
  \label{eq:legendre:fourier:dert_C10}\\[1.em]
  \textrm{for~}n=2, k=0:\qquad & 
  \frac{d\Cs_{2,0}}{\dt} =  \frac{\qs}{\ms}\left[\Es\STAR\Big(\sigma_{2,1}\Cs_{1}-\gamma^s\Diffv{\fs\phi_2}_{\va}^{\vb}\Big)\right]_{0}.
  \label{eq:legendre:fourier:dert_C20}
\end{align}
To derive the conservation laws for mass, momentum and energy for the
fully discrete approximation, we note that the analog of
equation~\eqref{eq:legendre:fourier:reduced} for $k=0$ becomes:
\begin{align}
  \frac{\Cs_{n,0}(t^{\tau+1})-\Cs_{n,0}(t^{\tau})}{\Delta t} 
  =
  \frac{\qs}{4\ms}
  \Bigg[\Big(\Es(\cdot,t^{\tau+1})+\Es(\cdot,t^{\tau})\Big)
  &\STAR
  \Bigg(\sum_{i=0}^{n-1}\sigma_{n,i}\Big(\Cs_{i}(t^{\tau+1})+\Cs_{i}(t^{\tau})\Big)
  \nonumber\\[0.5em]
  & -\gamma^s\Diffv{\big(\fs(\cdot,\vs,t^{\tau+1})+\fs(\cdot,\vs,t^{\tau})\big)\phi_n}_{\va}^{\vb}
  \Bigg)\Bigg]_{0}
  \label{eq:legendre:fourier:reduced:full}
\end{align}
as the collisional term is zero, and where $\Es(\cdot,t)$ and $\fs(\cdot,\vs,t)$ are the
electric field and the distribution function, respectively, as functions of $x$ for a given 
value of $\vs$ and $t$. 
By setting $n=0,1,2$ in~\eqref{eq:legendre:fourier:reduced:full} we
can also derive the analog of
equations~\eqref{eq:legendre:fourier:dert_C00}-\eqref{eq:legendre:fourier:dert_C20}
for the fully discrete approximation, which we omit.
In the following developments we consider the boundary term:
\begin{align}
  \mathcal{B}^{s;\tau,\tau+1}_{n,0} = 
    -\gamma^s\frac{\qs}{4}(\vb-\va)L\,\left[
    \big( \Es(\cdot,t^{\tau+1})+\Es(\cdot,t^{\tau}) \big)\STAR
    \Diffv{\big(\fs(\cdot,\vs,t^{\tau+1})+\fs(\cdot,\vs,t^{\tau})\big)\phi_n}_{\va}^{\vb}
  \right]_{0}.
\end{align}
Note that $\mathcal{B}^{s;\tau,\tau+1}_{n,0}=0$ when
$\fs(\cdot,\va,t^{\tau+z})=\fs(\cdot,\vb,t^{\tau+z})=0$ for $z=0,1$.

\subsection{Conservation of mass}
Using the Legendre-Fourier expansion of $\fs$ and the orthogonality
relations~\eqref{eq:Legendre:orthogonality} and~\eqref{eq:Fourier:orthogonality},
the total mass of the species $s$ is given by
\begin{align}
  \Ms(t) 
  = \ms\int_{0}^{L}\int_{\va}^{\vb}\fs(x,v,t)\,dv\,dx
  = \ms(\vb-\va)L\,C_{0,0}^s(t).
  \label{eq:legendre:fourier:M2}
\end{align}
By taking the time derivative of Eq~\eqref{eq:legendre:fourier:M2} and
using~\eqref{eq:legendre:fourier:dert_C00} it follows that
\begin{align}
  \frac{d\Ms}{\dt}
  =  \ms(\vb-\va)L\,\frac{d C_{0,0}^s(t)}{\dt}
  = -\gamma^s\,(\vb-\va)L\,\qs\left[\Es\STAR\Diffv{\fs\phi_0}_{\va}^{\vb}\right]_{0}.
 \label{eq:legendre:fourier:M3}
\end{align}
The conservation of the total mass per species includes a boundary
term that is zero if $\fs(x,\va,t)=\fs(x,\vb,t)=0$ (see
Remark~\ref{remark:zero:velocity:bcs}).

From~\eqref{eq:legendre:fourier:M2} and
using~\eqref{eq:legendre:fourier:reduced:full} with $n=0$, we derive
the conservation of the total mass per species in the full discrete model:
\begin{align}
  \Ms(t^{\tau+1})-\Ms(t^{\tau}) 
  = \ms(\vb-\va)L\,\big( C_{0,0}^s(t^{\tau+1}) - C_{0,0}^s(t^{\tau}) \big)
  = \Delta t\mathcal{B}^{s;\tau,\tau+1}_{0,0}.
  \label{eq:full:mass}
\end{align}
Equation~\eqref{eq:full:mass} states that the mass variation between
times $t^{\tau+1}$ and $t^{\tau}$ is balanced by the boundary term in
the right-hand side.

\subsection{Conservation of momentum}
\label{subsect:momentum:conservation}
The total momentum of the plasma is defined as
\begin{align}
  P(t) 
  = \sum_{s\in\{e,i\}}\Ps(t)
  = \sum_{s\in\{e,i\}}\ms\int_{0}^{L}\int_{\va}^{\vb}v\fs(x,v,t)\,dv\,dx,
  \label{eq:legendre:fourier:Pb}
\end{align}
where $\Ps(t)$ is the total momentum of the species $s$.
Introducing the Legendre-Fourier expansion of $\fs$, using the
integrated recursive formula~\eqref{eq:intg:vf}, orthogonality
relations~\eqref{eq:Legendre:orthogonality}
and~\eqref{eq:Fourier:orthogonality}, and mass
equation~\eqref{eq:legendre:fourier:M2} yield
\begin{align}
  \Ps(t)
  =\ms(\vb-\va)L\,\sigma_{1}\Cs_{1,0}(t)+\overline{\sigma}\Ms(t).
  \label{eq:legendre:fourier:P2}
\end{align}
Taking the time derivative of equation~\eqref{eq:legendre:fourier:P2}
and using~\eqref{eq:legendre:fourier:dert_C10} it follows that
\begin{align}
  \frac{d P}{\dt}
  &=\sum_{s\in\{e,i\}}\frac{d\Ps}{\dt}
  = \sum_{s\in\{e,i\}}\ms(\vb-\va)L\,\sigma_{1}\frac{d\Cs_{1,0}}{\dt}+\overline{\sigma}\frac{d\Ms}{\dt}\nonumber\\[0.5em]
  &= 
  \sum_{s\in\{e,i\}}
  \qs(\vb-\va)L\,\left(
    \sigma_{1}\sigma_{1,0}\,\big[\Es\STAR\Cs_0\big]_{0} -
    \gamma^s
    \left[\,\Es\STAR\left( 
        \overline{\sigma}\Diffv{\fs\phi_0}_{\va}^{\vb} + 
        \sigma_{1}\Diffv{\fs\phi_1}_{\va}^{\vb}\right)\,\right]_{0}
  \right).
  \label{eq:legendre:fourier:P3}
\end{align}
Using the Poisson equation the first term in the last right-hand side
is zero because the summation on the convolution index is on a
symmetric range of indices and the argument of the summation is
anti-symmetric:
\begin{align*}
  (\vb-\va)L\sum_{s\in\{e,i\}}\qs\big[\Es\STAR\Cs_0\big]_{0} 
  &= (\vb-\va)L\,\sum_{k=-\NF}^{\NF}\Es_{k}(t)\,\sum_{s\in\{e,i\}}\qs\Cs_{0,-k}(t)  \nonumber\\[0.5em]
  &= -2\pi\img\epsilon_0\,\sum_{k=-\NF}^{\NF}k\Es_{k}(t)\Es_{-k}(t)
  = 0.
\end{align*}
Consequently, equation~\eqref{eq:legendre:fourier:P3} becomes:
\begin{align}
  \frac{d P}{\dt} = -(\vb-\va)L\,\sum_{s\in\{e,i\}}\qs\gamma^s
  \left[
    \,\Es\STAR
    \left(
      \overline{\sigma}\Diffv{\fs\phi_0}_{\va}^{\vb}+\sigma_{1}\Diffv{\fs\phi_1}_{\va}^{\vb}
    \right)
    \,
  \right]_{0}.
  \label{eq:legendre:fourier:P4}
\end{align}
The conservation of the total momentum includes a boundary term that
is zero if $\fs(x,\va,t)=\fs(x,\vb,t)=0$ (see
Remark~\ref{remark:zero:velocity:bcs}).

From~\eqref{eq:legendre:fourier:P2} and
using~\eqref{eq:legendre:fourier:reduced:full} with $n=1$, we derive
the variation of momentum per species $s$ between times
$t^{\tau}$ and $t^{\tau+1}$:
\begin{align*}
  \Ps(t^{\tau+1})-\Ps(t^{\tau}) 
  &= \ms(\vb-\va)L\,\sigma_{1}\Big(\Cs_{1,0}(t^{\tau+1})-\Cs_{1,0}(t^{\tau})\Big)
  +\overline{\sigma}\Big(\Ms(t^{\tau+1})-\Ms(t^{\tau})\Big)
  \nonumber\\[0.5em]
  &= \frac{\qs}{4}(\vb-\va)L\,\sigma_{1}\Delta t\,
  \Bigg[\Big(\Es(\cdot,t^{\tau+1})+\Es(\cdot,t^{\tau})\Big)\STAR
  \Big(\sigma_{1,0}\big(\Cs_{0}(t^{\tau+1})+\Cs_{0}(t^{\tau})\big)\Big)\Bigg]_{0}
  \nonumber\\[0.5em]
  &\quad+\Delta t\left(
  \sigma_{1}\mathcal{B}^{s;\tau,\tau+1}_{1,0}
  +\overline{\sigma}\mathcal{B}^{s;\tau,\tau+1}_{0,0}
  \right).
\end{align*}
Furthermore, summing over all the species, taking the zero-th Fourier
mode of the convolution product, and using the Poisson equation yield:
\begin{align*}
  &(\vb-\va)L\,\sum_{s\in\{e,i\}}\qs
  \Bigg[\Big(\Es(\cdot,t^{\tau+1})+\Es(\cdot,t^{\tau})\Big)\STAR
  \Big(\sigma_{1,0}\big(\Cs_{0}(t^{\tau+1})+\Cs_{0}(t^{\tau})\big)\Big)\Bigg]_{0} 
  \nonumber\\ 
  &\qquad\qquad=
  \sum_{k=-\NF}^{\NF}\big(E(\cdot,t^{\tau+1})+E(\cdot,t^{\tau})\big)_{k}
  \sum_{s\in\{e,i\}}\qs\,(\vb-\va)L\,\big(\Cs_{0}(t^{\tau+1})+\Cs_{0}(t^{\tau})\big)_{-k}
  \nonumber\\ 
  &\qquad\qquad=
  -2\pi i\epsilon_0
  \sum_{k=-\NF}^{\NF}k
  \big(E(\cdot,t^{\tau+1})+E(\cdot,t^{\tau})\big)_{k}
  \big(E(\cdot,t^{\tau+1})+E(\cdot,t^{\tau})\big)_{-k}
  = 0.
\end{align*}
Therefore, in the full discrete model the \emph{conservation of the
  total momentum} holds in the form:
\begin{align}
  P(t^{\tau+1})-P(t^{\tau}) = 
  \sum_{s\in\{e,i\}}\big(\Ps(t^{\tau+1})-\Ps(t^{\tau})\big) = 
  \Delta t\left(
  \sigma_{1}\mathcal{B}^{s;\tau,\tau+1}_{1,0}
  +\overline{\sigma}\mathcal{B}^{s;\tau,\tau+1}_{0,0}
  \right),
  \label{eq:full:momentum}
\end{align}
which states that the variation of the total momentum between times $t^{\tau}$
and $t^{\tau+1}$ is balanced by the boundary terms in the right-hand side
of~\eqref{eq:full:momentum}.

\subsection{Conservation of energy}
\label{subsect:energy:conservation}
The total energy of the plasma is defined as
\begin{align}
  \ENERGY_{tot}(t)
  = \sum_{s\in\{e,i\}}\Ensk(t) + \Epot(t)
  = \sum_{s\in\{e,i\}}\frac{m_s}{2}\int_{0}^{L}\int_{\va}^{\vb}v^2\fs(x,v,t)\,dv\,dx
  + \frac{\epsilon_0}{2}\int_{0}^{L}E(x,t)^2\,dx,
\end{align}
where $\Ensk(t)$ and $\Epot(t)$ are the kinetic energy of the species
$s$ and the potential energy at time $t$, respectively.
Introducing the Legendre-Fourier expansion of $\fs$ and using the
orthogonality
relations~\eqref{eq:Legendre:orthogonality} 
and~\eqref{eq:Fourier:orthogonality}, the kinetic energy of
species $s$ is reformulated as:
\begin{align}
  \Ensk(t)
  =
  \frac{\ms}{2}(\vb-\va)L\,\Big(
  \sigma_{2}\sigma_{1}\,\Cs_{2,0}(t) +
  2\sigma_{1}\overline{\sigma}\,\Cs_{1,0}(t) +
  \big( \sigma_{1}^2+\sigma_{0}^2+\overline{\sigma}^2 \big)\Cs_{0,0}(t)
  \Big).
  \label{eq:legendre:fourier:E2}
\end{align}
 
We take the derivative in time of $\Ensk(t)$ and
use~\eqref{eq:legendre:fourier:dert_C00}-\eqref{eq:legendre:fourier:dert_C20}
to obtain
\begin{align}
  \frac{d\Ensk}{\dt}
  &=
  \frac{\ms}{2}(\vb-\va)L\,\Big(
  \sigma_{2}\sigma_{1}\,\frac{d\Cs_{2,0}}{\dt} +
  2\sigma_{1}\overline{\sigma}\,\frac{d\Cs_{1,0}}{\dt} +
  \big( \sigma_{1}^2+\sigma_{0}^2+\overline{\sigma}^2 \big)\frac{d\Cs_{0,0}}{\dt}
  \Big)
  \nonumber\\[0.5em]
  &=
  \frac{\qs}{2}(\vb-\va)L\,\left[\Es\STAR\Big(
  \sigma_{1}\sigma_{2}\sigma_{2,1}\,\Cs_{1} +
  2\overline{\sigma}\sigma_{1}\sigma_{1,0}\,\Cs_{0}
  \Big)\right]_{0}
  +\mathcal{B}^s_{kin}
  \label{eq:kinetic:energy:00}
\end{align}
where we introduced the ``kinetic'' boundary term per species $s$:
\begin{align*}
  \mathcal{B}^{s}_{kin} = -\gamma^s\frac{\qs}{2}(\vb-\va)L\,
  \left[\Es\STAR\Big(\,
    \sigma_1\sigma_2\Diffv{\fs\phi_2}_{\va}^{\vb} +
    2\sigma_1\overline{\sigma}\Diffv{\fs\phi_1}_{\va}^{\vb} +
    (\sigma_1^2+\sigma_0^2+\overline{\sigma}^2)\Diffv{\fs\phi_0}_{\va}^{\vb}
    \,\Big)\right]_{0}.
\end{align*}
As $\sigma_{2}\sigma_{2,1}=2\sigma_{1}\sigma_{1,0}=2$, and
applying~\eqref{eq:legendre:fourier:current:density:b}
to~\eqref{eq:kinetic:energy:00}, we obtain:
\begin{align}
  \frac{d\Ensk}{\dt}
  = \qs(\vb-\va)L\,\Big[\Es\STAR\Big(\sigma_{1}\Cs_{1}+\overline{\sigma}\Cs_{0}\Big)\Big]_{0} 
  + \mathcal{B}^{s}_{kin}
  = \Big[\Es\STAR\Js\Big]_{0} + \mathcal{B}^{s}_{kin}.
  \label{eq:legendre:fourier:E3}
\end{align}
Using~\eqref{eq:E:Fourier}, the orthogonality
relation~\eqref{eq:Fourier:orthogonality} and the
convolution notation, the potential energy of the electric field is
given by:
\begin{align}
  \Epot(t) 
  = \frac{\epsilon_0}{2}\int_{0}^{L}\Es(x,t)^2\dx
  = \frac{\epsilon_0}{2}\sum_{k=-\NF}^{\NF}E_k(t)E_{-k}(t)
  = \frac{\epsilon_0}{2}\big[\Es\STAR\Es\big]_{0}.
  \label{eq:Epot:Fourier}
\end{align}
Then, we take the time derivative of the equation above, use
Ampere's equation~\eqref{eq:legendre:fourier:Ampere} and note that
$[\,\Es\,\STAR\,C_{A}\,]_0=0$ as the average of $E$ on $[0,L]$ is zero to
obtain:
\begin{align}
  \frac{d\Epot}{\dt} 
  = \epsilon_0\left[\Es\STAR\frac{\partial\Es}{\partial t}\right]_{0}
  = -\Bigg[\Es\STAR\Bigg(
  \sum_{s\in\{e,i\}}\big(\Js+\gamma^sQ^s\big)
  +C_{A}
  \Bigg)
  \Bigg]_{0}
  = -\Bigg[\Es\STAR\sum_{s\in\{e,i\}}\Js\Bigg]_{0} + \mathcal{B}_{pot}
\end{align}
where, after expanding the convolution product, we introduced the
symbol 
\begin{align*}
  \mathcal{B}_{pot}=-\sum_{k=-\NF}^{\NF}\Es_{k}\sum_{s\in\{e,i\}}\gamma^sQ^s_{-k},
\end{align*}
for the ``potential'' boundary term, $Q^s_{k}$ being the boundary term defined in 
Ampere's equation~\eqref{eq:legendre:fouriere:Ampere:boundary:condition}.
Adding the total kinetic energy for all species and the potential energy gives:
\begin{align*}
  \frac{d\mathcal{E}_{tot}}{\dt} = 
  \sum_{s\in\{e,i\}}\frac{d\Ensk(t)}{\dt} + \frac{d\Epot}{\dt} = 
  \sum_{s\in\{e,i\}}\mathcal{B}^{s}_{kin} + \mathcal{B}_{pot}.
\end{align*}
The conservation of the total energy includes a boundary term that is
zero if $\fs(x,\va,t)=\fs(x,\vb,t)=0$ (see
Remark~\ref{remark:zero:velocity:bcs}).

\medskip
From~\eqref{eq:legendre:fourier:E2}, the variation of the kinetic
energy $\Ensk(t)$ between times $t^{\tau}$ and $t^{\tau+1}$ reads as:
\begin{align*}
  \Ensk(t^{\tau+1})-\Ensk(t^{\tau})
  &=
  \frac{\ms}{2}(\vb-\va)L\,\Big(
  \sigma_{2}\sigma_{1}\,\big(\Cs_{2,0}(t^{\tau+1}) - \Cs_{2,0}(t^{\tau})\big)
  \nonumber\\[0.5em]
  &\qquad
  +2\sigma_{1}\overline{\sigma}\,\big(\Cs_{1,0}(t^{\tau+1}) - \Cs_{1,0}(t^{\tau})\big) +
  \big( \sigma_{1}^2+\sigma_{0}^2+\overline{\sigma}^2 \big)\,\big(\Cs_{0,0}(t^{\tau+1}) - \Cs_{0,0}(t^{\tau})\big)
  \Big).
\end{align*}
Using~\eqref{eq:legendre:fourier:reduced:full} with $n=2$ yields:
\begin{align*}
  \Ensk(t^{\tau+1})-\Ensk(t^{\tau})
  &=
  \frac{\qs}{8}(\vb-\va)L\Delta t\,
  \Bigg[\,
  \Big(\Es(\cdot,t^{\tau+1})+\Es(\cdot,t^{\tau}\Big)\STAR
  \Big(
  \sigma_{21}\sigma_{2}\sigma_{1}\big(\Cs_{1}(t^{\tau+1})+\Cs_{1}(t^{\tau})\big)
  \nonumber \\[0.5em]
  &
  +2\sigma_{10}\sigma_{1}\overline{\sigma}\big(\Cs_{0}(t^{\tau+1})+\Cs_{0}(t^{\tau})\big)
  \Big)
  \Bigg]_{0} + 
  \Delta t\mathcal{B}^{s;\tau,\tau+1}_{kin},
\end{align*}
where
\begin{align*}
  \mathcal{B}^{s;\tau,\tau+1}_{kin} =
  \sigma_{2}\sigma_{1}\mathcal{B}^{s;\tau,\tau+1}_{2,0} +
  2\sigma_{1}\overline{\sigma}\mathcal{B}^{s;\tau,\tau+1}_{1,0}
  +\big( \sigma_{1}^2+\sigma_{0}^2+\overline{\sigma}^2 \big)\mathcal{B}^{s;\tau,\tau+1}_{0,0}.
\end{align*}
Noting that $\sigma_{21}\sigma_{2}=2\sigma_{10}\sigma_{1}=2$,
using the definition of the convolution product $\STAR$, the
Fourier decomposition of the electric field and the Legendre coefficients,
and the definition of the Fourier coefficients of the current density $\Jsk(t)$ given 
in~\eqref{eq:legendre:fourier:current:density:b}
yield:
\begin{align}
  \Ensk(t^{\tau+1})-\Ensk(t^{\tau}) =
  \frac{\Delta t}{4}\sum_{k=-\NF}^{\NF}\Big(\Ett_{-k}+\Et_{-k}\Big)
  \Big(\Jsk(t^{\tau+1})+\Jsk(t^{\tau})\Big) 
  +\Delta t\mathcal{B}^{s;\tau,\tau+1}_{kin}.
  \label{eq:full:energy:00}
\end{align}
From~\eqref{eq:Epot:Fourier}, the variation of the potential energy
between times $t^{\tau+1}$ and $t^{\tau}$ is given by:
\begin{align*}
  \Epot(t^{\tau+1})-\Epot(t^{\tau}) 
  &= 
  \frac{\epsilon_0}{2}
  \big[\Es(\cdot,t^{\tau+1})\STAR\Es(\cdot,t^{\tau+1})\big]_{0} -
  \frac{\epsilon_0}{2}
  \big[\Es(\cdot,t^{\tau})\STAR\Es(\cdot,t^{\tau})\big]_{0}
  \nonumber\\[0.5em]
  &= \frac{\epsilon_0}{2}\Big[
  \big(\Es(\cdot,t^{\tau+1})-\Es(\cdot,t^{\tau})\big)\STAR\big(\Es(\cdot,t^{\tau+1})+\Es(\cdot,\tau)\big)
  \Big]_{0}
  \nonumber\\[0.5em]
  &
  = \frac{\epsilon_0}{2}
  \sum_{k=-\NF}^{\NF}
  \Big(\Es^{\tau+1}_{k}-\Es^{\tau}_{k}\Big)\,
  \Big(\Es^{\tau+1}_{-k}+\Es^{\tau}_{-k}\Big).
\end{align*}
Using the discrete analog of Ampere's equation given by~\eqref{eq:full:discrete:Ampere} 
and~\eqref{eq:full:discrete:Ampere:BC} yields:
\begin{align}
  \Epot(t^{\tau+1})-\Epot(t^{\tau})
  =
  -\frac{\Delta t}{4}
  \sum_{k=-\NF}^{\NF}
  \Big(\Ett_{-k}+\Et_{-k}\Big)
  \sum_{s\in\{e,i\}}\Big(\Jsk(t^{\tau+1})+\Jsk(t^{\tau})\Big)
  +\Delta t\mathcal{B}^{\tau,\tau+1}_{pot}
\end{align}
where
\begin{align}
\mathcal{B}^{\tau,\tau+1}_{pot} =
  \frac{1}{2}
  \sum_{k=-\NF}^{\NF}
  \Big(\Ett_{-k}+\Et_{-k}\Big)
  \Big( 
  \mathcal{B}^{Amp}_{k}(t^{\tau+1}) +
  \mathcal{B}^{Amp}_{k}(t^{\tau})
  \Big).
  \label{eq:full:energy:05}
\end{align}
Finally, we add the kinetic energy terms for $s\in\{e,i\}$
in~\eqref{eq:full:energy:00} and the potential
energy~\eqref{eq:full:energy:05} to find the relation
expressing the
\emph{total energy conservation} for the full discrete approximation:
\begin{align}
  \ENERGY_{tot}(t^{\tau+1})-\ENERGY_{tot}(t^{\tau})
  =
  \Delta t\sum_{s\in\{e,i\}}\mathcal{B}^{s;\tau,\tau+1}_{kin}
  +\Delta t\mathcal{B}^{\tau,\tau+1}_{pot}.
  \label{eq:full:energy:10}
\end{align}
Equation~\eqref{eq:full:energy:10} 
states that the variation of the total energy between times $t^{\tau}$
and $t^{\tau+1}$ is balanced by the proper combination of kinetic and potential boundary terms 
in the right-hand side and
expresses the \emph{conservation of the total energy} for the full discretization
of the Vlasov-Poisson system.

\section{Numerical experiments}
\label{sec:numerical}

In this section we assess the computational performance of the
Legendre-Fourier method by solving the Landau damping, two-stream
instability and ion acoustic wave problems.
These test cases are classical problems in plasma physics and are
routinely used to benchmark kinetic codes.
In our numerical experiments, we are mainly interested in showing the
conservation properties of the method, i.e., the discrepancy between
the initial value of mass, momentum and energy, and their value at
successive instants in time during the simulation.
We also investigate the stability of the method, i.e., how the
$L^2$-norm of the distribution function defined as
in~\eqref{eq:L2:stability:fs} changes during the time evolution of the
system.
The penalty $\gamma^e$ is applied to all Legendre modes except the
first three and the stability of the Legendre-Fourier method is
ensured by the artificial collisional term when $\nu^e=1$.
This strategy, which is discussed at the end of
section~\ref{sec:L2:stability}, is very effective in providing a
stable method with good conservation properties.
In the two-stream instability problem, we also investigate the effect
of applying penalty $\gamma^e$ on all the moment equations on the
conservation of the total energy.

In the first two test problems, the ions constitute a fixed background with
density $\rho^i(x,t)=1$.

We also introduce the following normalization: time is normalized on
the electron plasma frequency $\omega_{pe}$; position $x$ on the
electron Debye length $\lambda_D$; velocity $v$ on the electron
thermal velocity $v_{te}=\sqrt{kT_e/m_e}$ where $k$ is the Boltzmann
constant, $T_e$ the electron temperature and $m_e$ the electron mass;
the electric field $E$ on $m_ev_{te}\omega_{pe}\slash{e}$, where $e$
is the elementary charge; species densities on a reference density
$n_0$; and the species distribution function on $v_{te}/n_0$.

\subsection{Landau damping}

Landau damping is a classical kinetic effect in warm plasmas, due to
particles in resonance with an initial wave perturbation.
This interaction leads to an exponential decay of the electric field perturbation.
This problem is particularly challenging for kinetic codes because of
the continuous filamentation in velocity space, which is a
characteristic feature of the collision-less plasma described by
the Vlasov equation.
Filamentation is controlled by the artificial collisional operator
introduced in~\eqref{eq:collisional:operator}.

The initial distribution of the electrons is given by
\begin{equation}
  f^e(x,\vs,t=0) = \frac{1}{\sqrt{2\pi}}e^{-\frac{v^2}{2}}\,\left[1+\varepsilon\cos\Big(\frac{2\pi}{L}kx\Big)\right],
  \label{finit}
\end{equation}
with $k=1$ and $\varepsilon=10^{-3}$.
The Legendre-Fourier expansion of Eq. (\ref{finit}) implies that the 
modes $C_{n,0}$, $C_{n,k}$ and $C_{n,-k}$ are excited at $t=0$.

In this test case, the final simulation time is $T=100$ with time
step $\Delta t=0.05$, $\NL=201$ Legendre modes and $2\NF+1=51$ Fourier
modes.
The domain of integration is set to $L=2\pi$, $\vb=-\va=5$.

\begin{figure}
  \begin{center}
    \begin{overpic}[scale=0.8]{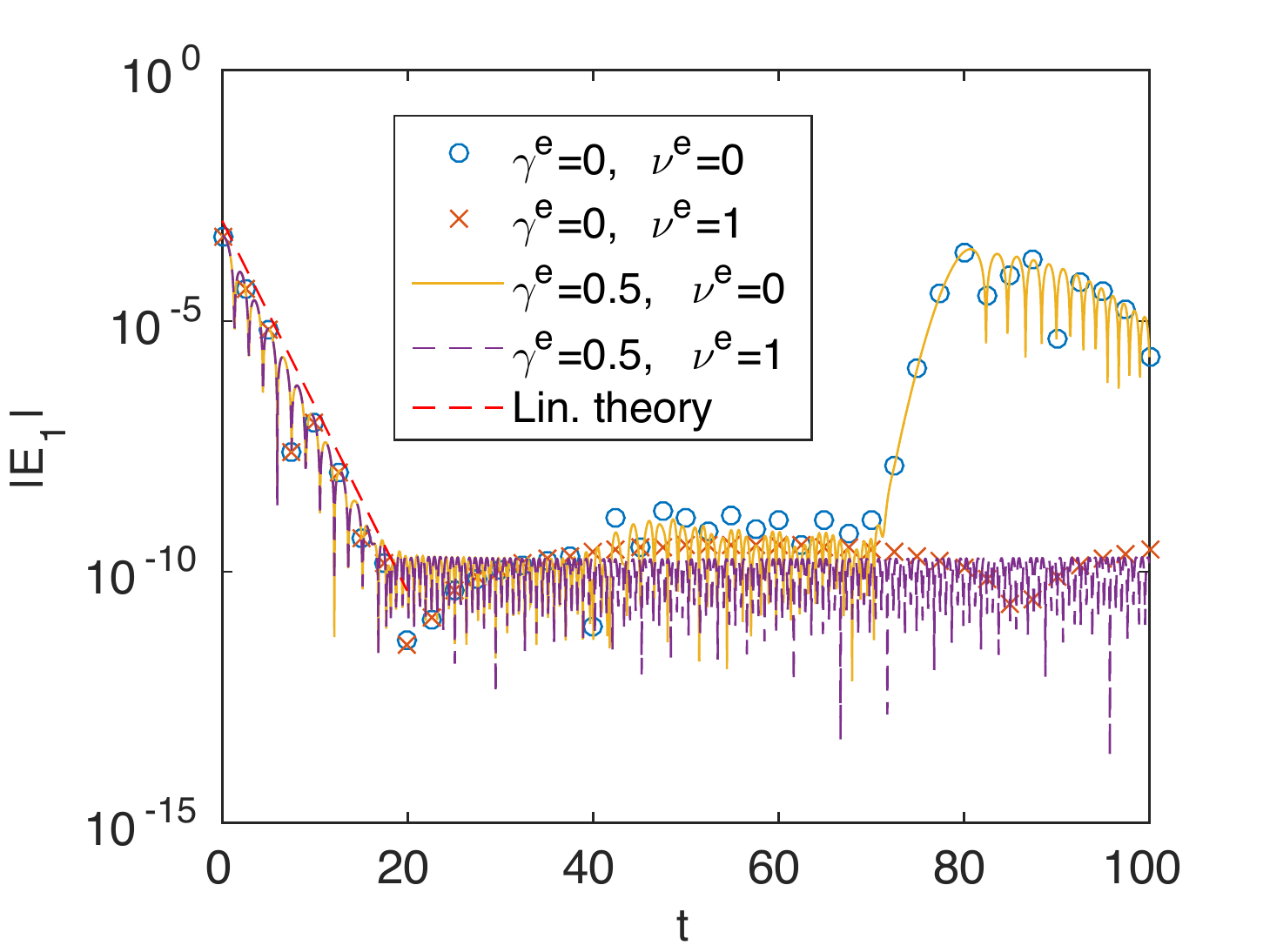}
      \put(48,-4){\begin{large}\textbf{Time}\end{large}}
      \put(-5,36){\begin{sideways}\begin{large}$\mathbf{\abs{E_1}}$\end{large}\end{sideways}}
    \end{overpic}
    \vspace{0.5cm}
    \caption{Landau damping test: first Fourier mode of the electric
      field versus time for the four combinations of
      $\gamma^e\in\{0,\,0.5\}$ and $\nu^e\in\{0,1\}$. Penalty
      $\gamma^e$ is not applied to the equations of the first three
      Legendre modes.}
  \label{fig:LD:00}
  \end{center}
\end{figure}

Figure~\ref{fig:LD:00} shows the first mode of the electric field
$|E_1|$ versus time for two different values of the stabilization
parameter ($\gamma^e\in\{0,\,0.\}$) and the collisional frequency
($\nu^e\in\{0,\,1\}$).
For all cases the damping rate is in good agreement with the Landau
damping theory, which predicts $\gamma_{damp}=-0.85$.
One can also notice that for all cases the simulation is stable,
regardless of the value of $\gamma^e$, and that $\gamma^e$ does not
really affect much the dynamics.
As expected, when $\nu^e=0$ the system exhibits recursive behavior.
The collisional operator with $\nu^e=1$ is however sufficient to
remove the recurrence effect and $|E_1|$ stabilizes around $10^{-10}$
for $t>20$.

\begin{figure}
  \begin{center}
    \begin{overpic}[scale=0.5]{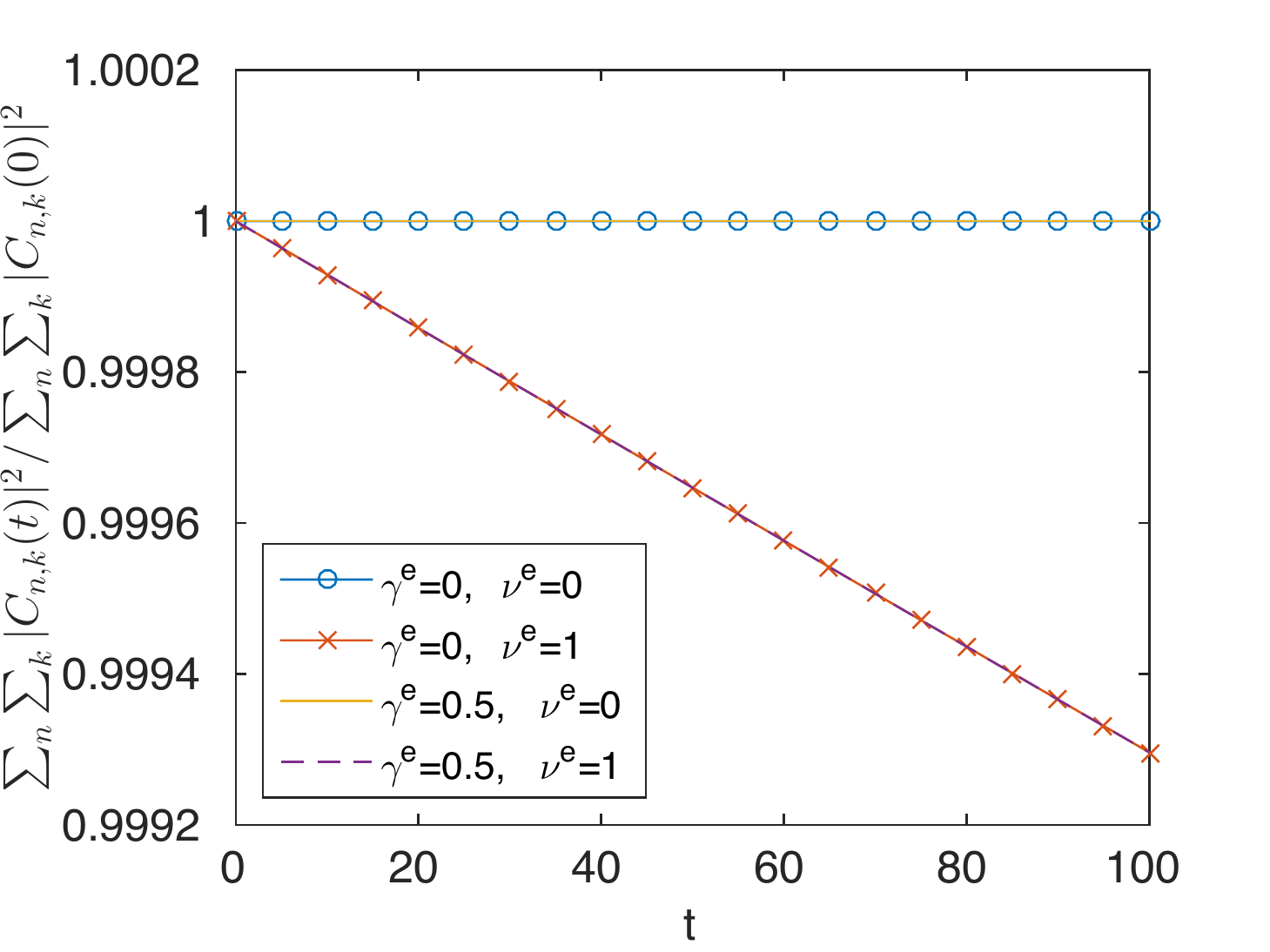}
      \put(48,-6){\textbf{Time}}
      \put(-8,21){\begin{sideways}\textbf{$L^2$ norm of $\fs$}\end{sideways}}
    \end{overpic}
    \hspace{0.75cm}
    \begin{overpic}[scale=0.5]{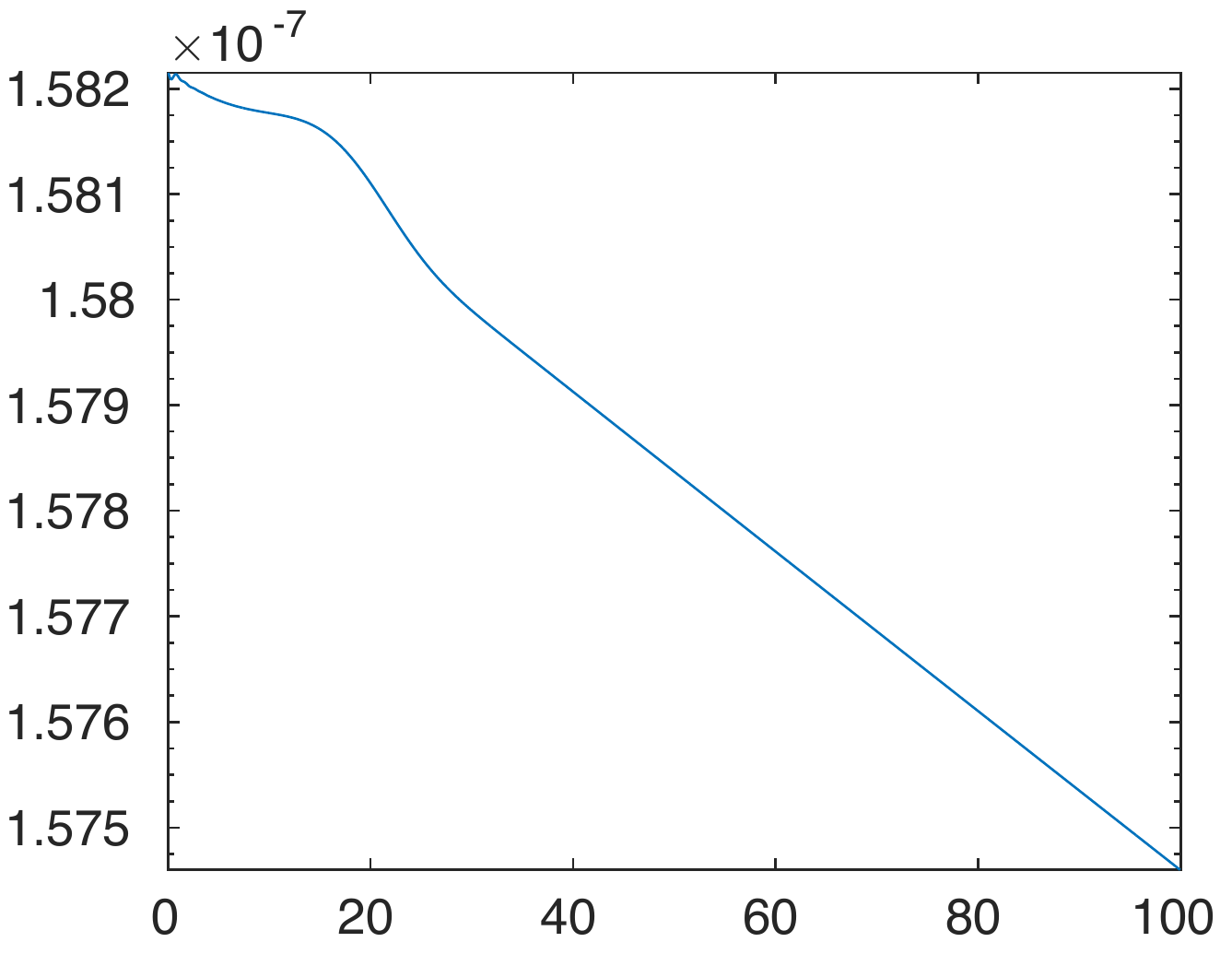}
      \put(48,-6){\textbf{Time}}
      \put(-8,13){\begin{sideways}
      \textbf{Variation of $\norm{\fs}_{L^2(\Omega)}$}
      \end{sideways}}
    \end{overpic}
    \vspace{0.5cm}
    \caption{Landau damping test: the left panel shows the $L2$ norm of 
      the electron distribution function $f^e$, cf.~\eqref{eq:relative:L2-norm:fs},
      versus time for the four combinations of
      $\gamma^e\in\{0,\,0.5\}$ and $\nu^e\in\{0,1\}$.
      The right panel shows the time variation of the same quantity
      versus time as predicted by Theorem~\ref{theo:L2stab} for
      $\gamma^e=0.5$ and $\nu^e=1$. 
      The time derivative is approximated by central finite
      differences. Penalty $\gamma^e$ is not applied to the equations
      of the first three Legendre modes.}
  \label{fig:LD:02}
  \end{center}
\end{figure}

Figure~\ref{fig:LD:02} (left) shows the time evolution of 
$\norm{f^e(\cdot,\cdot,t)}_{L^2(\Omega)}^2/\norm{f^e(\cdot,\cdot,0)}_{L^2(\Omega)}^2$, which is
normalized to its value at time
$t=0$, for the same cases of Fig.~\ref{fig:LD:00}.
According to\eqref{eq:L2:stability:fs}, this quantity
is computed as
\begin{align}
  \frac{\norm{f^e(\cdot,\cdot,t)}_{L^2(\Omega)}^2}{\norm{f^e(\cdot,\cdot,0)}_{L^2(\Omega)}^2} =
  \frac{ \sum_{n=0}^{\NL-1}\sum_{k=-\NF}^{\NF}\abs{C^e_{n,k}(t)}^2 }{ \sum_{n=0}^{\NL-1}\sum_{k=-\NF}^{\NF}\abs{C^e_{n,k}(0)}^2 }.
  \label{eq:relative:L2-norm:fs}
\end{align}

When $\nu^e=0$, the $L^2$ norm of $f^e$ is constant on the scale of
the plot and the boundary term in~\eqref{eq:L2stab} has a rather
negligible effect.
Instead, when $\nu^e=1$, the $L^2$ norm of $f^e$
decreases with an almost constant slope since the collisional term
in~\eqref{eq:L2stab} is dominant.
Figure~\ref{fig:LD:02} (right) shows that Theorem~\ref{theo:L2stab}
[equation~\eqref{eq:L2stab}] is indeed satisfied numerically.
In Fig.~\ref{fig:LD:02} (right) the time derivative is computed by
central finite differences.

\begin{figure}
  \begin{center}
    \begin{overpic}[scale=0.8]{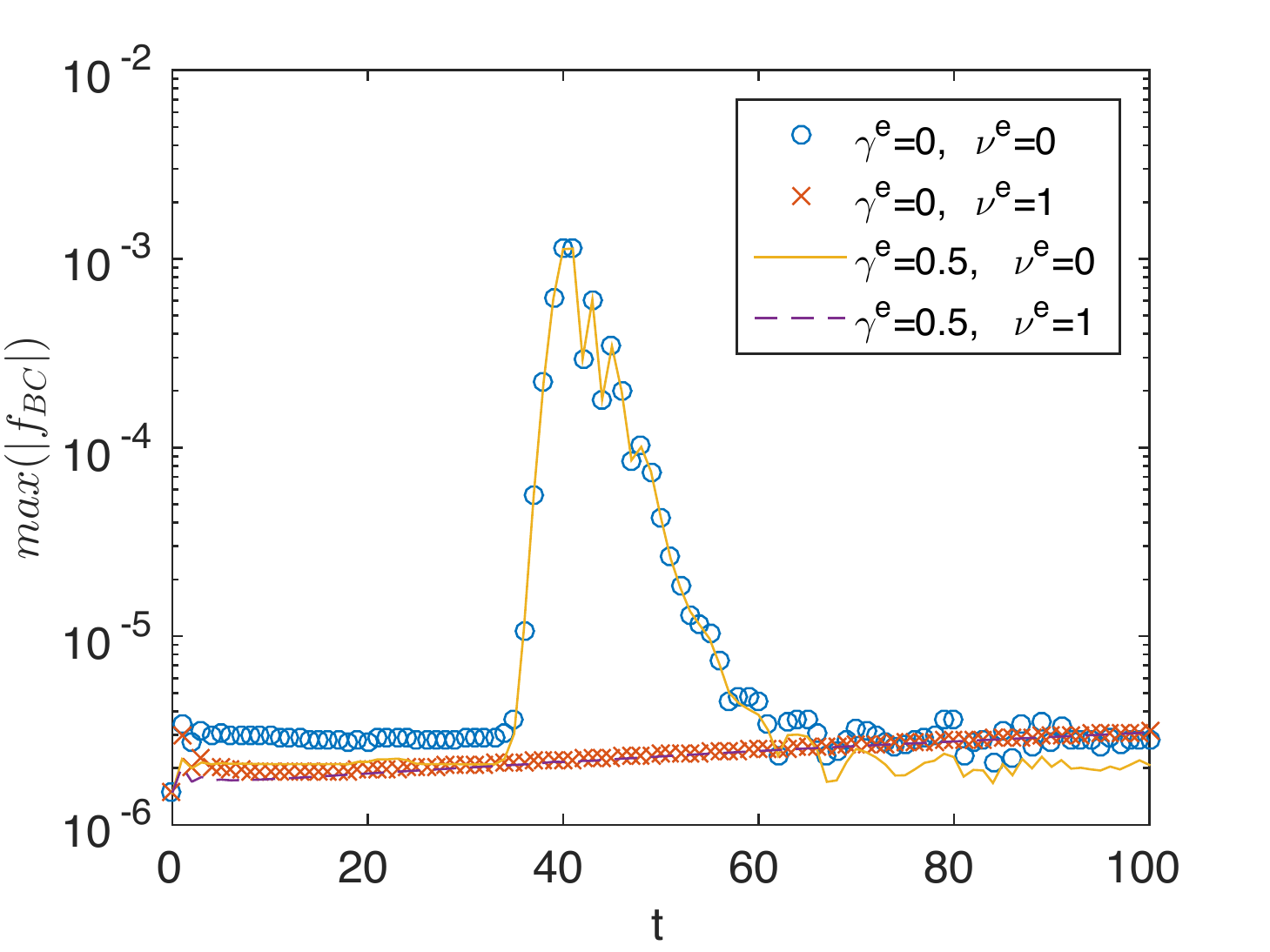}
      \put(48,-5){\FLARGE{\textbf{Time}}}
      \put(-5,29){\begin{sideways}
      \FLARGE{\textbf{max($\abs{\mathbf{f^e_{BC}}}$)}}
      \end{sideways}}
    \end{overpic}
    \vspace{0.5cm}
    \caption{Landau damping test: maximum value of the electron
      distribution function $f^e$ at the boundaries of the velocity
      range for the four combinations of $\gamma^e\in\{0,\,0.5\}$ and
      $\nu^e\in\{0,\,1\}$.  Penalty $\gamma^e$ is not applied to the
      equations of the first three Legendre modes.}
  \label{fig:LD:04}
  \end{center}
\end{figure}

Figure~\ref{fig:LD:04} shows the time evolution of the maximum value
of the distribution function at the boundary of the system
$v=\vb,\,\va$: $\max\abs{f_{BC}}=\max(\abs{f^e(\vb)},\abs{f^e(\va)})$,
with the same format of Fig.~\ref{fig:LD:00}.
One can notice the beneficial effect of the collisional operator: when
$\nu^e=0$ there is a sharp increase of $\max\abs{f_{BC}}$ around $t\sim
40$, while for $\nu^e=1$ it holds that $\max\abs{f_{BC}}$ approximately
$10^{-7}$ throughout the whole simulation.

\begin{figure}
  \begin{center}
    \begin{overpic}[scale=0.8]{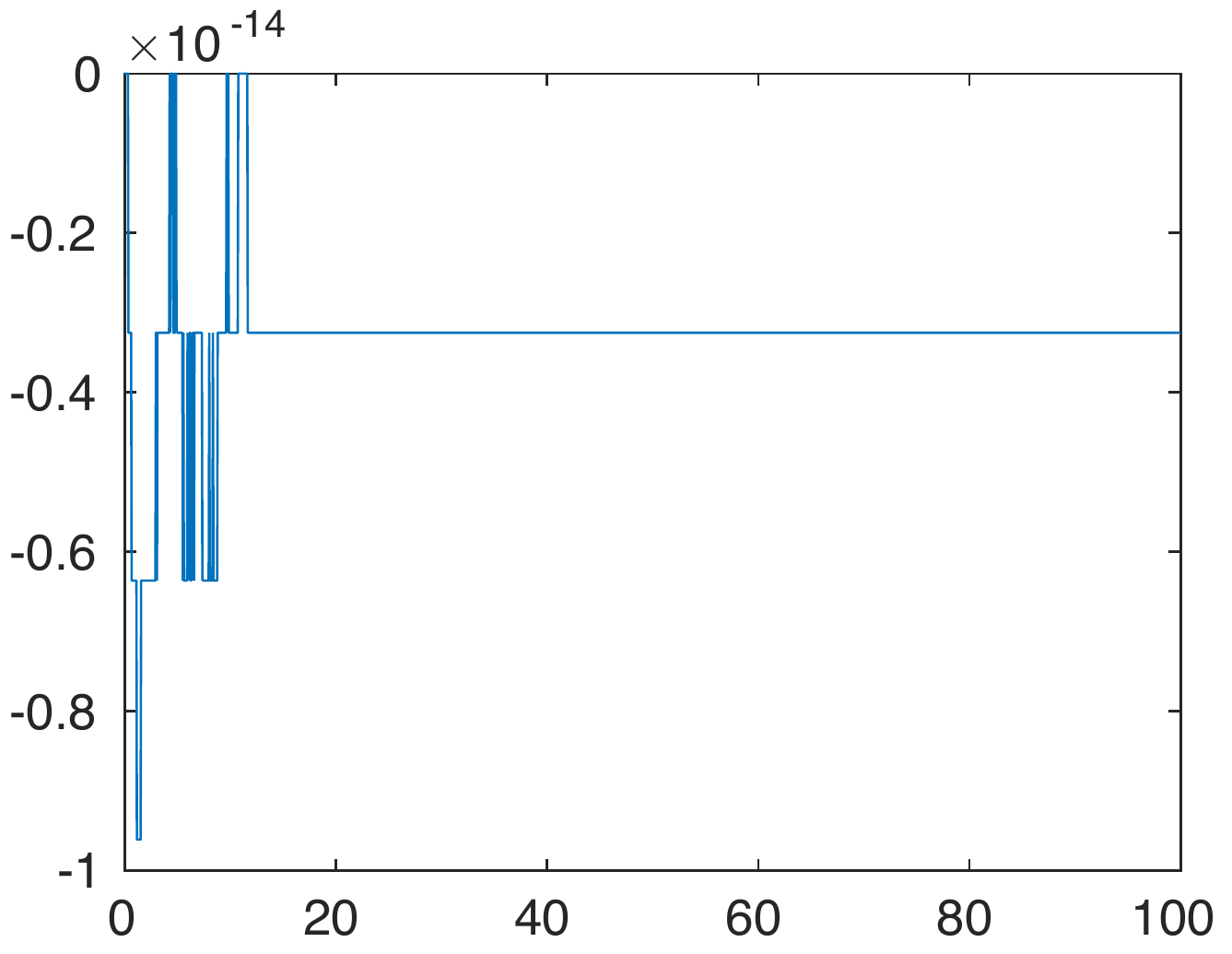}
      \put(48,-5){\FLARGE{\textbf{Time}}}
      \put(-5,15){\begin{sideways}
      \FLARGE{\textbf{Variation of total energy}}
      \end{sideways}}
    \end{overpic}
    \vspace{0.5cm}
    \caption{Landau damping test: relative variation of the total
      energy versus time for $\gamma^e=0.5$ and $\nu^e=1$. Penalty
      $\gamma^e$ is not applied to the equations of the first three
      Legendre modes.}
  \label{fig:LD:06}
  \end{center}
\end{figure}

Finally, the Legendre-Fourier method presented in this work provides
exact conservation laws.
The relative discrepancy of the mass, defined as
$(M(t^{\tau})-M(0))\slash{M(0)}$, and the discrepancy of momentum,
defined as $P(t^{\tau})-P(0)$, are \emph{exactly zero} at any discrete
time step $0\leq t^{\tau}\leq 200$ in our double precision
implementation and are therefore not shown.
The relative discrepancy of the total energy, defined as
$(\mathcal{E}_{tot}(t^{\tau})-\mathcal{E}_{tot}(0))\slash{\mathcal{E}_{tot}(0)}$
is shown in Figure~\ref{fig:LD:06} and is smaller than $10^{-14}$.

\subsection{Two-stream instability}

\begin{figure}
  \begin{center}
    \begin{overpic}[scale=0.5]{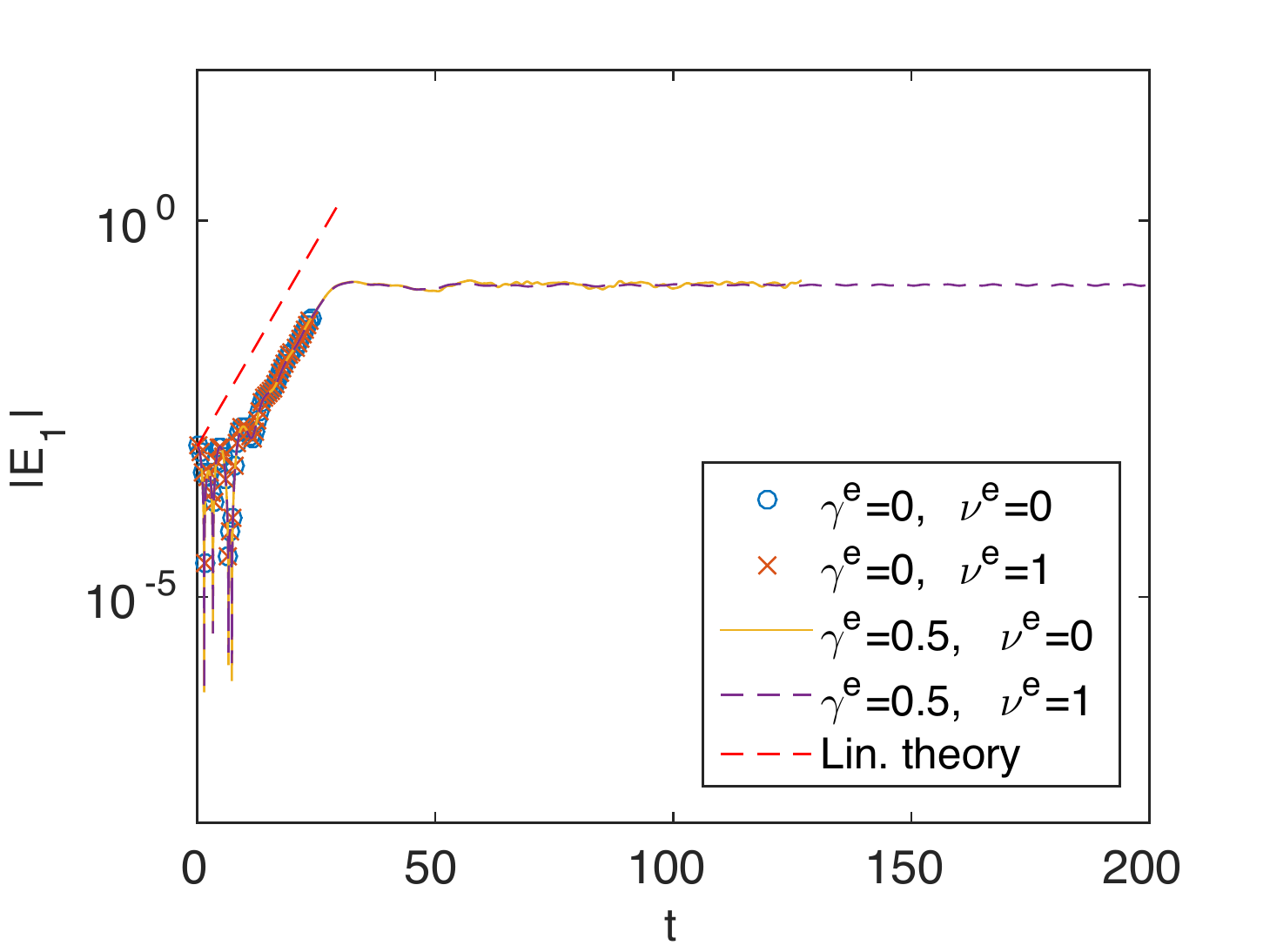}
      \put(48,-6){\textbf{Time}}
      \put(-5,36){\begin{sideways}$\mathbf{\abs{E_1}}$\end{sideways}}
    \end{overpic}
    \hspace{0.6cm}
    \begin{overpic}[scale=0.5]{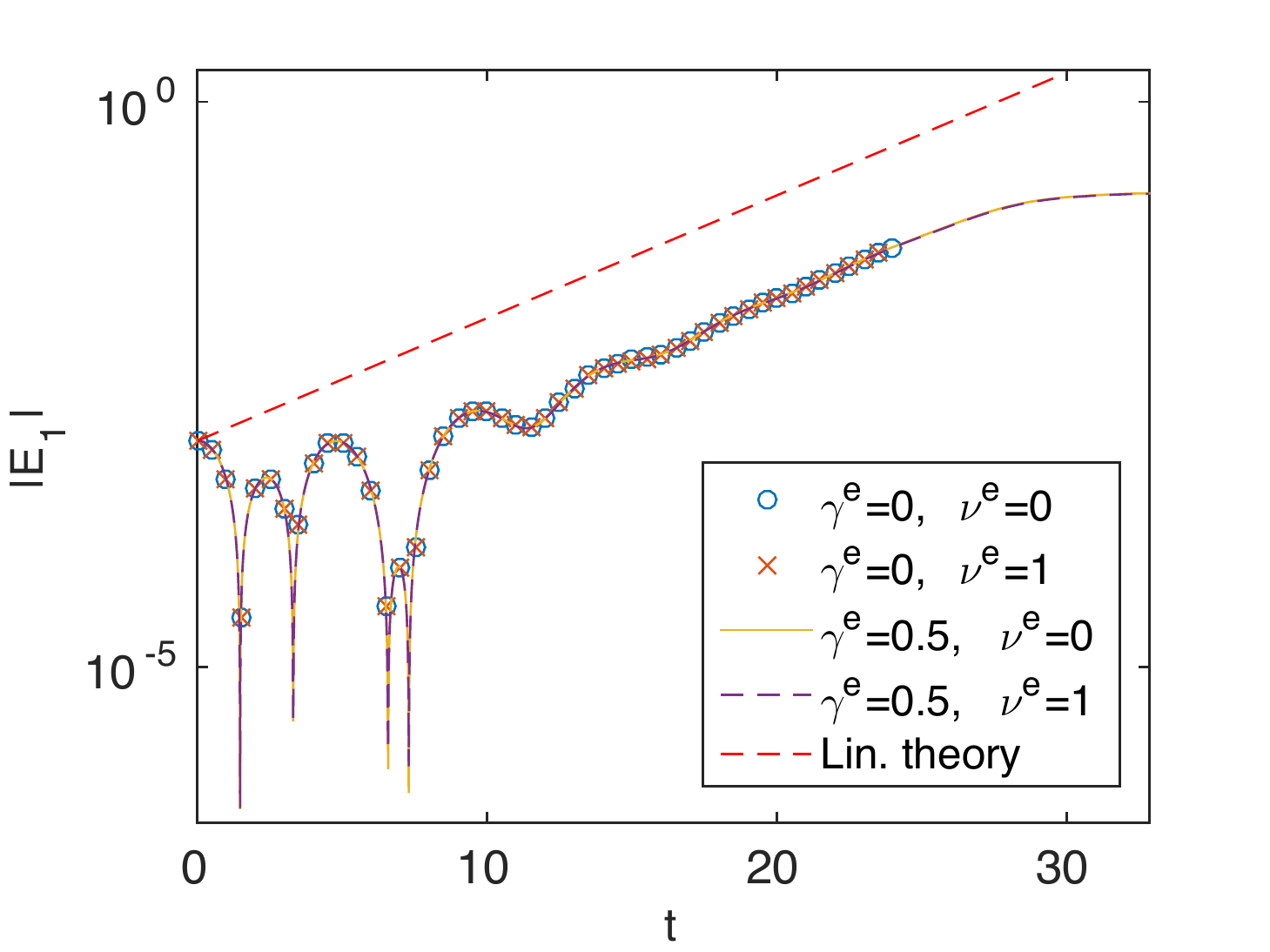}
      \put(48,-6){\textbf{Time}}
      \put(-5,36){\begin{sideways}$\mathbf{\abs{E_1}}$\end{sideways}}
    \end{overpic}
    \vspace{0.5cm}
    \caption{Two-stream instability test: first Fourier mode of the electric
      field versus time for the four combinations of
      $\gamma^e\in\{0,\,0.5\}$ and $\nu^e\in\{0,1\}$. Penalty
      $\gamma^e$ is not applied to the equations of the first three
      Legendre modes.}
  \label{fig:2S:00}
  \end{center}
\end{figure}

\begin{figure}
  \begin{center}
    \begin{overpic}[scale=0.5]{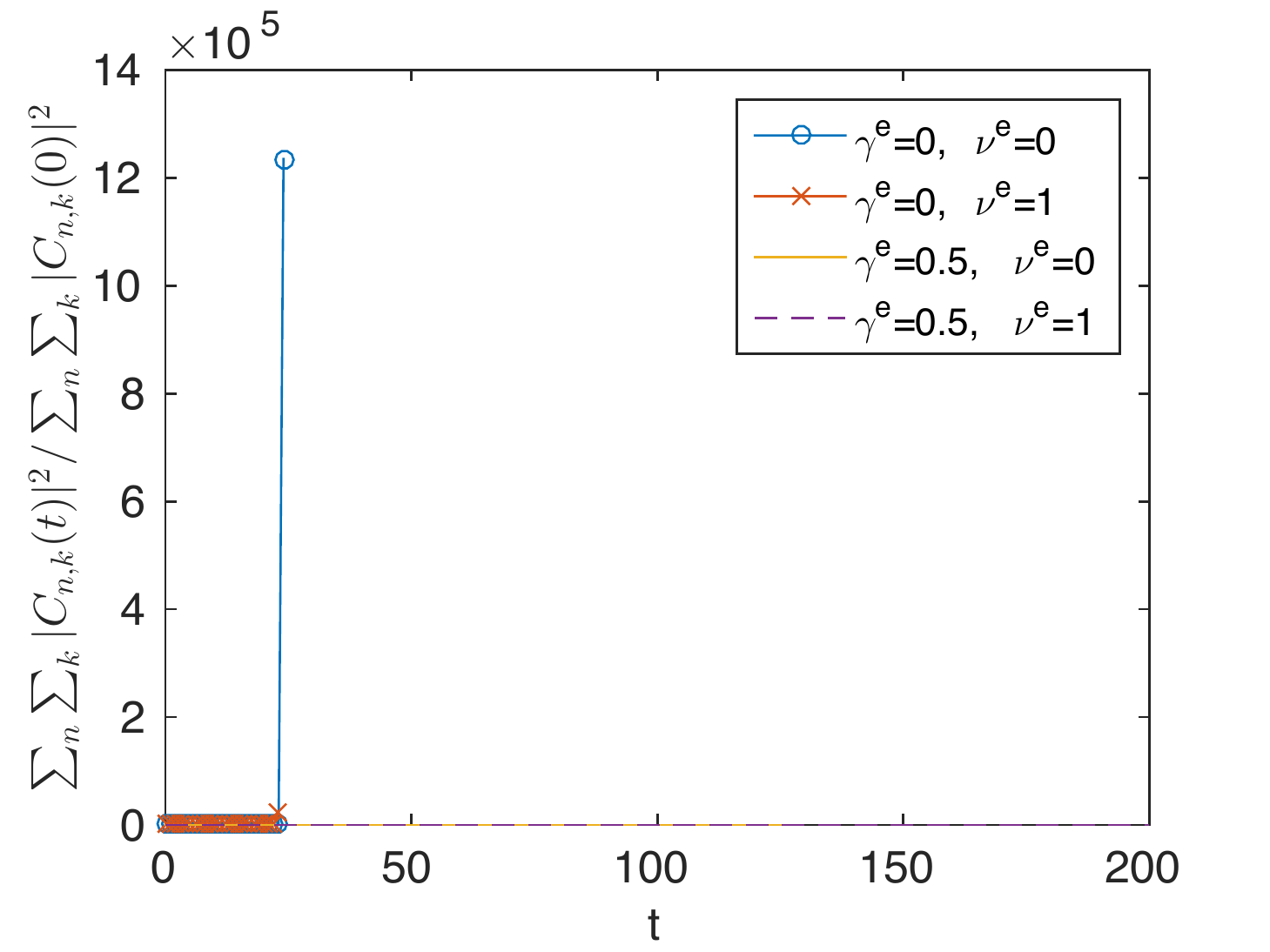}
      \put(48,-6){\textbf{Time}}
      \put(-8,21){\begin{sideways}\textbf{$L^2$ norm of $\fs$}\end{sideways}}
    \end{overpic}
    \hspace{0.75cm}
    \begin{overpic}[scale=0.5]{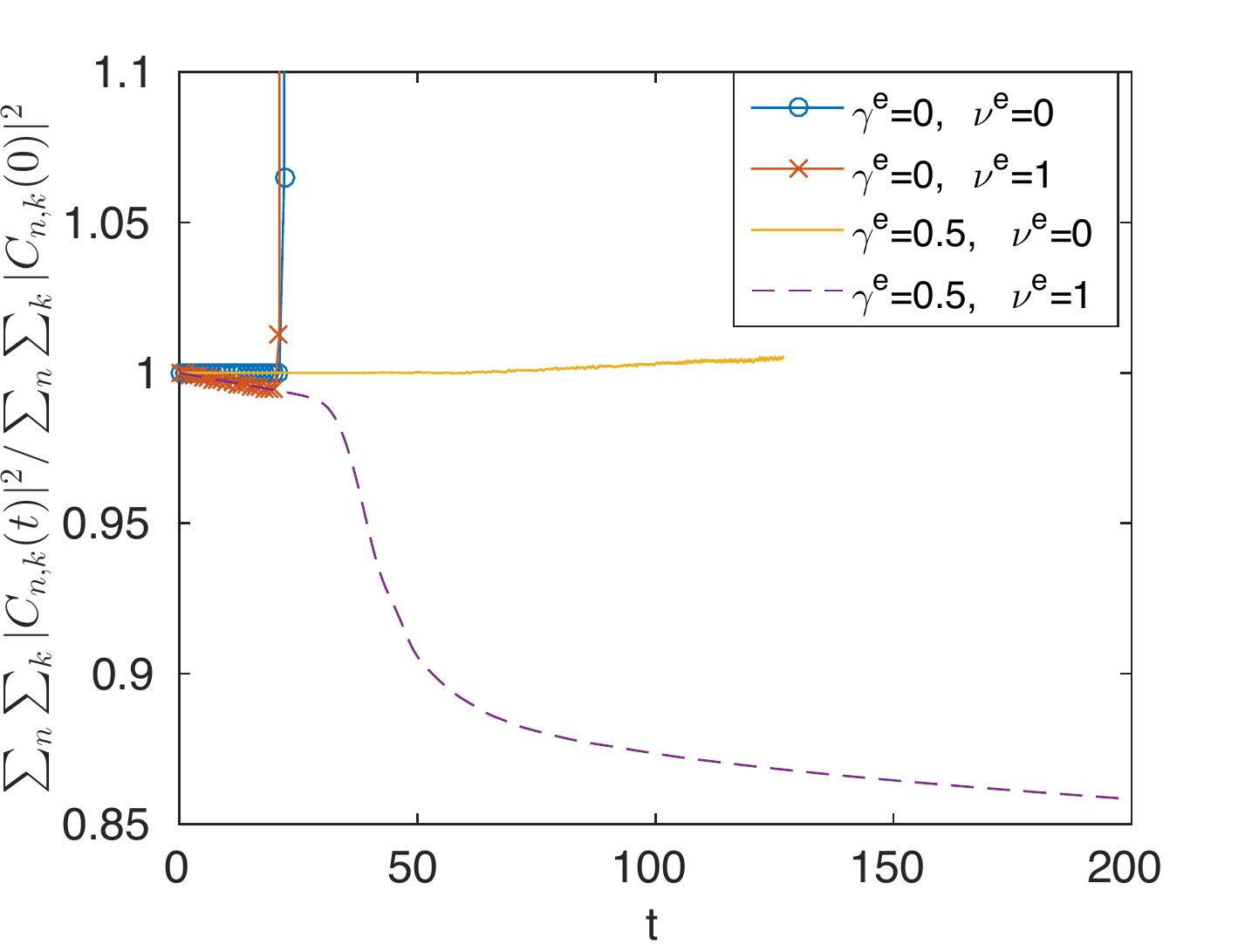}
      \put(48,-6){\textbf{Time}}
      \put(-8,21){\begin{sideways}\textbf{$L^2$ norm of $\fs$}\end{sideways}}
    \end{overpic}
    \vspace{0.5cm}
    \caption{Two-stream instability test: the left panel shows the $L2$ norm of 
      the electron distribution function $f^e$, cf.~\eqref{eq:relative:L2-norm:fs},
      versus time for the four combinations of
      $\gamma^e\in\{0,\,0.5\}$ and $\nu^e\in\{0,1\}$.
      The right panel is a zoom around $1$. Penalty $\gamma^e$ is not
      applied to the equations of the first three Legendre modes.}
  \label{fig:2S:02}
  \end{center}
\end{figure}

\begin{figure}
  \begin{center}
    \begin{overpic}[scale=0.8]{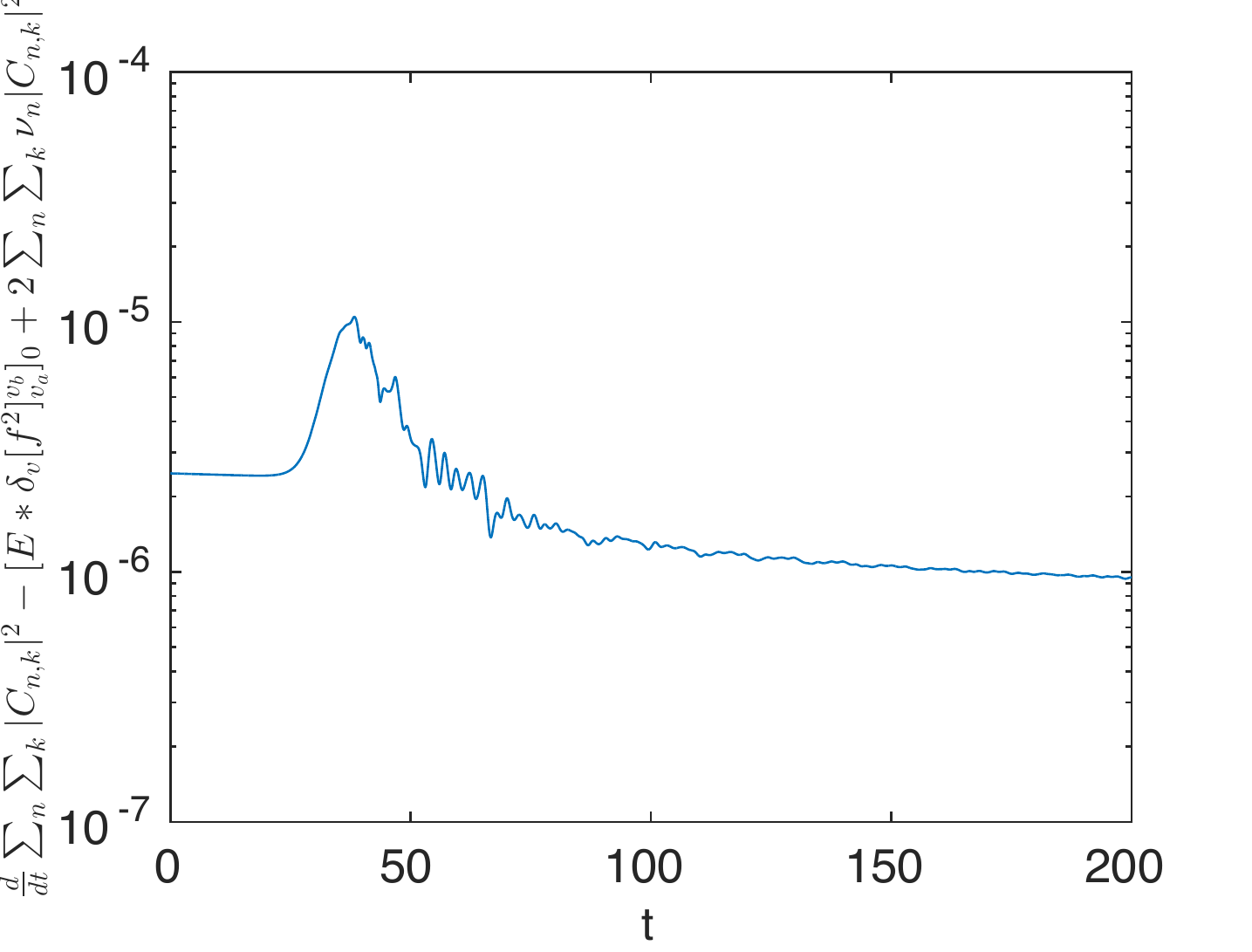}
      \put(48,-6){\FLARGE{\textbf{Time}}}
      \put(-8,20){\begin{sideways}
      \FLARGE{\textbf{Variation of $\norm{\fs}_{L^2(\Omega)}$}}
      \end{sideways}}
    \end{overpic}
    \vspace{0.75cm}
    \caption{Two-stream instability test: time variation of the $L2$
      norm of the electron distribution function $f^e$ as predicted by
      Theorem~\ref{theo:L2stab} versus time for $\gamma^e=0.5$ and
      $\nu^e=1$.
      The time derivative is approximated by central finite
      differences. Penalty $\gamma^e$ is not applied to the
      equations of the first three Legendre modes.}
  \label{fig:2S:03}
  \end{center}
\end{figure}

\begin{figure}
  \begin{center}
    \begin{overpic}[scale=0.8]{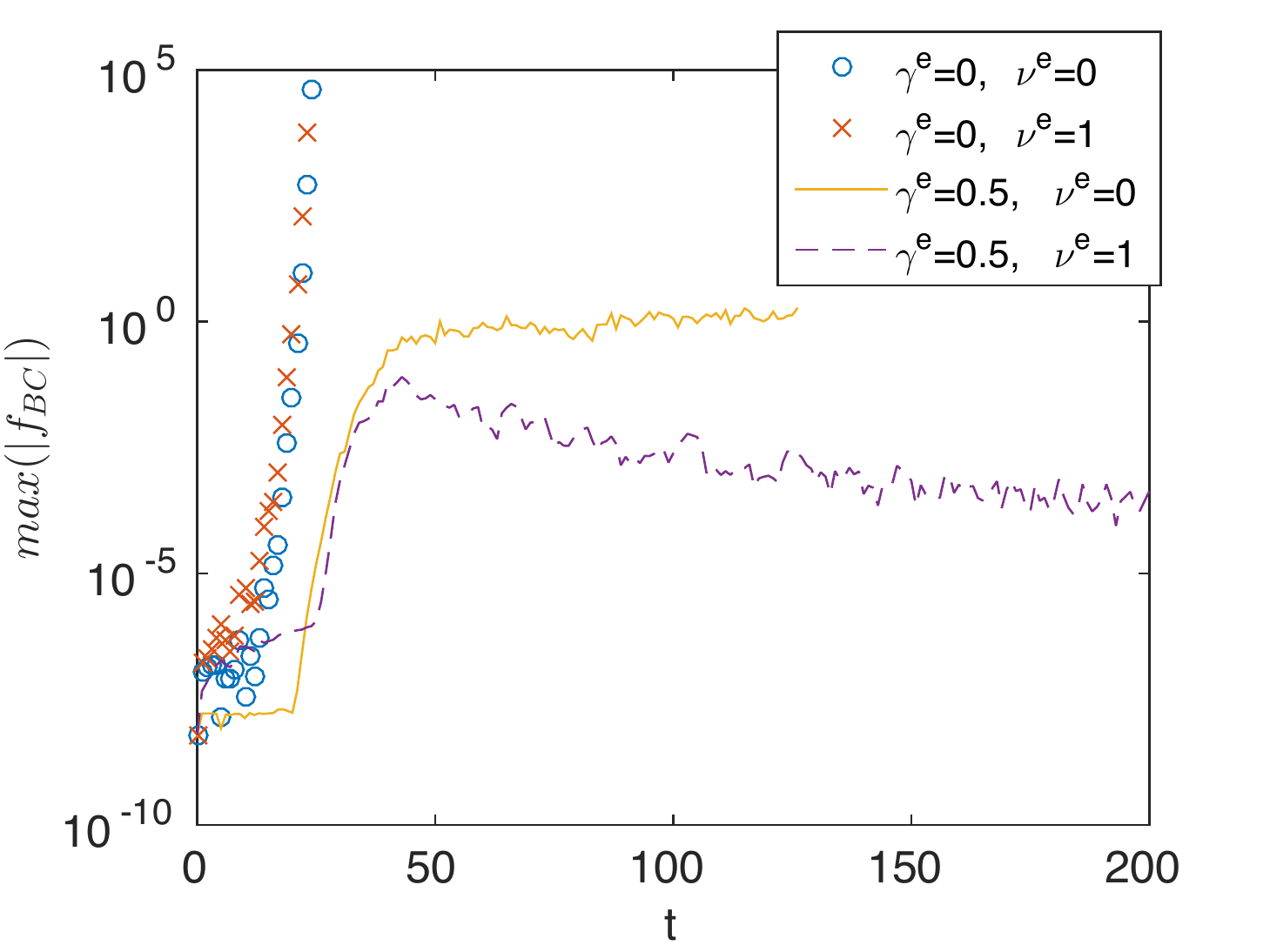}
      \put(48,-5){\FLARGE{\textbf{Time}}}
      \put(-5,29){\begin{sideways}
      \FLARGE{\textbf{max($\abs{\mathbf{f^e_{BC}}}$)}}
      \end{sideways}}
    \end{overpic}
    \vspace{0.75cm}
    \caption{Two-stream instability test: maximum value of the electron
      distribution function $f^e$ at the boundaries of the velocity
      range for the four combinations of $\gamma^e\in\{0,\,0.5\}$ and
      $\nu^e\in\{0,\,1\}$.  Penalty $\gamma^e$ is not applied to the
      equations of the first three Legendre modes.}
  \label{fig:2S:04}
  \end{center}
\end{figure}

\begin{figure}
  \begin{center}
    \begin{overpic}[scale=0.5]{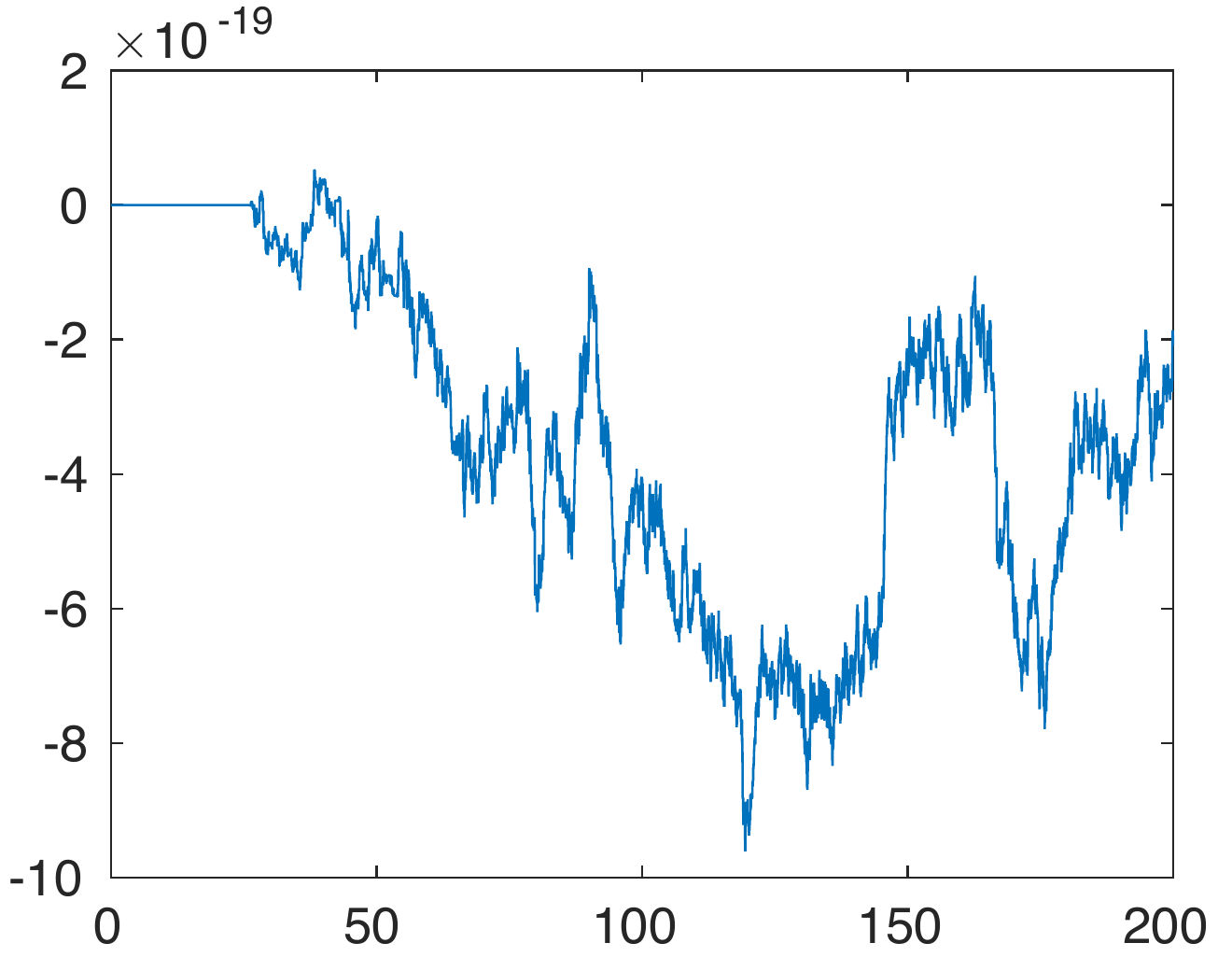}
      \put(48,-6){\textbf{Time}}
      \put(-8,0){\begin{sideways}\textbf{Variation of total momentum}\end{sideways}}
    \end{overpic}
    \hspace{0.75cm}
    \begin{overpic}[scale=0.5]{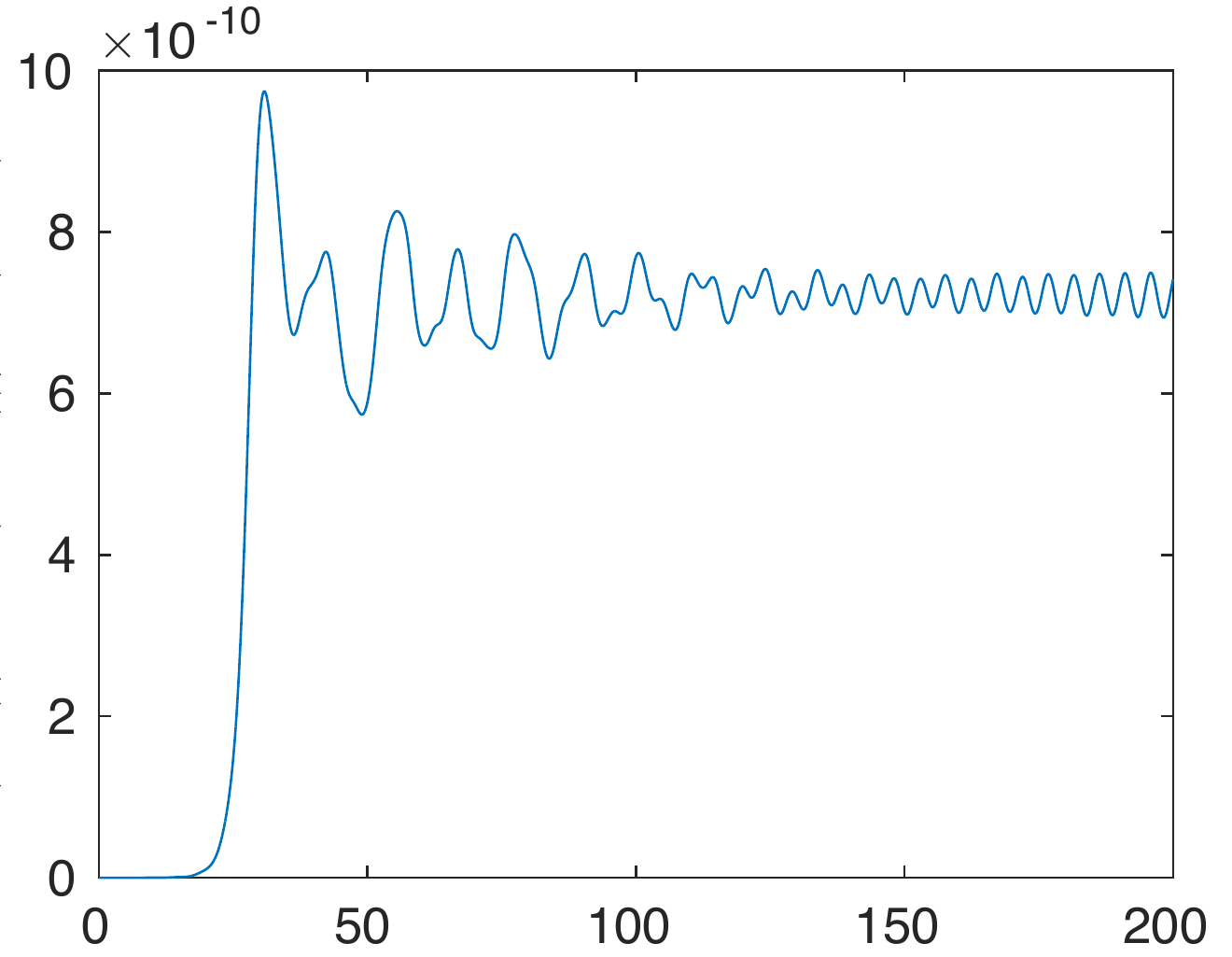}
      \put(48,-6){\textbf{Time}}
      \put(-8,6){\begin{sideways}\textbf{Variation of total energy}\end{sideways}}
    \end{overpic}
    \vspace{0.75cm}
    \caption{Two-stream instability test: the left panel shows the
      variation of momentum with respect to its initial value versus
      time, i.e., $P(t^{\tau})-P(0)$; the right panel shows the
      relative variation of energy with respect to its initial value
      versus time, i.e.,
      $(\mathcal{E}_{tot}(t^{\tau})-\mathcal{E}_{tot}(0))\slash{\mathcal{E}_{tot}(0)}$)
      for $\gamma^e=0.5$ and $\nu^e=1$. Penalty $\gamma^e$ is not
      applied to the equations of the first three Legendre modes. }
  \label{fig:2S:05}
  \end{center}
\end{figure}

\begin{figure}
  \begin{center}
    \begin{tabular}{cc}
    \includegraphics[scale=0.35]{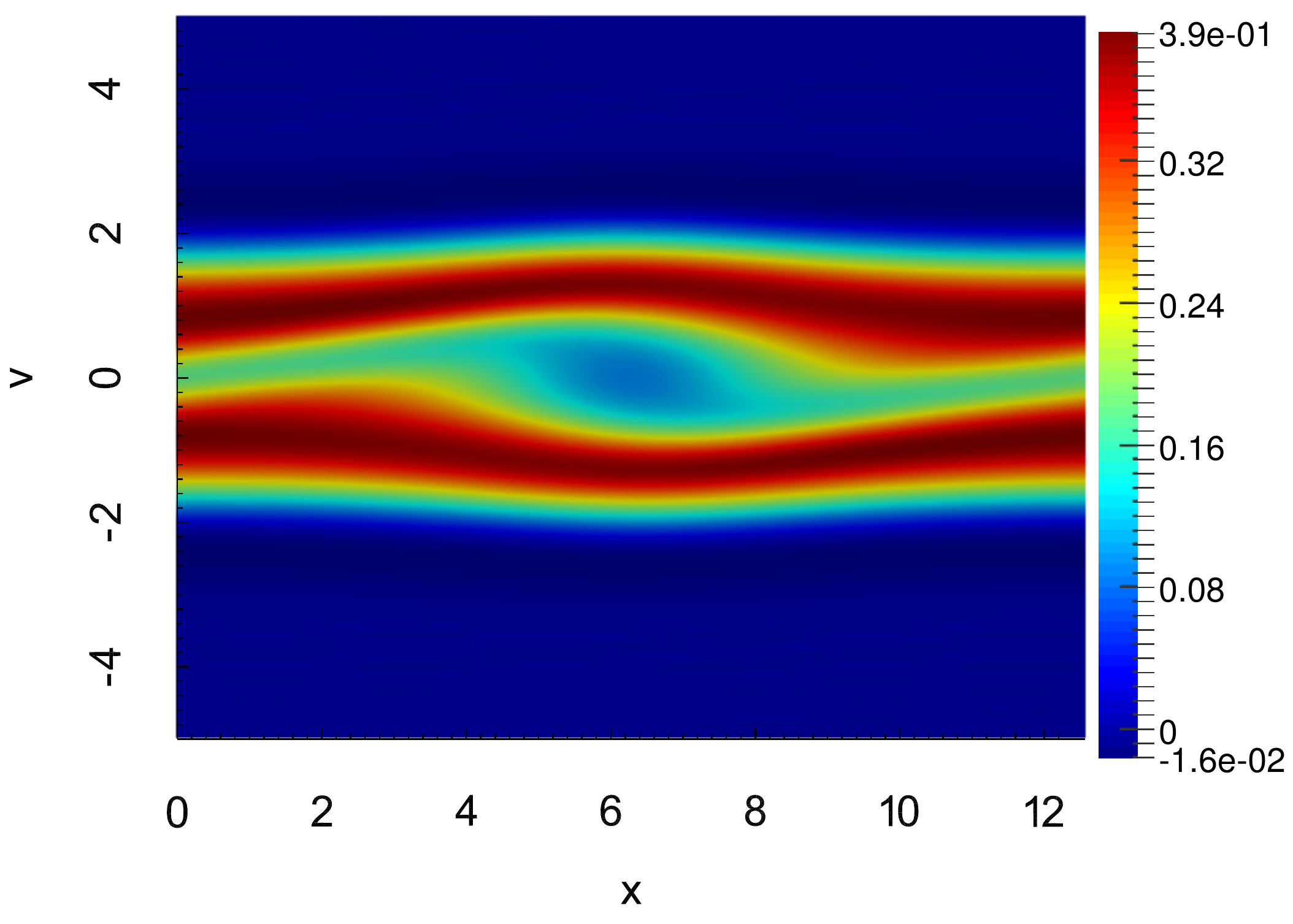}&
    \includegraphics[scale=0.35]{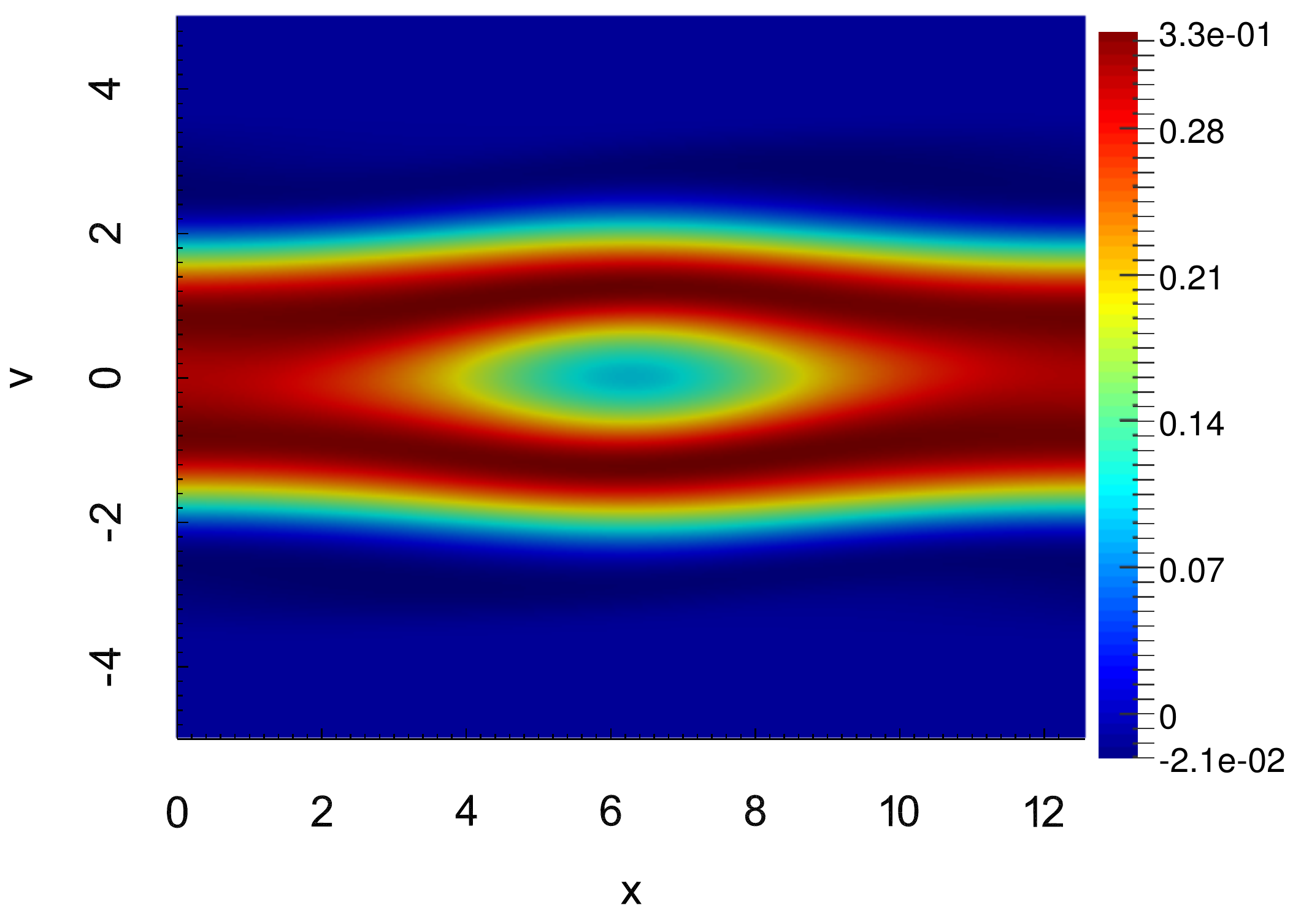}\\[-0.5em]
    \multicolumn{2}{c}{$\NL=50$, $[\va,\vb]=[-5,5]$}\\[1.5em]
    \includegraphics[scale=0.35]{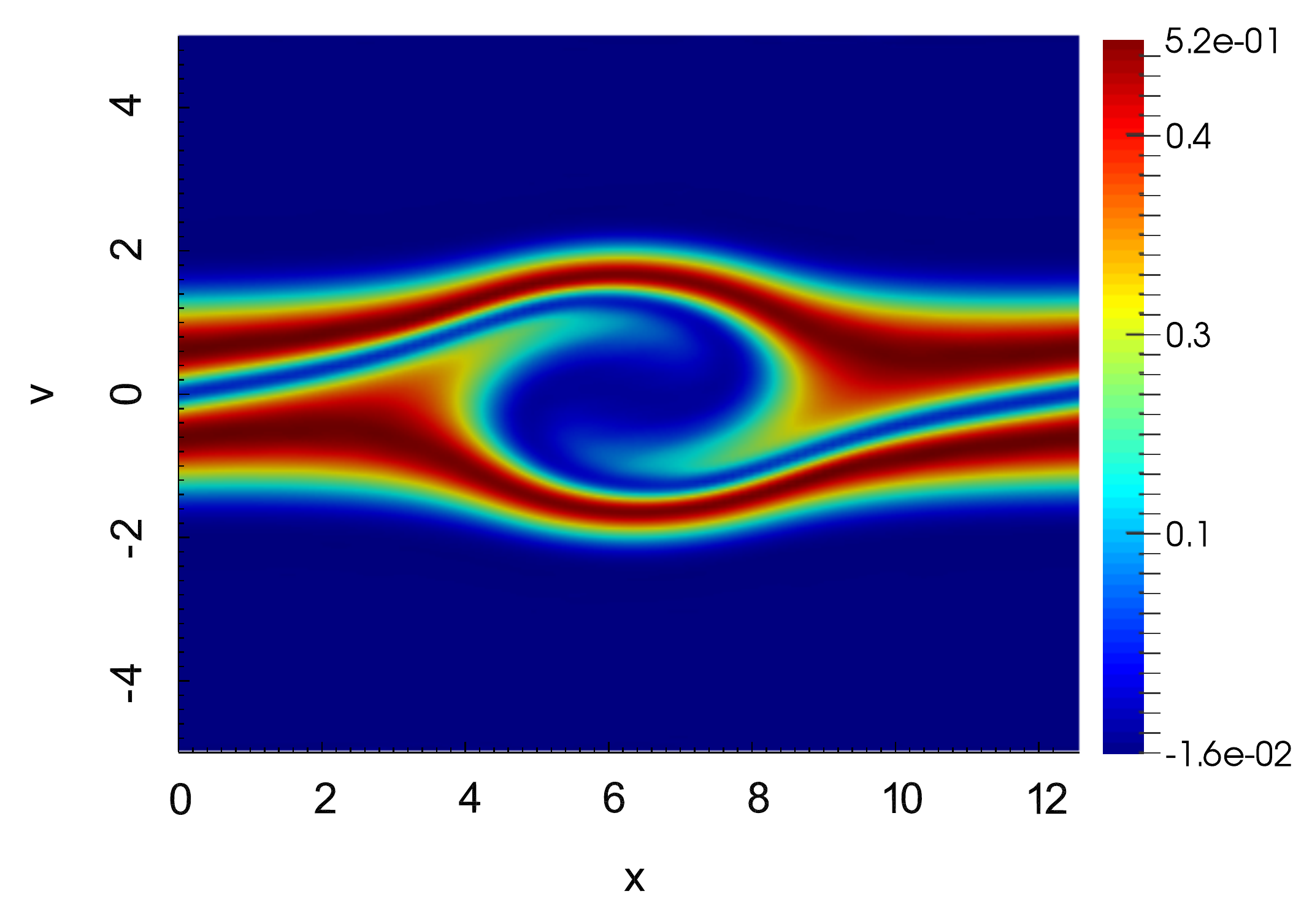}&
    \includegraphics[scale=0.35]{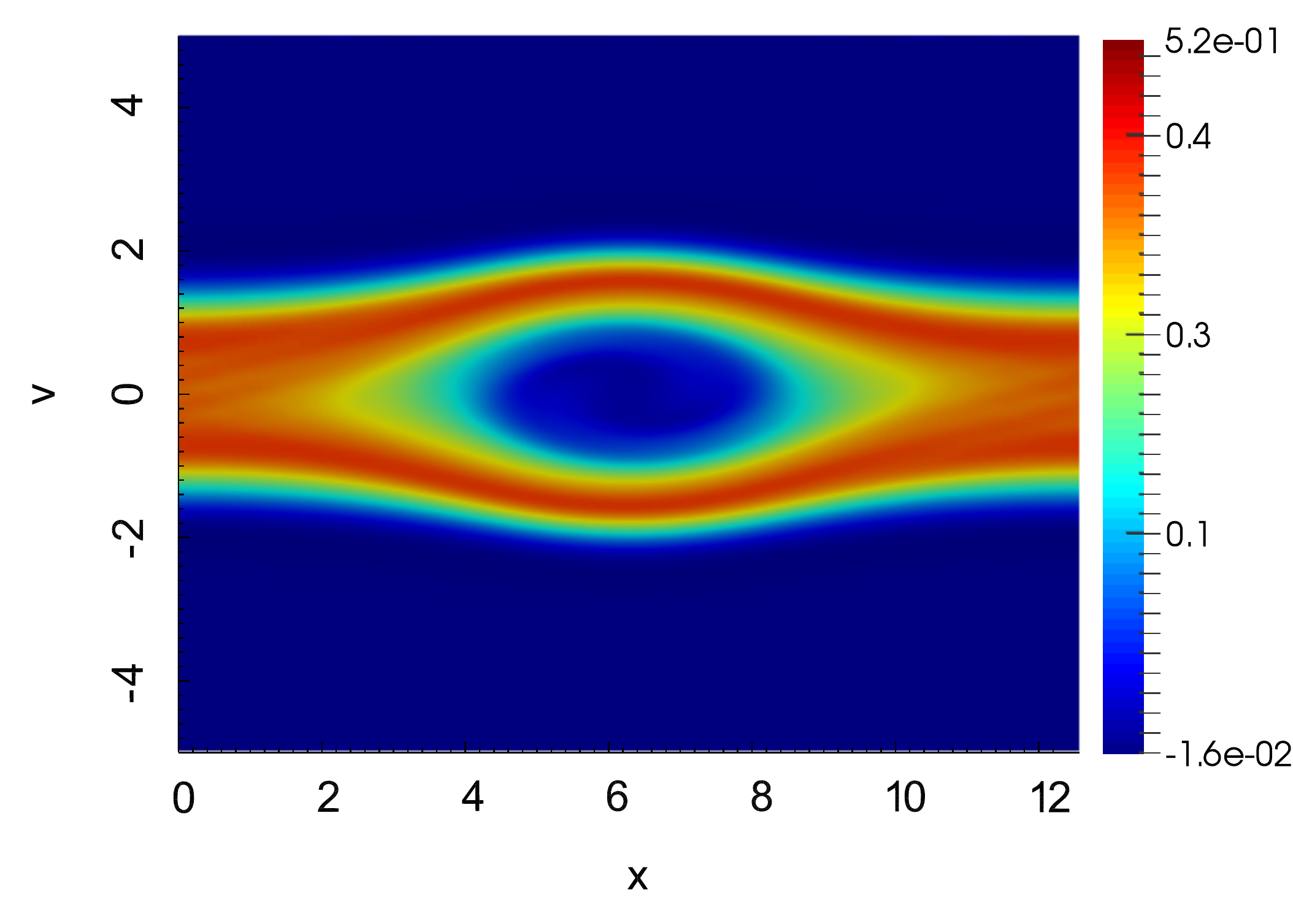}\\[-0.5em]
    \multicolumn{2}{c}{$\NL=100$, $[\va,\vb]=[-5,5]$}\\[1.5em]
    \includegraphics[scale=0.35]{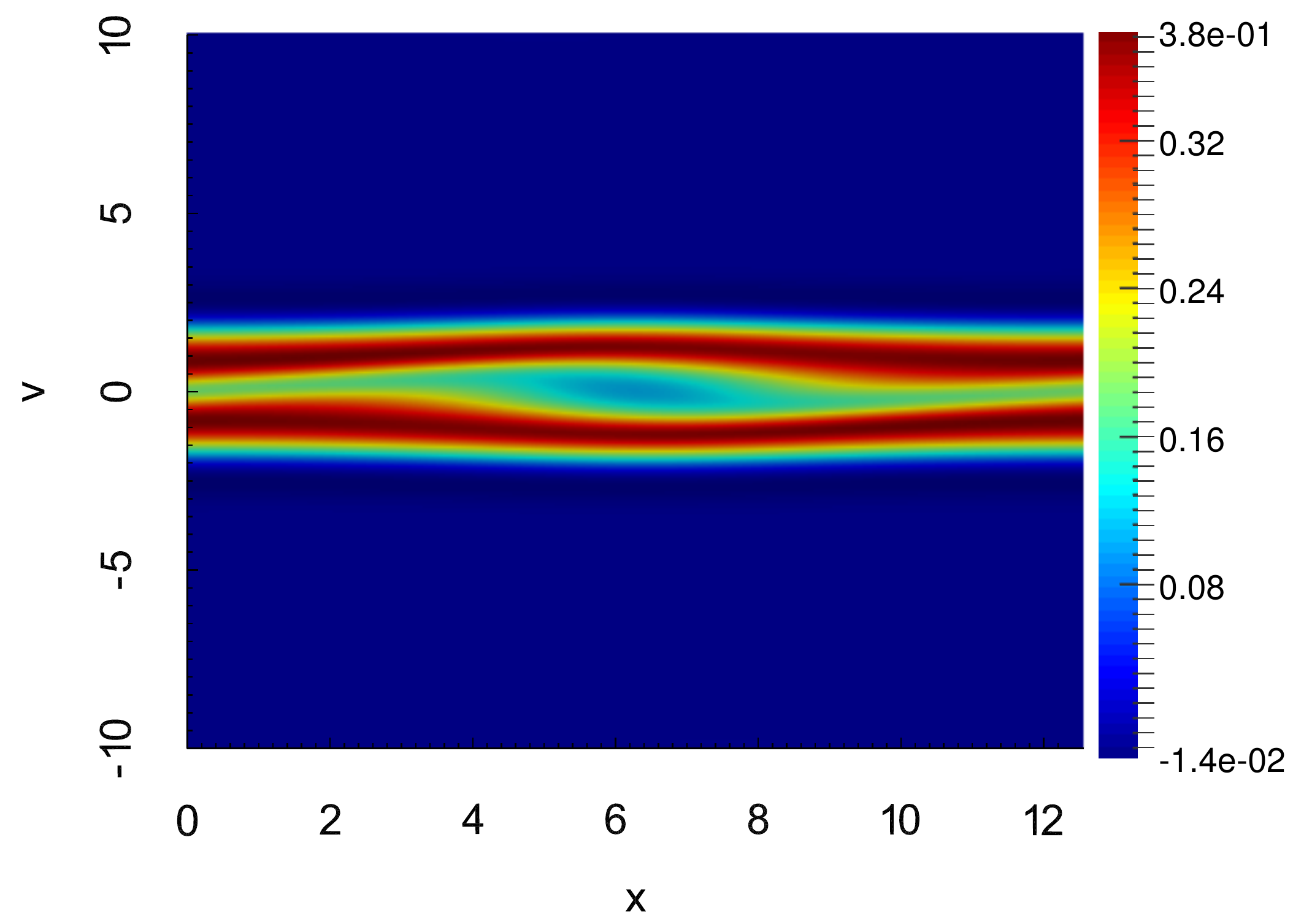}&
    \includegraphics[scale=0.35]{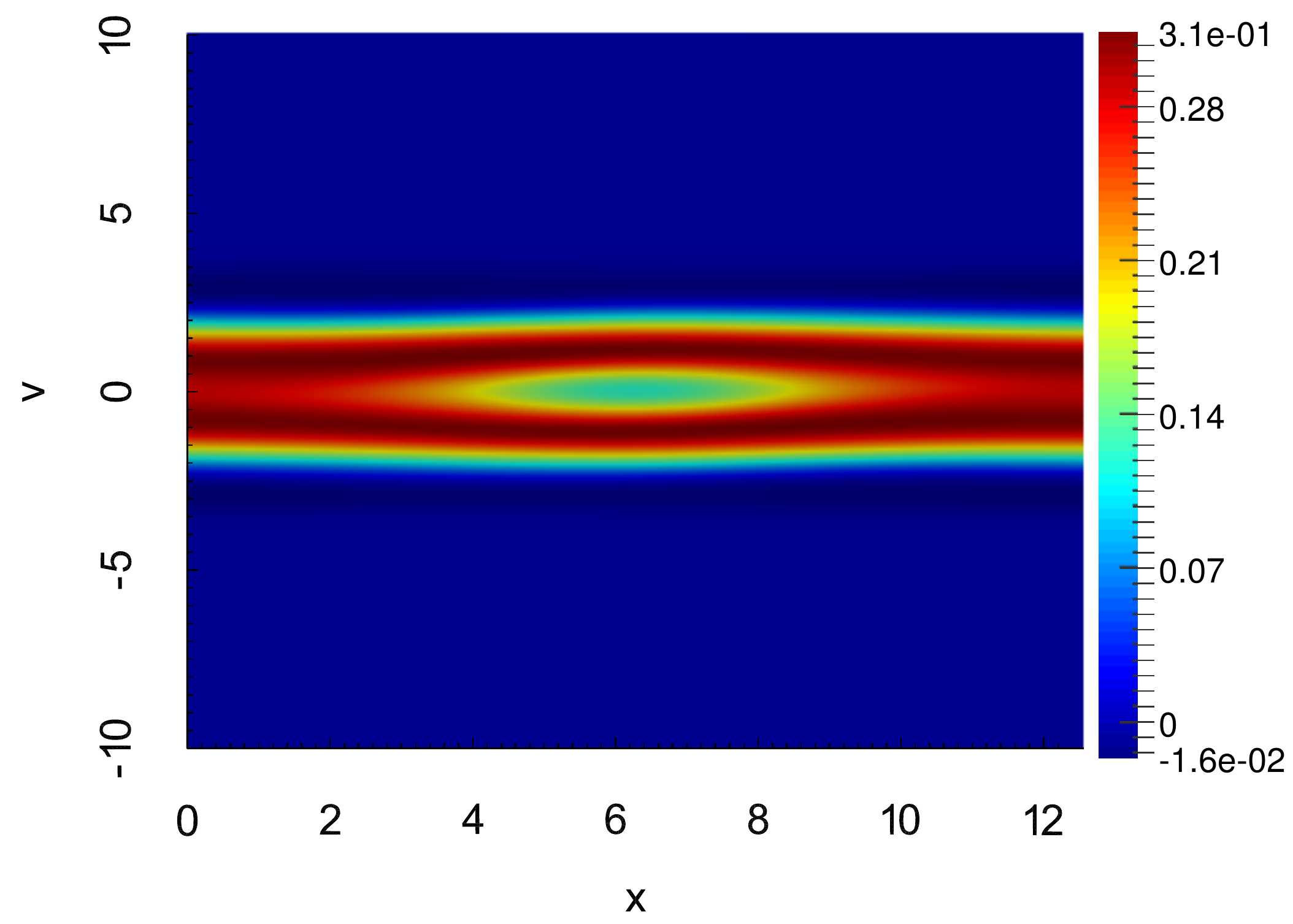}\\[-0.5em]
    \multicolumn{2}{c}{$\NL=100$, $[\va,\vb]=[-10,10]$}
    \end{tabular}
    \caption{Two stream instability test: the plots show the distribution function 
      $f^e$ in the phase space at the intermediate time $t=30$ (left panel) and the 
      final time $t=60$ (right panel).
      In these calculations we use $\nu^{e}=1$,
      $\gamma^{e}=0.5$, $\NF=25$, $L=4\pi$, and three different combinations
      of the number of Legendre modes $\NL$ and the velocity range $[\va,\vb]$, which are 
      displayed under the corresponding plots.
      The resolution of $f^e$ depends on such combinations.
      In particular, by comparing the plots on top and middle rows we see that
      the resolution of $f^e$ improves by increasing $\NL$ in a fixed velocity range.
      On the other hand, by comparing the plots on middle and bottom
      rows we see that the resolution of $f^e$ worsens by increasing
      the domain size with a fixed $\NL$.}
    \label{fig:Two-Stream:fs:phase-space}
  \end{center}
\end{figure}

The two-stream instability is excited when the distribution function
of a species consists of two populations of particles streaming in
opposite directions with a large enough relative drift velocity.
We initialize the electron distribution function with two
counter-streaming Maxwellians with equal temperature:
\begin{align}
  f^{e}(x,\vs,t=0) = 
  \frac{1}{2\sqrt{2\pi}\alpha}\left[ 
    e^{-\left(\frac{\vs-u_e}{\sqrt{2} \alpha}\right)^2} + 
    e^{-\left(\frac{\vs+u_e}{\sqrt{2} \alpha}\right)^2}
  \right]
  \,\left[1+\varepsilon\cos\Big(\frac{2\pi}{L}kx\Big)\right]
  \label{eq:two-stream:init-sol}
\end{align}
where $u_e$ is the drift velocity.
For this test case, we have chosen the following parameters:
$\alpha=1\slash{\sqrt{8}}$, $u_e=1$, $\varepsilon=10^{-3}$, $k=1$.
We integrate the Vlasov-Poisson system by using the time step 
$\Delta t=0.01$, $\NL=201$ Legendre modes, and $2\NF+1=51$ 
Fourier modes.
The domain of integration in phase space is set to $L=4\pi$,
$\vb=-\va=5$ for all the calculations shown in
Figures~\ref{fig:2S:00}-\ref{fig:2S:05},
while in Figure~\ref{fig:Two-Stream:fs:phase-space} we show the
distribution function of electrons that is computed for three
different combinations of $\NL$ and velocity range $[\va,\vb]$.
This example was also considered with similar input parameters as in
Ref.~in~\cite{camporeale1}, where the Vlasov equation was discretized
using $\NH = 100$ Hermite modes.
The electron distribution function was inizialized by combining two
drifting Maxwellians centered at two different velocities, each
expanded in the Hermite basis. Since the discretization was based on
the Asymmetrically Weighted Hermite basis
functions~\cite{holloway96,schumer}, the two Maxwellians were
completely described by setting only the first mode of each expansion.
The remaining modes were needed to describe the non-Maxwellian
evolution of the solution.
When using the Legendre-Fourier discretization proposed in this work,
there is no correspondence between the first mode and the Maxwellian
distribution.
Thus, in order to have sufficient accuracy, the spectral expansion
requires to consider all the polynomial modes from the beginning.

In Figure~\ref{fig:2S:00} we show the first Fourier mode of the
electric field $\abs{E_1(t)}$ versus time for the four combinations of
$\gamma^e\in\{0,\,0.5\}$ and $\nu^e\in\{0,\,1\}$.
The initial part of the dynamics is the same for the four curves and
one can see the development of the two-stream instability.
The slope of the numerical curves matches well the theoretical slope
predicted by the linear theory, which is shown as a dashed line in the
plot.
When $\gamma^e=0$, the two curves for $\nu^s\in\{0,\,1\}$ stop at
$t\sim 25$ (slightly prior to the end of the linear phase) because of
the development of a numerical instability.

When $\gamma^e=0.5$ and $\nu^e=1$ the scheme is numerically stable and
reaches the final time of the simulation, $T=200$, without problems.
Instead, the case $\nu^e=0$ stops converging at around $t\sim 130$
because of problems related to the behavior of $f^e$ at the boundary
(as documented below).

Figure~\ref{fig:2S:02} (left) shows the time evolution of the $L^2$
norm of the distribution function $f^e$ normalized with respect to
initial value according to~\eqref{eq:relative:L2-norm:fs} the cases
presented in Fig.~\ref{fig:2S:00}.
Figure~\ref{fig:2S:02} (right) shows a zoom around $1$. 
One can clearly see that the $L^2$ norm of $f^e$ grows unboundedly
when $\gamma^e=0$, indicating that the first term on the right hand
side of equation~\eqref{eq:L2stab} provides a positive feedback that
is not even compensated by the collisional term when $\nu^e=1$.
Hence, the scheme is numerically unstable.

%
  When $\gamma^e=0.5$ the scheme is numerically stable.
  Indeed, by applying $\gamma^e$ to all the moment equation, we have
  verified numerically that the $L^2$ norm of $f^e$ is constant in
  time for $\nu^e=0$ and damps for $\nu^e=1$ as predicted by Theorem
  $3.2$, cf. equation~\eqref{eq:stability:time:variation}.
  If $\gamma^e$ is applied to all the Legendre modes \emph{except the
    first three} we obtain the behavior shown in Figure~\ref{fig:2S:02},
  where a slow growth of the $L^2$ norm of $f^e$ is visible for
  $\nu^e=0$.
  Figure~\ref{fig:2S:03} shows the numerical representation of
  equation~\eqref{eq:L2stab}, where the time derivative is approximated
  by central finite differences for the case $\gamma^e=0.5$ and
  $\nu^e=1$.
  From this figure, we deduce that Theorem~\ref{theo:L2stab} and
  equation~\eqref{eq:L2stab}, are verified numerically to a good degree
  of accuracy.

The behavior of the maximum value of the distribution function on the
domain boundary, $\abs{f_{BC}}$, is shown in Fig.~\ref{fig:2S:04}.
As expected, for $\gamma^e=0$ the simulation is unstable and $f^e$
grows unbounded on the boundary.
The stabilization provided by $\gamma^e=0.5$ is effective and limits
the value of $f$ there.
However, when $\nu^e=0$ one can see that $\abs{f_{BC}}$ still grows
sizably and becomes of order unity (i.e. of the same order of the
initial distribution function) at around $t\approx 50$.
Clearly this signals that the simulation is not accurate anymore. 
When $\nu^e=1$, on the other hand, $\abs{f_{BC}}$ remains reasonably
small throughout the simulation.

\begin{figure}
 \begin{center}
   \begin{overpic}[scale=0.5]{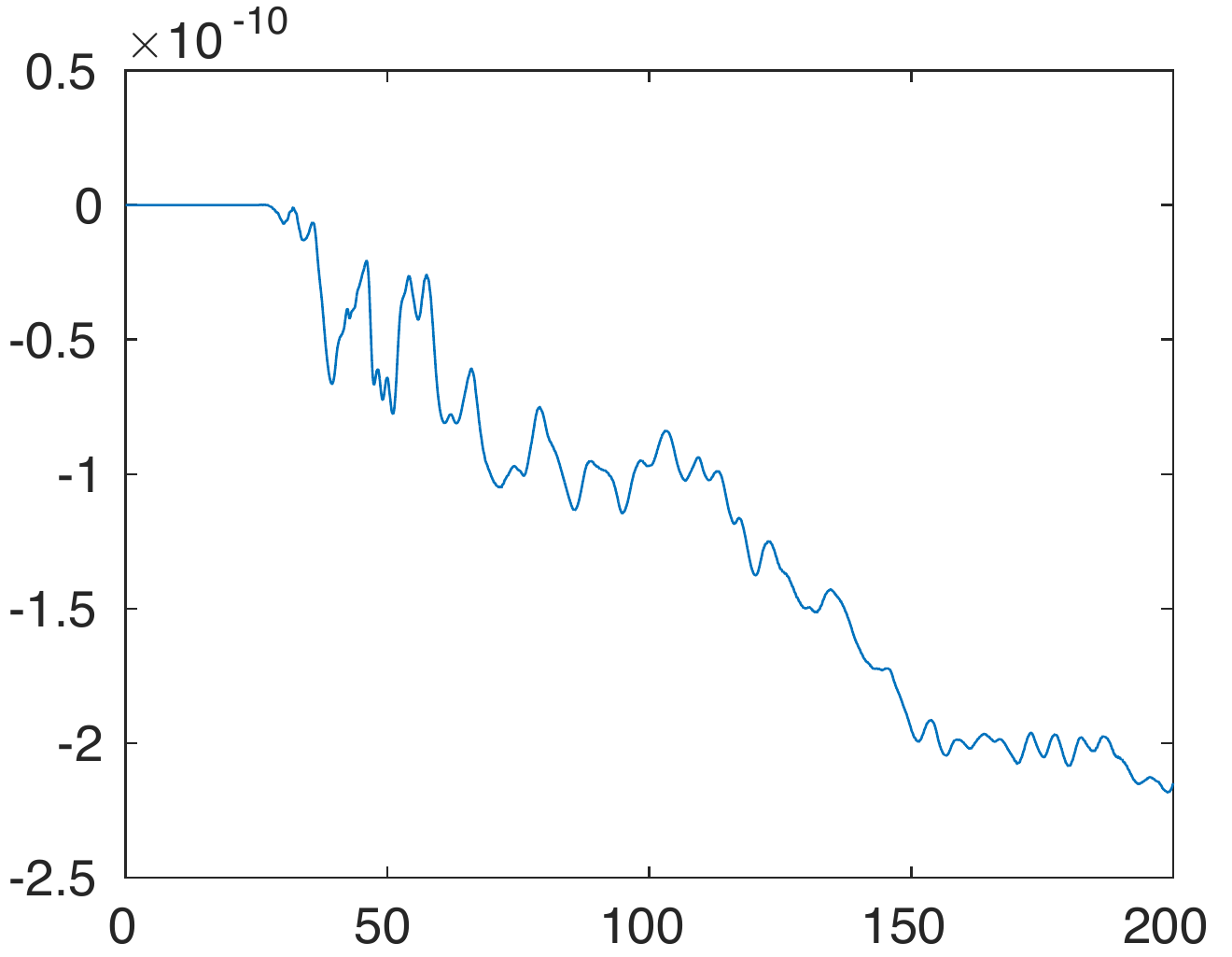}
      \put(48,-6){\textbf{Time}}
      \put(-8,0){\begin{sideways}\textbf{Variation of total momentum}\end{sideways}}
   \end{overpic}
   \hspace{0.75cm}
   \begin{overpic}[scale=0.5]{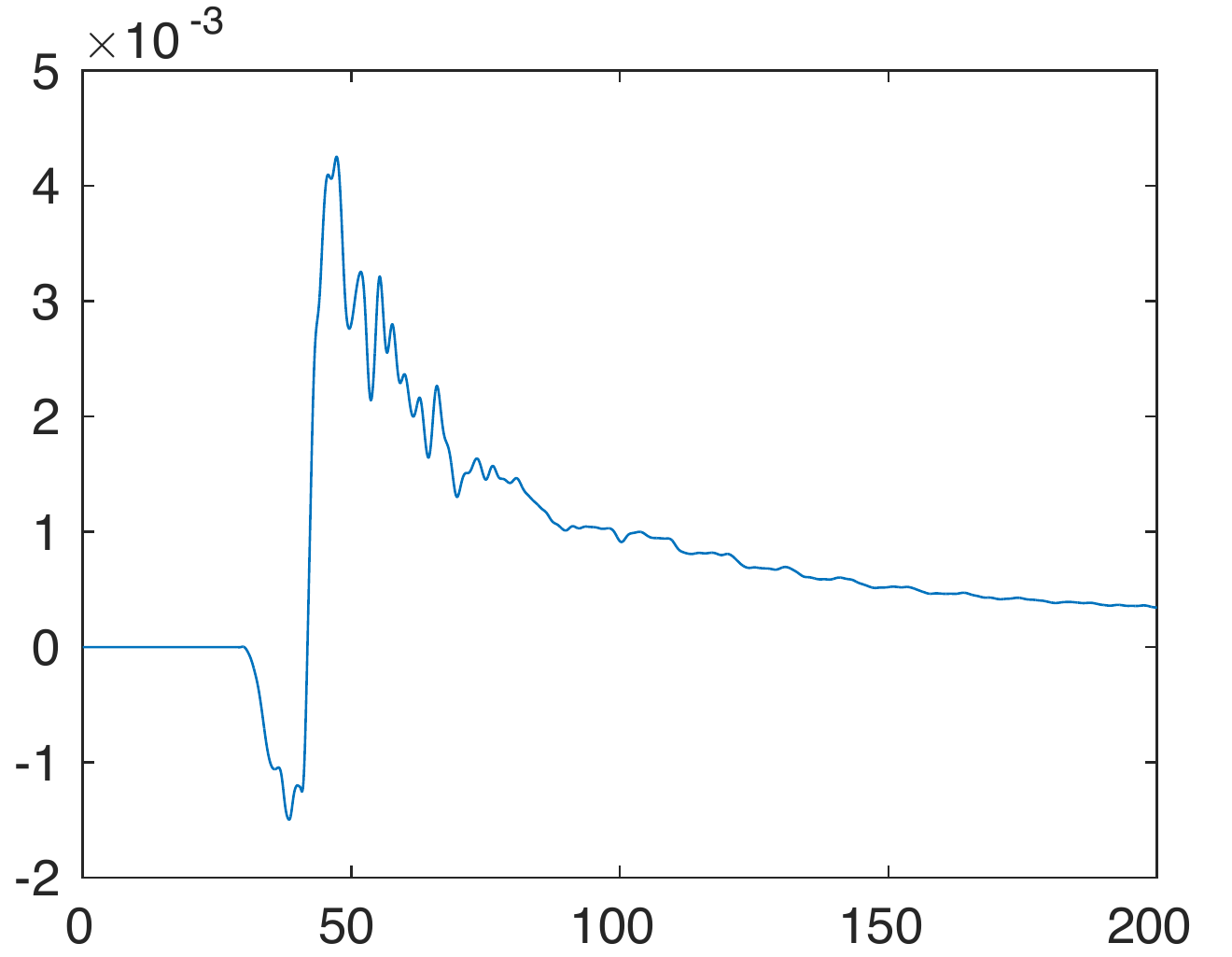}
      \put(48,-6){\textbf{Time}}
      \put(-8,6){\begin{sideways}\textbf{Variation of total energy}\end{sideways}}
   \end{overpic}
   \vspace{0.75cm}
   \caption{Two-stream instability: variation of total energy versus
     time by applying $\gamma^e$ to all modes; the collisional
     coefficient is $\nu=1$. Penalty $\gamma^e$ is applied to the
     equations of all Legendre modes.  }
   \label{fig:2S:09-all}
 \end{center}
\end{figure}

In Figures~\ref{fig:2S:05} and~\ref{fig:2S:09-all} we show the
variation in time of momentum (left plot, $P(t^{\tau})-P(0)$ although
$P(0)=0$ in this case) and relative variation in time of total energy
(right plot,
$(\mathcal{E}_{tot}(t^{\tau})-\mathcal{E}_{tot}(0))\slash{\mathcal{E}_{tot}(0)}$)
with respect to the initial value.
As for all the previous figures, the plots shown in
Figure~\ref{fig:2S:05} are obtained by applying the penalty $\gamma^e$
to all the moment equations except the first three.
In this case, total momentum and total energy, as well as mass which
is not shown, are conserved extremely well in the simulations, as
predicted by the analysis of Sections~\ref{sec:conservation:laws}
and~\ref{sec:time:integration}.
Instead, the results of Figure~\ref{fig:2S:09-all} are obtained by
applying penalty $\gamma^e$ to all the moment equations.
In this case, the total momentum variation that is visible is
of the order of magnitude of $10^{-10}$ and total energy variation is
of the order of magnitude of $10^{-3}$.
These results are still in accord with the analysis of
Sections~\ref{sec:conservation:laws} because we know from sections
\ref{subsect:momentum:conservation}
and~\ref{subsect:energy:conservation} that both momentum and energy
variation contain boundary terms that are not included in this
diagnostics.
It is worth noting that these boundary terms explicitly contain
$\gamma^e$, and are zero if $\gamma^e=0$ in their expression.
Also note that in the two-stream instability problem, momentum is
symmetric and that these results show that the symmetry of the problem
is not violated by the Legendre-Fourier method.

In Figure \ref{fig:Two-Stream:fs:phase-space} we show the electron
distribution function in phase space that is computed by using three
different combinations of $\NL$, the number of Legendre modes, and
velocity range $[\va,\vb]$ for $\gamma^s=0.5$ and $\nu^s=1$.
In particular, the plots on top are obtained by using $\NL=50$ and
integrating over the velocity range $[-5,5]$; the plots in the middle
are obtained by using $\NL=100$ and the velocity range $[-5,5]$; the
plots on bottom are obtained by using $\NL=100$ and the velocity range
$[-10,10]$.
The plots on the left show the distribution function at $t=30$, the
plots on the right at $t=60$.
The resolution of $f^e$ clearly depends on the combination that is
chosen: it improves by increasing $\NL$ in a fixed velocity range and
it worsen by increasing the domain size with a fixed $\NL$.

\subsection{Ion acoustic wave}

Last, we consider the evolution of an ion acoustic wave.
This is a truly multiscale example, occurring on the slow time scales
associated with the ions but where the electron motion concurs in
defining the properties of the wave.
Following~\cite{camporeale1}, we initialize a perturbation in the ion
distribution function at t=0
\begin{equation}
  f^i(x,v,t=0)=\frac{1}{\sqrt{2\pi}\alpha_i} 
  e^{-\left(\frac{v}{\sqrt{2} \alpha_i}\right)^2}
  \,\left[1+\varepsilon\cos\Big(\frac{2\pi}{L}x\Big)\right]
\end{equation}
while the electrons are Maxwellian and unperturbed
\begin{equation}
  f^e(x,v,t=0)=\frac{1}{\sqrt{2\pi}\alpha_e} 
  e^{-\left(\frac{v}{\sqrt{2} \alpha_e}\right)^2}
\end{equation}
Other parameters are $L=10$, $N_L=101$, $2N_F+1=51$,
$\varepsilon=0.01$, $\alpha_i=1/135$, $\alpha_e=1$,
$\gamma^e=0.5$, while $\Delta t$ and $\nu^s$ are varied
parametrically.
Although we only present results with a smaller perturbation
$\varepsilon=0.01$, we have also tried larger perturbations and
essentially successfully reproduced the results of
Ref.~\cite{camporeale1} for $\varepsilon=0.2$.

\begin{figure}
  \begin{center}
    \begin{overpic}[scale=0.8]{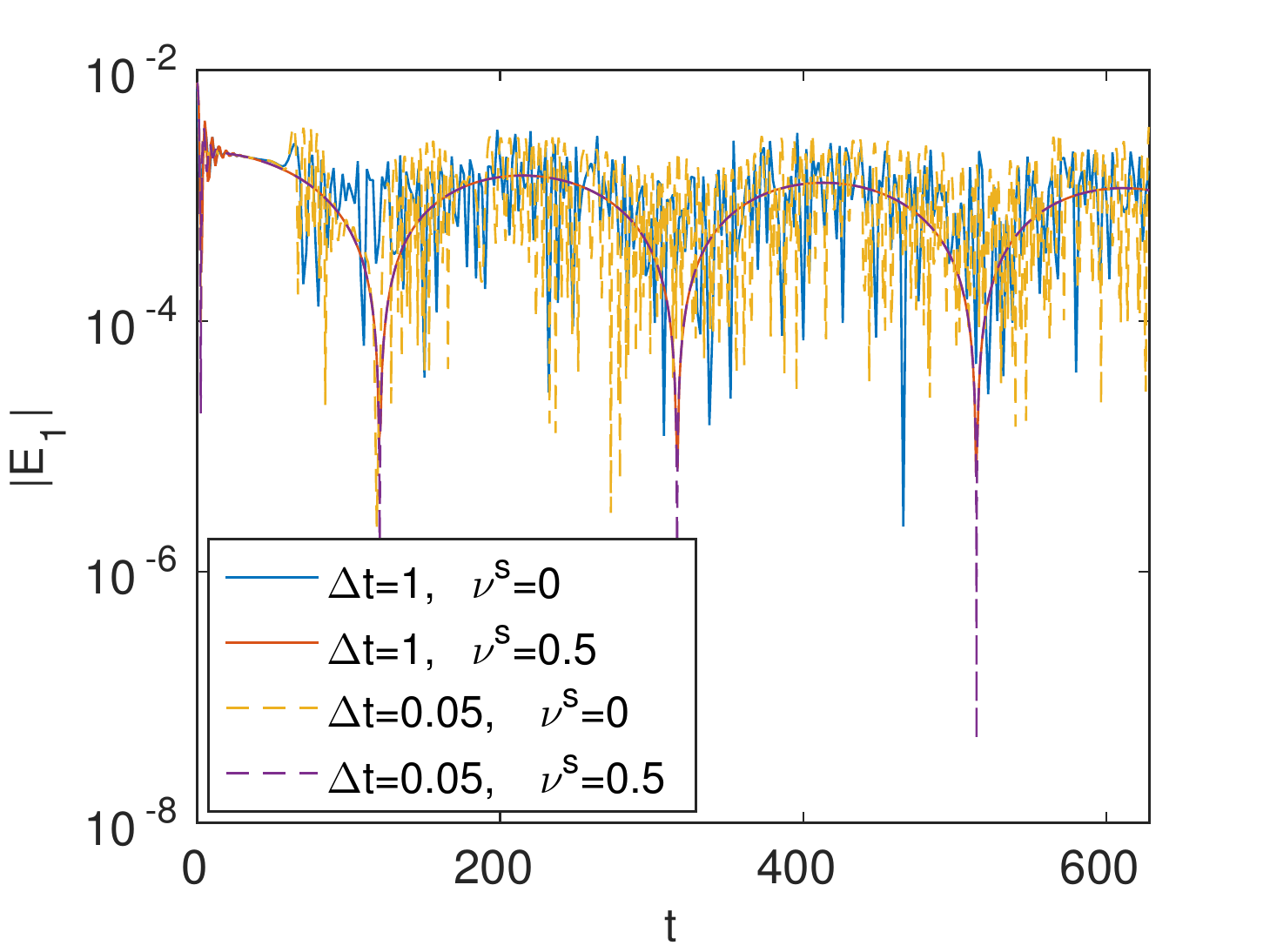}
      \put(48,-4){\begin{large}\textbf{Time}\end{large}}
      \put(-5,36){\begin{sideways}\begin{large}$\mathbf{\abs{E_1}}$\end{large}\end{sideways}}
    \end{overpic}
    \vspace{0.5cm}
    \caption{Ion acoustic wave test: first Fourier mode of the electric
      field versus time for the four combinations of
      $\gamma^e\in\{0,\,0.5\}$ and $\nu^e\in\{0,1\}$. Penalty
      $\gamma^e$ is not applied to the equations of the first three
      Legendre modes.}
  \label{fig:IA:00}
  \end{center}
\end{figure}

Figure \ref{fig:IA:00} shows the amplitude of the electric field for
the first Fourier mode initially excited at $t=0$.
Four curves are plotted, corresponding to $\Delta t\in\{0.05,\,1\}$
and $\nu^s\in\{0,\,0.5\}$.
The initial evolution of the system is the same for all the curves and
one can see some electron oscillations.
However, when $\nu^s=0$ the simulations are corrupted by a large
amount of unphysical oscillations (quite irrespective of $\Delta t$).
When $\nu^s=0.5$, on the other hand, the ion acoustic wave signal is
recovered well: the period of $|E_1|$ obtained from the simulations is
$197$, in good agreement with the theoretical value of $201$.
We note that the curves obtained with $\nu^s=0.5$ and $\Delta t=0.05$
and $\Delta t=1$ are virtually indistinguishable, showing the ability
of our numerical scheme to step over the faster frequency in the
system, the electron plasma frequency, without any sign of numerical
instability.
We have also performed simulations with larger $\Delta t$ (up
  to $\Delta t=10$, not shown).
  The ion acoustic wave becomes progressively less accurate but, as
  expected, there is no sign of numerical instabilities.

\begin{figure}
  \begin{center}
    \begin{overpic}[scale=0.8]{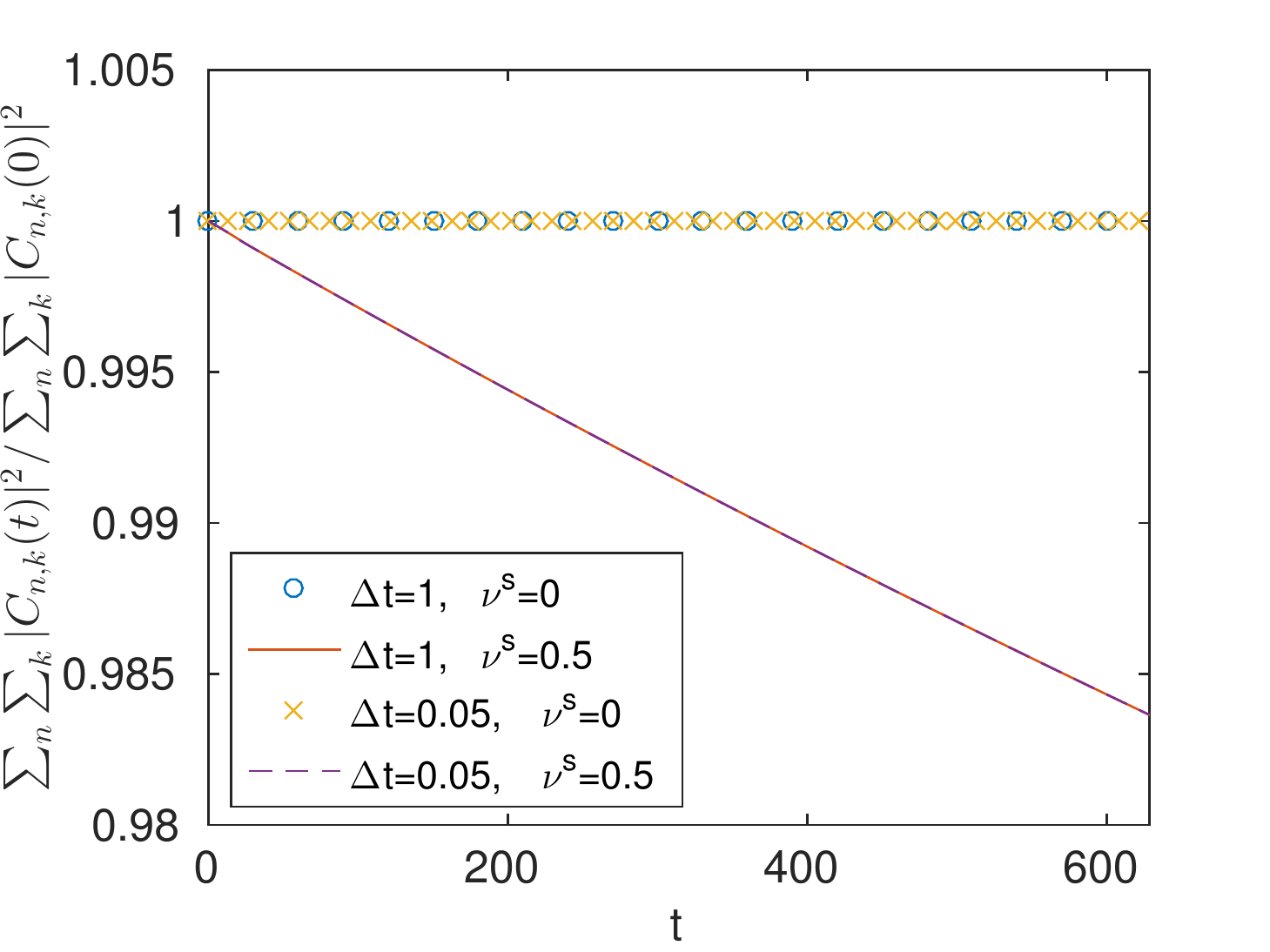}
      \put(48,-6){\FLARGE{\textbf{Time}}}
      \put(-8,24){\begin{sideways}\FLARGE{\textbf{$L^2$ norm of $\fs$}}\end{sideways}}
    \end{overpic}
    \vspace{0.5cm}
    \caption{Ion acoustic wave test: time variation of the $L2$ norm
      of the electron distribution function as predicted by
      Theorem~\ref{theo:L2stab} versus time for the four combinations
      of $\gamma^e\in\{0,\,0.5\}$ and $\nu^e\in\{0,1\}$.  Penalty
      $\gamma^s$ is not applied to the equations of the first three
      Legendre modes.}
  \label{fig:IA:01}
  \end{center}
\end{figure}

Figure \ref{fig:IA:01} shows the time evolution of the $L^2$ norm of
the distribution function normalized as
in~\eqref{eq:relative:L2-norm:fs} for the four simulations of
Fig.~\ref{fig:IA:00}.
As for the Landau damping case, when $\nu^s=0$ the $L^2$ norm of the
distribution function is flat (on the scale of the plot), indicating a
minimal contribution of the boundary terms in~\eqref{eq:L2stab}.
When $\nu^s=0.5$, the $L^2$ norm of the distribution function
decreases in time due to the dominant contribution of the collisional
term.

\begin{figure}
  \begin{center}
    \begin{overpic}[scale=0.8]{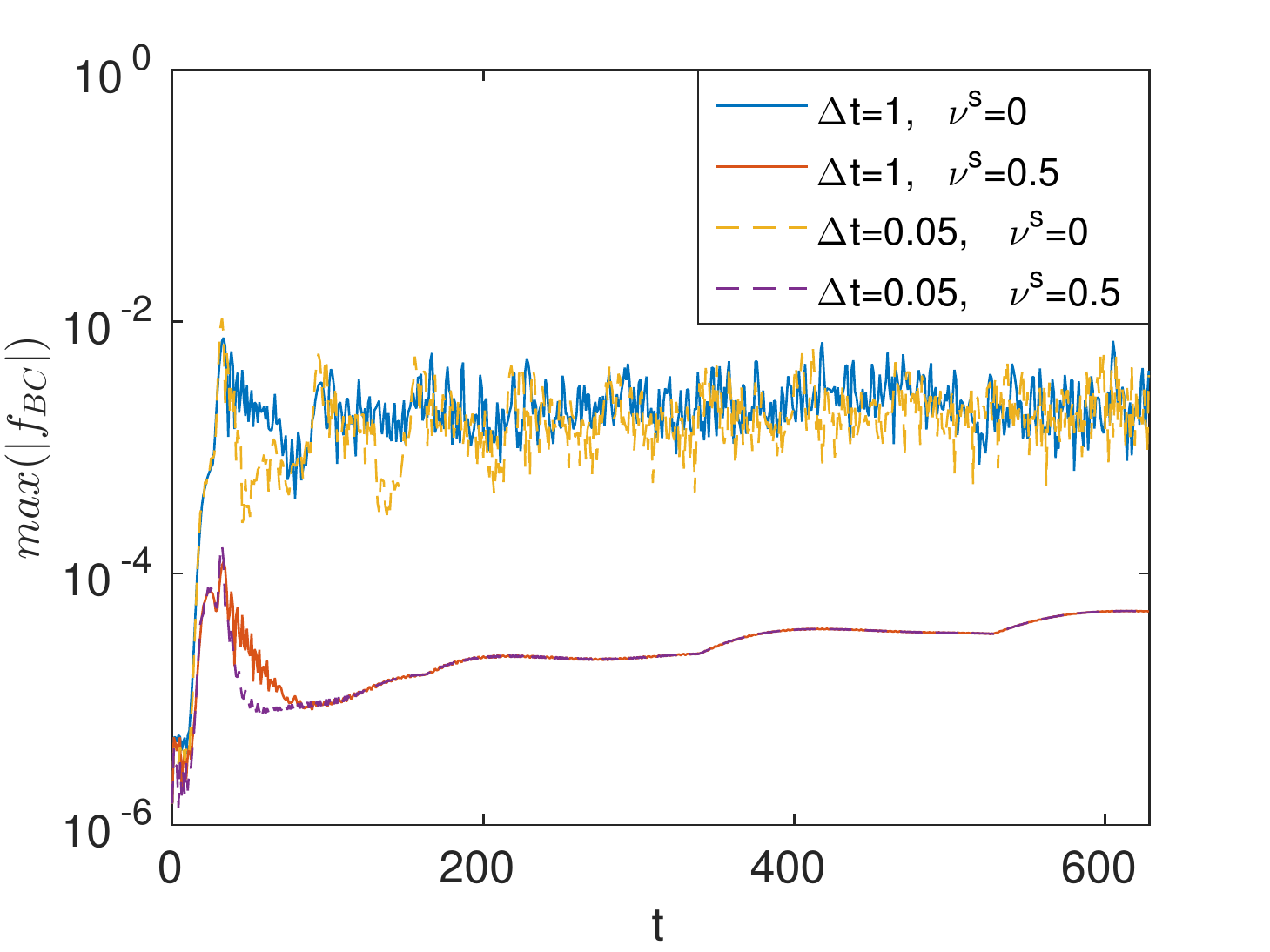}
      \put(48,-5){\FLARGE{\textbf{Time}}}
      \put(-5,29){\begin{sideways}
      \FLARGE{\textbf{max($\abs{\mathbf{f_{BC}}}$)}}
      \end{sideways}}
    \end{overpic}
    \vspace{0.5cm}
    \caption{Ion acoustic wave test: maximum value of the distribution
      function at the boundaries of the velocity range for
      $\gamma^e=0.5$ and the four combinations of $\Delta t=0.05,\,1$
      and $\nu^s=0,0.5$.  Penalty $\gamma^s$ is not applied to the
      equations of the first three Legendre modes.}
  \label{fig:IA:02}
  \end{center}
\end{figure}

Figure \ref{fig:IA:02} shows the maximum of $\abs{f}$ on the
boundaries of the velocity space, with the same format of
Fig. \ref{fig:IA:00}.
Although $\max\abs{f_{BC}}$ remains fairly small for all the cases,
once again one can see the beneficial effect of the collisional
operator: for $\nu^s=0.5$ it holds that $\max\abs{f_{BC}}$ is more
than an order of magnitude smaller than for $\nu^s=0$.

\begin{figure}
  \begin{center}
    \begin{overpic}[scale=0.5]{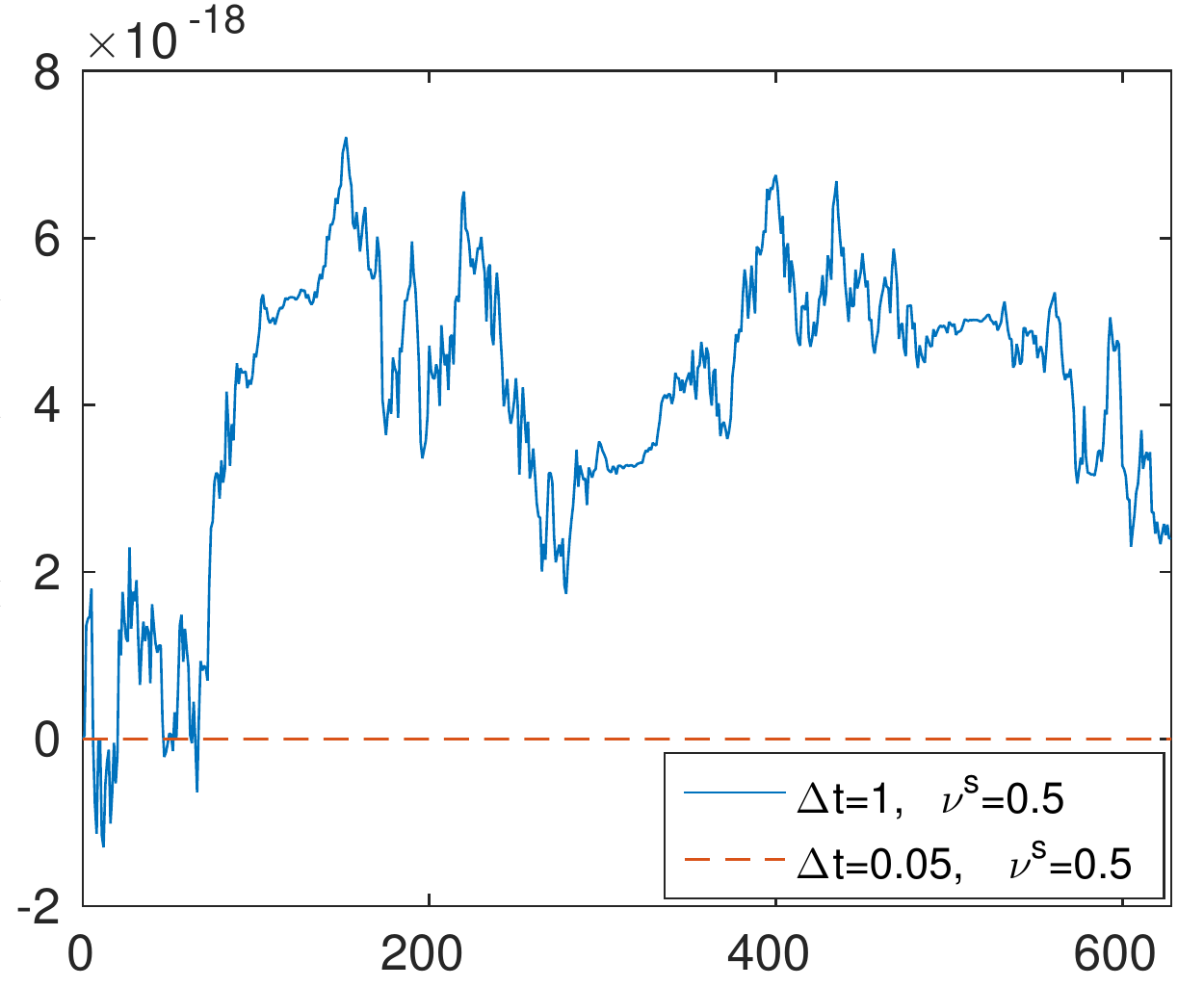}
      \put(48,-6){\textbf{Time}}
      \put(-8,0){\begin{sideways}\textbf{Variation of total momentum}\end{sideways}}
    \end{overpic}
    \hspace{0.75cm}
    \begin{overpic}[scale=0.5]{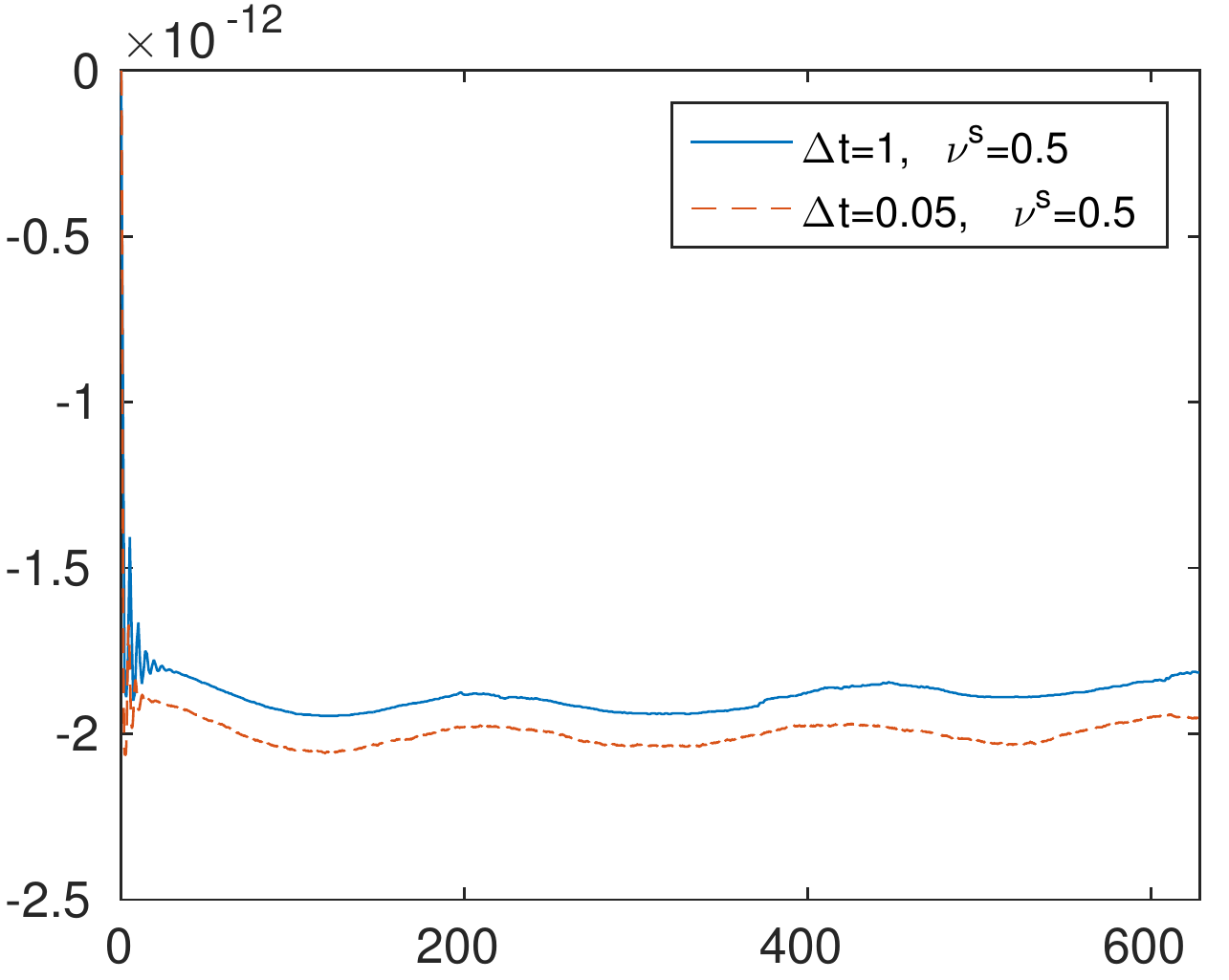}
      \put(48,-6){\textbf{Time}}
      \put(-8,6){\begin{sideways}\textbf{Variation of total energy}\end{sideways}}
    \end{overpic}
    \vspace{0.75cm}
    \caption{Ion acoustic wave test: the left panel shows the
      variation of momentum with respect to its initial value versus
      time, i.e., $P(t^{\tau})-P(0)$; the right panel shows the
      relative variation of energy with respect to its initial value
      versus time, i.e.,
      $(\mathcal{E}_{tot}(t^{\tau})-\mathcal{E}_{tot}(0))\slash{\mathcal{E}_{tot}(0)}$)
      for $\gamma^e=0.5$, $\Delta t\in\{0.05,\,1\}$ and $\nu^s=0.5$.
      Penalty $\gamma^s$ is not applied to the equations of the first
      three Legendre modes.}
  \label{fig:IA:03}
  \end{center}
\end{figure}

Finally, Fig. \ref{fig:IA:03} shows the time evolution of the total
momentum and the relative variation of the total energy for the
simulations with $\nu^s=0.5$ and $\Delta t=0.05$ and $\Delta t=1$
(total mass is not shown since it is conserved exactly).
In general, as expected, both quantities are conserved well.
One can notice that the error in the total momentum is controlled by
the time step, while this is not the case for the total energy.

\section{Conclusions}
\label{sec:conclusions}
In this paper a spectral method for the numerical solution of the
Vlasov-Poisson equations of a plasma has been presented.
The plasma distribution function is decomposed in Legendre polynomials
applied directly on a finite domain in velocity space.
The resulting set of moment equations is further discretized spatially
by a Fourier decomposition (periodic boundary conditions are assumed)
and in time by a fully-implicit, second order accurate Crank-Nicolson
scheme.
A collisional term is also considered in the discrete model to
  control the filamentation effect, but does not affect the
  conservation properties of the method.
A Jacobian-Free Newton-Krylov method (with the GMRES solver for the
inner linear iterations) is used to solve the discrete non-linear
equations.

The most significant aspects of our work are three.
First, the method is formulated in such a way that the boundary
conditions in velocity space ($f^s=0$ at the boundary of the velocity
domain) are applied in weak form.
That is, they are not enforced exactly through an expansion basis
obtained by a linear combination of the Legendre polynomials.
Instead, the boundary conditions are satisfied approximately via an
integration by parts once the Vlasov equation is projected onto the
Legendre basis functions (see Sec. \ref{sec:vlasov}).
Second, introducing a penalty on the weak form of the boundary
conditions allows the formulation of the numerical scheme to be
$L^2$-stable.
Third, the numerical scheme features conservation laws for total mass,
momentum and energy in weak form.
The numerical experiments performed in Sec. \ref{sec:numerical} on
Landau damping, two-stream instability and ion acoustic wave test
cases confirm both the stability of the method and the validity of the
conservation laws.



\clearpage

\clearpage

\clearpage

\clearpage
\ack
\vspace{-0.5\baselineskip}

The authors gratefully acknowledge discussions with C. La Cognata and
J. Nordstrom (University of Linkoping, Sweden); and L. Chac\'{o}n and
D. Moulton (Los Alamos National Laboratory); D. Funaro (University of
Modena and Reggio Emilia, Italy).

This work was partially funded by the Laboratory Directed Research and
Development program (LDRD), under the auspices of the National Nuclear
Security Administration of the U.S. Department of Energy by Los Alamos
National Laboratory, operated by Los Alamos National Security LLC
under contract DE-AC52-06NA25396.

The content of this article is also published in the journal paper 
of Reference~\cite{manzini2016legendre}.




\clearpage
\appendix
\section{Legendre polynomials: recursive relations}

\vspace{-0.25\baselineskip}
Consider the set of Legendre polynomials $\{L_n(s)\}_{n=0}^{\infty}$
that are recursively defined in $[-1,1]$ by~\eqref{eq:legendre:def}.
The two following recursion formulas hold:
\begin{subequations}
  \begin{align}
    \vs\phi_{n}(\vs)
    &= \sigma_{n+1}\,\phi_{n+1}( \vs )
    + \sigma_{n}  \,\phi_{n-1}( \vs )
    + \overline{\sigma}\,\phi_{n}  ( \vs ),
    \label{eq:recursion:formula:A}\\[0.5em]
    \vs^2\phi_{n}(\vs)
    &= 
    \sigma_{n+2}\sigma_{n+1}\,\phi_{n+2}( \vs )
    + 2\sigma_{n+1}\overline{\sigma}\,\phi_{n+1}( \vs )
    + \big( \sigma_{n+1}^2+\sigma_{n}^2+\overline{\sigma}^2 \big)\,\phi_{n}( \vs )
    \nonumber\\[0.5em]&\quad 
    + 2\sigma_{n}\overline{\sigma}\,\phi_{n-1}( \vs )
    + \sigma_{n}\sigma_{n-1}\,\phi_{n-2}( \vs ),
    \label{eq:recursion:formula:B}
  \end{align}
\end{subequations}
where $\sigma_{n}$ and $\overline{\sigma}$ are defined
in~\eqref{eq:legendre:sigma:def}.
To prove~\eqref{eq:recursion:formula:A}, note that the left-hand side
term and the two right-hand side terms of the recursion formula for
$n\geq 2$ can be rewritten as
\begin{align*}
  (2n+1) s L_{n}(s) 
  &= \frac{2}{\vb-\va} \left( \vs - \frac{\va+\vb}{2} \right)
  \sqrt{2n+1}\,\sqrt{2n+1}\,L_{n}( s(\vs) ) 
  = \frac{2\sqrt{2n+1}}{\vb-\va} \left( \vs - \frac{\va+\vb}{2} \right)
  \,\phi_{n}(\vs)
  \\[0.5em] 
  (n+1) L_{n+1}(s)
  &= \frac{n+1}{\sqrt{2(n+1)+1}}
  \sqrt{2(n+1)+1}\,L_{n+1}( s(\vs) )
  = \frac{n+1}{\sqrt{2(n+1)+1}}\phi_{n+1}(\vs)
  \\[0.5em]
  n L_{n-1}(s)
  &= \frac{n}{\sqrt{2(n-1)+1}}
  \sqrt{2(n-1)+1}\,L_{n-1}( s(\vs) )
  = \frac{n}{\sqrt{2(n-1)+1}}\,\phi_{n-1}( \vs )
\end{align*}
Collecting together and rearranging the three terms yields:
\begin{align*}
  \vs\phi_{n}(\vs) 
  = \frac{\vb-\va}{2}\left[
    \,\frac{n+1}{\sqrt{(2n+3)(2n+1)}}\,\phi_{n+1}( \vs )
    + \frac{n}  {\sqrt{(2n+1)(2n-1)}}\,\phi_{n-1}( \vs )
  \right]
  + \frac{\va+\vb}{2}\,\phi_{n}(\vs),
\end{align*}
which has the same form as~\eqref{eq:recursion:formula:A} where
$\sigma_{n}$ and $\overline{\sigma}$ can be readily determined by
comparison.
To prove~\eqref{eq:recursion:formula:B} just consider
$\vs^2\phi_n(\vs)=\vs\big( \vs\phi_n(\vs) \big)$ and
apply~\eqref{eq:recursion:formula:A} twice.
Moreover, a straightforward calculation yields
\begin{align}
  \sigma_{n}\sigma_{n,i} = \prty{n}{i}\frac{n\sqrt{2i+1}}{2n-1},
  \label{eq:sigma:sigma}
\end{align}
and in particular we have that
$\sigma_{2}\sigma_{2,1}=2\sigma_{1}\sigma_{1,0}=2$.

Integrating $\phi_n(\vs)$, $\vs\phi_n(\vs)$, $\vs^2\phi_n(\vs)$ and
using~\eqref{eq:recursion:formula:B}-\eqref{eq:recursion:formula:B}
give other three useful recurrence formulas:
\begin{subequations}
  \begin{align}
    \int_{\va}^{\vb}\phi_{n}(\vs)\d\vs 
    &= (\vb-\va)\delta_{n,0} 
    \label{eq:intg:f}\\[0.5em]
    \int_{\va}^{\vb}\vs\phi_{n}(\vs)\d\vs 
    &= (\vb-\va)\big( \sigma_{1}\delta_{n,1} + \overline{\sigma}\delta_{n,0} \big) 
    \label{eq:intg:vf}\\[0.5em]
    \int_{\va}^{\vb}\vs^2\phi_{n}(\vs)\d\vs 
    &= (\vb-\va)\Big(
    + \big( \sigma_{1}^2+\overline{\sigma}^2 \big)\,\delta_{n,0}
    + 2\sigma_{1}\overline{\sigma}\,\delta_{n,1}
    + \sigma_{2}\sigma_{1}\,\delta_{n,2}
    \Big).
    \label{eq:intg:v2f}
  \end{align}
\end{subequations}
All these three relations follows by noting that $\phi_{0}(\vs)=1$ and
applying the orthogonality property~\eqref{eq:Legendre:orthogonality}.
Relation~\eqref{eq:intg:f} is obvious.
To derive \eqref{eq:intg:vf} and \eqref{eq:intg:v2f} we also note that
we can remove the terms containing $\delta_{n+2,0}$ and
$\delta_{n+1,0}$ since $n\geq 0$.
Moreover, we can substitute $n=0,1,2$ in the $\sigma$-coefficients of
$\delta_{n,0}$, $\delta_{n-1,0}$, and $\delta_{n-2,0}$, and note that
the effect of $\delta_{n-1,0}$ and $\delta_{n-2,0}$ is respectively
equivalent to $\delta_{n,1}$ and $\delta_{n,2}$.
Finally, we note that $\sigma_{0}=0$.
Relation~\eqref{eq:intg:vf} follows from
\begin{align*}
  \int_{\va}^{\vb}\vs\phi_{n}(\vs)\d\vs
  &= \sigma_{n+1}\,\int_{\va}^{\vb}\phi_{n+1}( \vs )\d\vs
  + \sigma_{n}  \,\int_{\va}^{\vb}\phi_{n-1}( \vs )\d\vs
  + \overline{\sigma}\,\int_{\va}^{\vb}\phi_{n}  ( \vs )\d\vs,
  \\[0.5em]
  &= (\vb-\va)\big( \sigma_{n+1}\delta_{n+1,0} + \sigma_{n}\delta_{n-1,0} + \overline{\sigma}\delta_{n,0} \big).
\end{align*}
Relation~\eqref{eq:intg:v2f} follows from
\begin{align*}
  &\int_{\va}^{\vb}\vs^2\phi_{n}(\vs)\d\vs
  = 
  \sigma_{n+2}\sigma_{n+1}\,\int_{\va}^{\vb}\phi_{n+2}( \vs )\d\vs
  + 2\sigma_{n+1}\overline{\sigma}\,\int_{\va}^{\vb}\phi_{n+1}( \vs )\d\vs
  \hspace{4cm}
\end{align*}
\begin{align*}
  &+ \big( \sigma_{n+1}^2+\sigma_{n}^2+\overline{\sigma}^2 \big)\,\int_{\va}^{\vb}\phi_{n}( \vs )\d\vs
  + 2\sigma_{n}\overline{\sigma}\,\phi_{n+1}( \vs )\,\int_{\va}^{\vb}\phi_{n-1}( \vs )\d\vs
  + \sigma_{n}\sigma_{n-1}\,\int_{\va}^{\vb}\phi_{n-2}( \vs )\d\vs,
  \\[0.5em]
  &= 
  (\vb-\va)\big( \sigma_{n+2}\sigma_{n+1}\,\delta_{n+2,0}
  + 2\sigma_{n+1}\overline{\sigma}\,\delta_{n+1,0}
  + \big( \sigma_{n+1}^2+\sigma_{n}^2+\overline{\sigma}^2 \big)\,\delta_{n,0}
  \\[0.5em]
  &\qquad 
  + 2\sigma_{n}\overline{\sigma}\,\delta_{n+1,0}+\,\delta_{n-1,0}
  + \sigma_{n}\sigma_{n-1}\,\delta_{n-2,0} \big).
\end{align*}


\section{Proof of~\eqref{eq:L2:stability:fs}.}
\label{sec:proof:L2:stability:fs}

\vspace{-0.25\baselineskip}
The proof of equation~\eqref{eq:L2:stability:fs} starts by applying
expansion~\eqref{eq:legendre:decomposition}, Legendre orthogonality
property~\eqref{eq:Legendre:orthogonality},
expansion~\eqref{eq:fourier:decomposition} and Fourier orthogonality
property~\eqref{eq:Fourier:orthogonality}:
\begin{align*}
  &\int_{0}^{L}\int_{\va}^{\vb}\abs{\fs(x,\vs,t)}^2\d\vs\dx
  = \sum_{m,n=0}^{\NL-1}\int_{0}^{L}\Cs_{m}(x,t)\Cs_{n}(x,t)\int_{\va}^{\vb}\phi_m(\vs)\phi_n(\vs)\d\vs\dx
  \\[0.5em]
  &\qquad= (\vb-\va)\sum_{m,n=0}^{\NL-1}\int_{0}^{L}\Cs_{m}(x,t)\Csn(x,t)\delta_{m,n}\dx
  = (\vb-\va)\sum_{n=0}^{\NL-1}\int_{0}^{L}\abs{\Csn(x,t)}^2\dx
  \\[0.5em]
  &\qquad= (\vb-\va)\sum_{n=0}^{\NL-1}\sum_{k,k'=-\NF}^{\NF}\big(\Cs_{n,k}(t)\big)^{\dagger}\Cs_{n,k'}(t)\int_{0}^{L}\psi_{-k}(x)\psi_{k}(x)\dx
  \\[0.5em]
  &\qquad= (\vb-\va)L\,\sum_{n=0}^{\NL-1}\sum_{k,k'=-\NF}^{\NF}\big(\Cs_{n,k}(t)\big)^{\dagger}\Cs_{n,k'}(t)\delta_{-k+k',0}
  = (\vb-\va)L\,\sum_{n=0}^{\NL-1}\sum_{k=-\NF}^{\NF}\abs{\Cs_{n,k}(t)}^2.
\end{align*}

\section{Proof of Lemma~\ref{lemma:L2stab:useful}.}

\vspace{-0.25\baselineskip}
To prove the left-most equality
in~\eqref{eq:legendre:fourier:L2stab:useful}, we first note that:
\begin{align}
  L\sum_{k=-\Mk}^{\Mk}(\Cvk)^{\dagger}\Big[\Es\STAR\matB\Cvs\Big]_k
  =L\sum_{n=0}^{\NL-1}\sum_{k=-\Mk}^{\Mk}(\Cs_{n,k})^{\dagger}\Big[\Es\STAR(\matB\Cvs)_{n}\Big]_k
  \label{eq:lemma:proof:10:a}
\end{align}
Using the definition of the discrete Fourier expansion of the
electric field $\Es$, the Legendre coefficients $\Csn(x,t)$, and
$(\matB\Cvs)_{n}$, we obtain:
\begin{align}
  &L\sum_{k=-\Mk}^{\Mk}(\Cs_{n,k})^{\dagger}\Big[\Es\STAR(\matB\Cvs)_{n}\Big]_k
  =L\sum_{k,k'=-\Mk}^{\Mk}(\Cs_{n,k})^{\dagger}\Es_{k'}\big(\matB\Cvs\big)_{n,k-k'}\nonumber\\[0.5em]
  &\qquad\qquad= \sum_{k,k',k''=-\Mk}^{\Mk}(\Cs_{n,k})^{\dagger}\Es_{k'}\big(\matB\Cvs\big)_{n,k''}\,L\delta_{-k+k'+k'',0}\nonumber\\[0.5em]
  &\qquad\qquad= \int_{0}^{L}
  \left(\sum_{k=-\Mk}^{\Mk} (\Cs_{n,k})^{\dagger}\psi_{-k}(x)\right)
  \left(\sum_{k'=-\Mk}^{\Mk}\Es_{k'}\psi_{k'}(x)\right)
  \left(\sum_{k''=-\Mk}^{\Mk}\big(\matB\Cvs\big)_{n,k''}\psi_{k''}(x)\right)
  \,\dx\nonumber\\[0.5em]
  &\qquad\qquad= \int_{0}^{L}\Cs_{n}(x,t)\Es(x,t)\big(\matB\Cvs\big)_{n}.
  \label{eq:lemma:proof:10}
\end{align}
Then, we note that:
\begin{align*}
  \begin{array}{l}
    \displaystyle\sum_{n=0}^{\NL-1}\Csn(x,t)\,(\matB\Cvs(x,t))_{n}=
    \qquad\qquad\qquad\qquad\qquad\qquad\qquad\,\,\,\,\,\,\,
    \mbox{\big[use equation~\eqref{eq:legendre:matB:def}\big]}
    \qquad\qquad\qquad\qquad
    \\[1.75em]
  \end{array}
\end{align*}
\begin{align*}
  \begin{array}{rll}
    &\qquad\displaystyle=\sum_{n=0}^{\NL-1}\sum_{i=0}^{n-1}\sigma_{n,i}\Csn(x,t)\Cs_{i}(x,t)
    &\qquad\mbox{\big[use definition of $\Cs_{i}$ from~\eqref{eq:legendre:decomposition}\big]}\\[1.5em]
    &\qquad\displaystyle=\frac{1}{\vb-\va}\sum_{n=0}^{\NL-1}\sum_{i=0}^{n-1}\sigma_{n,i}\Csn(x,t)\int_{\va}^{\vb}\fs(x,\vs,t)\phi_{i}(\vs)\d\vs
    &\qquad\mbox{\big[use derivative formula~\eqref{eq:legendre:first:derivative}\big]}\\[1.5em]
    &\qquad\displaystyle=\frac{1}{\vb-\va}\sum_{n=0}^{\NL-1}\Csn(x,t)\int_{\va}^{\vb}\fs(x,\vs,t)\frac{d\phi_{n}(\vs)}{d\vs}\d\vs
    &\qquad\mbox{\big[use again decomposition~\eqref{eq:legendre:decomposition}\big]}\\[1.5em]
    &\qquad\displaystyle=\frac{1}{\vb-\va}\int_{\va}^{\vb}\fs(x,\vs,t)\frac{\partial\fs(x,\vs,t)}{\partial\vs}\d\vs
    &\qquad\mbox{\big[use the definition of the derivative\big]}\\[1.5em]
    &\qquad\displaystyle
    =\frac{1}{2(\vb-\va)}\int_{\va}^{\vb}\frac{\partial(\fs)^2}{\partial\vs}\d\vs
    = \displaystyle\frac{1}{2}\Diffv{(\fs)^2}_{\va}^{\vb}
  \end{array}
\end{align*}
Using the last relation above in~\eqref{eq:lemma:proof:10} and
tranforming back in Fourier space yield:
\begin{align}
  L\sum_{k=-\Mk}^{\Mk}(\Cvk)^{\dagger}\Big[\Es\STAR\matB\Cvs\Big]_k
  &=\sum_{n=0}^{\NL-1}\int_{0}^{L}\Cs_{n}(x,t)\Es(x,t)\big(\matB\Cvs\big)_{n}
  =\frac{1}{2}\int_{0}^{L}\Es(x,t)\Diffv{(\fs(x,\vs,t))^2}_{\va}^{\vb}\,\dx
  \nonumber\\[0.5em]
  &=\frac{1}{2}L\left[\Es\STAR\Diffv{(\fs(x,\vs,t))^2}_{\va}^{\vb}\right]_0,
  \label{eq:lemma:proof:15}
\end{align}
where $[\,\ldots\,]_{0}$ denotes the zero-th Fourier mode and which is the
first equality in~\eqref{eq:legendre:fourier:L2stab:useful}.

\medskip 
Applying again the definition of the discrete Fourier transform, the
right-most equality in~\eqref{eq:legendre:fourier:L2stab:useful} is
proved as follows:
\begin{align}
  &L\sum_{k=-\Mk}^{\Mk}(\Cvk)^{\dagger}\Big[\Es\STAR\delta_{\vs}\big[\fs\bm\phi\big]_{\va}^{\vb}\Big]_k
  =L\sum_{k=-\Mk}^{\Mk}\sum_{n=0}^{\NL-1}(\Cs_{n,k})^{\dagger}\Big[\Es\STAR\delta_{\vs}\big[\fs\phi_n\big]_{\va}^{\vb}\Big]_k
  \nonumber\\[0.5em]
  &\qquad
  =L\sum_{n=0}^{\NL-1}\sum_{k,k'=-\Mk}^{\Mk}(\Cs_{n,k})^{\dagger}\Es_{k'}\Big(\delta_{\vs}\big[\fs\phi_n\big]_{\va}^{\vb}\Big)_{k-k'}
  \nonumber\\[0.5em]
  &\qquad=\sum_{n=0}^{\NL-1}\sum_{k,k',k''=-\Mk}^{\Mk}(\Cs_{n,k})^{\dagger}\Es_{k'}\Big(\delta_{\vs}\big[\fs\phi_n\big]_{\va}^{\vb}\Big)_{k''}
  \,L\delta_{-k+k'+k'',0}\nonumber\\[0.5em]
  &\qquad=\sum_{n=0}^{\NL-1}\int_{0}^{L} 
  \left(\sum_{k=-\Mk}^{\Mk}(\Cs_{n,k})^{\dagger}\psi_{-k}(x)\right)
  \left(\sum_{k'=-\Mk}^{\Mk}\Es_{k'}\psi_{k'}(x)\right)
  \left(\sum_{k''=-\Mk}^{\Mk}\Big(\delta_{\vs}\big[\fs\phi_n\big]_{\va}^{\vb}\Big)_{k''}\psi_{k''}(x)\right)
  \,\dx
  \nonumber\\[0.5em]
  &\qquad
  =\sum_{n=0}^{\NL-1}\int_{0}^{L}\Cs_{n}(x,t)\Es(x,t)\delta_{\vs}\big[\fs(x,\vs,t)\phi_n(\vs)\big]_{\va}^{\vb}\dx
  \nonumber\\[0.5em]
  &\qquad
  =\int_{0}^{L}\Es(x,t)\delta_{\vs}\big[\fs(x,\vs,t)\,\sum_{n=0}^{\NL-1}\Cs_{n}(x,t)\phi_n(\vs)\big]_{\va}^{\vb}\dx
  \nonumber\\[0.5em]
  &\qquad
  =\int_{0}^{L}\Es(x,t)\delta_{\vs}\big[\fs(x,\vs,t)^2\big]_{\va}^{\vb}\dx
  =L\left[\Es\STAR\delta_{\vs}\big[(\fs)^2\big]_{\va}^{\vb}\right]_{0}.
  \label{eq:lemma:proof:20}
\end{align}
The three members of~\eqref{eq:legendre:fourier:L2stab:useful} are
real numbers since intermediate steps in the previous developments are
formed by real quantities.


\end{document}